\def\VERSION{7.12.2025}
\def\WHO{nbd} 
\def\users{us}    
\def\users{world} 
\documentclass[11pt]{article} 

\usepackage{color}
\usepackage{amsmath}
\usepackage{amsthm}
\usepackage{amssymb}
\usepackage{epsfig}
\usepackage{psfrag}
\usepackage{graphicx}
\usepackage{textcomp}
\usepackage{bm}
\usepackage{pgf}
\numberwithin{equation}{section}
\usepackage{upgreek} 
\newtheorem{theorem}{Theorem}[section]

\newtheorem{definition}[theorem]{Definition}
\newtheorem{example}[theorem]{Example}
\newtheorem{proposition}[theorem]{Proposition}

\newtheorem{remark}[theorem]{Remark}

\usepackage{mathrsfs,cite}
\usepackage{ifthen}
\usepackage{ulem}
\usepackage[colorlinks=true, pdfstartview=FitV, linkcolor=blue, 
            citecolor=blue, urlcolor=blue]{hyperref}
\usepackage{cancel}\ifthenelse{\equal{\users}{world}}
{
\newcommand{\REM}[1]{}

	\newcommand{\DELETE}[1]{}

        \newcommand{\COMMENT}[1]{}
        \newcommand{\TCOMMENT}[1]{}
    \newcommand{\MARGINOTE}[1]{}
}	
{
\usepackage[notcite,notref,color]{showkeys}

\definecolor{brown}{rgb}{0.6,0.2,0.2}
\newcommand{\REM}[1]{\marginpar{\bfseries\tiny{\color{blue}#1}}}

 \newcommand{\COMMENT}[1]{{\color{blue}\uuline{#1}\color{black}}} 
 \newcommand{\DELETE}[1]{{\color{brown}\cancel{#1}\color{black}}}

 \newcommand{\TCOMMENT}[1]{{\color{blue}{ #1}}}
\newcommand{\MARGINOTE}[1]{\marginpar{\color{red}\tiny\texttt{#1}}}
\usepackage{fancyhdr}
\pagestyle{fancy}
\headheight=28pt\headwidth=17cm
\definecolor{gray}{gray}{0.5}
\newcount\hour \newcount\minute
\hour=\time
\divide \hour by 60
\minute=\time
\loop \ifnum \minute > 59 \advance \minute by -60 \repeat

\rhead{\color{gray}time discretization Eulerian\\
by T.Roub\'\i\v cek}
\chead{}
\lhead{Version \VERSION, file: \jobname.tex, \underline{\Large\color{red}\WHO's working on} \\
compiled:
\number\day.\number\month.\number\year\ at
\the\hour:\ifnum\minute<10 0\fi\the\minute\ h }
}

\newcommand{\R}{\mathbb{R}}

\newcommand{\bbI}{\mathbb{I}}

\newcommand{\bbD}{\mathbb{D}}

\newcommand\DT[1]{\mathchoice
                 {{\buildrel{\hspace*{.1em}\text{\LARGE.}}\over{#1}}}
                 {{\buildrel{\hspace*{.1em}\text{\LARGE.}}\over{#1}}}
                 {{\buildrel{\hspace*{.1em}\text{\Large.}}\over{#1}}}
                 {{\buildrel{\hspace*{.1em}\text{\large.}}\over{#1}}}}

\newcommand{\lineunder}[2]{\LU{\begin{array}[t]{c}\underbrace{#1}\vspace*{.5em}\end{array}}{\mbox{\footnotesize\rm #2}}}
\newcommand{\linesunder}[3]{\LSU{\begin{array}[t]{c}\underbrace{#1}\vspace*{.5em}\end{array}}{\mbox{\footnotesize\rm #2}}{\mbox{\footnotesize\rm#3}}}
\newcommand{\LU}[2]{\begin{array}[t]{c}#1\vspace*{-1em}\\_{#2}\end{array}}
\newcommand{\LSU}[3]{\begin{array}[t]{c}#1\vspace*{-1em}\\_{#2}\vspace*{-.5em}\\_{#3}\end{array}}
\newcommand{\threelinesunder}[4]{\threeLSU{\begin{array}[t]{c}\underbrace{#1}\vspace*{.5em}\end{array}}{\mbox{\footnotesize\rm #2}}{\mbox{\footnotesize\rm #3}}{\mbox{\footnotesize\rm #4}}}
\newcommand{\threeLSU}[4]{\begin{array}[t]{c}#1\vspace*{-1em}\\_{#2}\vspace*{-.5em}\\_{#3}\vspace*{-.5em}\\_{#4}\end{array}}
\renewcommand{\d}{{\rm d}}

\newcommand{\divS}{\mathrm{div}_{\scriptscriptstyle\textrm{\hspace*{-.1em}S}}^{}}
\newcommand{\eq}[1]{(\ref{#1})}
\newcommand{\Cdot}{\hspace{-.1em}\cdot\hspace{-.1em}}
\newcommand{\Colon}{\hspace{-.15em}:\hspace{-.15em}}
\newcommand{\UUU}[3]{\begin{array}[b]{c}\vspace*{-1.2mm}_{\text{\scriptsize{#2}}}\vspace*{-.6mm}\\[-.0em]_{\text{\scriptsize{#3}}}\vspace*{.6mm}\\[-.5em]#1\end{array}}

\def\vv{{\bm v}}
\def\pp{{\bm p}}
\def\uu{{\bm u}}

\def\xx{{\bm x}}
\def\yy{{\bm y}}
\def\nn{{\bm n}}

\newcommand{\DD}{\bm D}
\newcommand{\TT}{\bm T}

\newcommand{\strain}{{\boldsymbol\varepsilon}}
\newcommand{\GRAVITY}{\bm g}

\newcommand\EE{{\bm e}}

\def\vvk{\vv_\etau^k}
\def\vvkk{\vv_\etau^{k-1}}
\def\overlineEetau{\hspace*{.2em}\overline{\hspace*{-.2em}\bm E}_\etau^{}}
\def\overlineFetau{\hspace*{.2em}\overline{\hspace*{-.2em}\bm F}_{\!\etau}^{}}

\def\overlinevvtau{\hspace*{.15em}\overline{\hspace*{-.15em}\vv}_{\etau}^{}}

\def\WW{{\bm W}}
\def\Lp{{\bm \varPi}}

\def\HYPER{\mu}

\def\jj{{\bm j}}
\def\ZJEp{{\bm\varPi}}
\def\Eng{e}
\def\ENG{\mbox{\small{${\mathscr E}$}}}
\def\Ent{u}
\def\ENT{\mbox{\small{${\mathscr U}$}}}
\def\intkappa{\varkappa}

\newcommand\ZJ[1]{\mathchoice
                 {{\buildrel{\hspace*{.1em}{_{\,\boldsymbol\circ}}}\over{#1}}}
                 {{\buildrel{\hspace*{.1em}{_{\,\boldsymbol\circ}}}\over{#1}}}
                 {{\buildrel{\hspace*{.1em}{\boldsymbol\circ}}\over{#1}}}
                 {{\buildrel{\hspace*{.1em}{\boldsymbol\circ}}\over{#1}}}}
\newcommand{\barOmega}{\,\overline{\!\varOmega}}

\newcommand{\Nabla}{\nabla}
\newcommand{\Rsym}{\mathbb R_{\rm sym}^{3\times 3}}
\def\Rdev{\R_{\rm dev}^{3\times 3}}

\newcommand{\Item}[2]{\parbox[t]{.05\textwidth}{#1}\hfill%
      \parbox[t]{.95\textwidth}{#2}\vspace*{.8mm}} 

\def\Vdots{\!\mbox{\setlength{\unitlength}{1em}
\begin{picture}(0,0)
\put(-.07,0){.}
\put(-.07,.3){.}
\put(-.07,.6){.}
\end{picture}
}
}

\newcommand{\wt}[1]{\mathchoice{\hspace*{-.09em}\text{\large$\hspace*{.09em}\tilde{\text{\normalsize$#1$}}\hspace*{.05em}$}\hspace*{-.05em}}
{\hspace*{-.09em}\text{\large$\hspace*{.09em}\tilde{\text{\normalsize$#1$}}\hspace*{.05em}$}\hspace*{-.05em}}
{\text{\normalsize$\hspace*{.08em}\tilde{\text{\scriptsize$#1$}}\hspace*{.06em}$}}
{\text{\small$\tilde{\text{\tiny$#1$}}$}}}
\newcommand{\wb}[1]{\mathchoice{\hspace*{-.09em}\text{\large$\hspace*{.09em}\bar{\text{\normalsize$#1$}}\hspace*{.05em}$}\hspace*{-.05em}}
  {\hspace*{-.09em}\text{\large$\hspace*{.09em}\bar{\text{\normalsize$#1$}}\hspace*{.05em}$}\hspace*{-.05em}}
  {\text{\normalsize$\hspace*{.08em}\bar{\text{\scriptsize$#1$}}\hspace*{.06em}$}}
  {\text{\small$\bar{\text{\tiny$#1$}}$}}}

\newcounter{myfigure}
\newenvironment{my-picture}[3]{\refstepcounter{myfigure}\label{#3}\setlength{\unitlength}{1em}\begin{picture}(#1,#2)}{\end{picture}}

\newcommand\DELETEDELETE[1]{}

\newcommand\pdt[1]{\frac{\partial{#1}}{\partial t}}
\newcommand\Ee{{\bm E}}              
\newcommand\Ep{{\bm P}}              

\begin{document}
\begin{sloppypar}

\allowdisplaybreaks

\def\EPS{\varepsilon}
\def\DELTA{\delta}
\def\etau{{\EPS\DELTA\tau}}
\def\EEps{{\EPS\DELTA}}

\def\TTtauk{\TT^k_{\!\etau}}
 \def\Eek{\Ee_\etau^k}
 \def\Eetau{\Ee_\etau^{}}
\def\overlineDetau{\hspace*{.2em}\overline{\hspace*{-.2em}\bm D}_\etau^{}}
\def\overlineTetau{\hspace*{.2em}\overline{\hspace*{-.1em}\bm T}_\etau^{}}
\def\overlineLpetau{\hspace*{.2em}\overline{\hspace*{-.3em}\bm\varPi\hspace*{-.1em}}_\etau^{}}
\def\overlineMetau{\hspace*{.4em}\overline{\hspace*{-.4em}\bm M\hspace*{-.1em}}_\etau^{}}
\def\overlineFetau{\hspace*{.2em}\overline{\hspace*{-.2em}\bm F}_{\!\etau}^{}}
\def\overlineHetau{\hspace*{.2em}\overline{\hspace*{-.2em}\bm H}_\etau^{}}
\def\overlineEetau{\hspace*{.2em}\overline{\hspace*{-.2em}\bm E}_\etau^{}}
\def\overlinevvtau{\hspace*{.1em}\overline{\hspace*{-.1em}\vv}_{\etau}^{}}
\def\overlineppetau{\hspace*{.15em}\overline{\hspace*{-.15em}\pp}_\etau^{}}
\def\overlinerhoetau{\hspace*{.15em}\overline{\hspace*{-.15em}\varrho}_\etau^{}}

\def\Ke{K_\text{\sc e}^{}}
\def\Ge{G_\text{\sc e}^{}}
\def\Kv{K_\text{\sc v}^{}}
\def\Gv{G_\text{\sc v}^{}}
\def\Gm{G_\text{\sc m}^{}}

\def\oDetau{\hspace*{.2em}\overline{\hspace*{-.2em}\bm D}_\etau^{}}
\def\overlineTetau{\hspace*{.2em}\overline{\hspace*{-.1em}\bm T}_\etau^{}}
\def\overlineLpetau{\hspace*{.2em}\overline{\hspace*{-.3em}\bm\varPi\hspace*{-.1em}}_\etau^{}}
\def\oMetau{\hspace*{.4em}\overline{\hspace*{-.4em}\bm M\hspace*{-.1em}}_\etau^{}}
\def\oFetau{\hspace*{.2em}\overline{\hspace*{-.2em}\bm F}_{\!\etau}^{}}
\def\oHetau{\hspace*{.2em}\overline{\hspace*{-.2em}\bm H}_\etau^{}}
\def\oEetau{\hspace*{.2em}\overline{\hspace*{-.2em}\bm E}_\etau^{}}
\def\othetatau{\hspace*{.2em}\overline{\hspace*{-.2em}\theta}_\etau^{}}
\def\oEpetau{\hspace*{.2em}\overline{\hspace*{-.2em}\bm\Pi}_\etau^{}}
\def\oeetau{\hspace*{.2em}\overline{\hspace*{-.05em}e}_\etau^{}}
\def\ocetau{\hspace*{.2em}\overline{\hspace*{-.05em}c}_\etau^{}}
\def\omuetau{\hspace*{.2em}\overline{\hspace*{-.05em}\mu}_\etau^{}}
\def\ouetau{\hspace*{.2em}\overline{\hspace*{-.05em}\Ent}_\etau^{}}
\def\overlinezetau{\hspace*{.05em}\overline{\hspace*{-.05em}\bm z}_\etau^{}}
\def\ovetau{\hspace*{.1em}\overline{\hspace*{-.1em}\vv}_{\etau}^{}}
\def\overlinethetatau{\hspace*{.1em}\overline{\hspace*{-.1em}\theta}_{\etau}^{}}
\def\underlinethetatau{\hspace*{.1em}\underline{\hspace*{-.1em}\theta}_{\etau}^{}}
\def\overlinethetatauTWO{\hspace*{.1em}\overline{\hspace*{-.1em}\theta}_{\etau}^{\,2}}
\def\overlinethetatauEXP{\hspace*{.1em}\overline{\hspace*{-.1em}\theta}_{\etau}^{\,\EXP}}
\def\overlinethetatauExp{\hspace*{.1em}\overline{\hspace*{-.1em}\theta}_{\etau}^{\,{\lambda}}}
\def\overlinethetatauEXp{\hspace*{.1em}\overline{\hspace*{-.1em}\theta}_{\etau}^{\,_{1+\lambda}}}
\def\opetau{\hspace*{.15em}\overline{\hspace*{-.15em}\pp}_\etau^{}}
\def\oretau{\hspace*{.15em}\overline{\hspace*{-.15em}\varrho}_\etau^{}}
\def\overlineZJEpetau{\hspace*{.15em}\overline{\hspace*{-.15em}\ZJEp}_\etau^{}}
\def\vvkk{\vv_\etau^{k-1}}
\def\Eetau{\Ee_\etau^{}}
\def\eetau{{\rm e},\etau} 
\def\petau{{\rm p},\etau} 
\def\oLpetau{\hspace*{.2em}\overline{\hspace*{-.2em}\bm L}_{\petau}^{}}
\def\oJetau{\hspace*{.2em}\overline{\hspace*{-.15em}J}_{\etau}^{}}
\def\oFeetau{\hspace*{.2em}\overline{\hspace*{-.2em}\bm F}_{\!\eetau}^{}}
\def\oHetau{\hspace*{.2em}\overline{\hspace*{-.2em}\mathfrak{H}}_{\etau}^{}}
\def\oDetau{\hspace*{.2em}\overline{\hspace*{-.2em}\bm D}_\etau^{}}

\def\ovT{\wb{\vv}_\TAU}
\def\opT{\wb{\pp}_\TAU}
\def\orT{\wb{\varrho}_\TAU}
\def\oEeT{\hspace*{.2em}\overline{\hspace*{-.2em}\bm E}_{\!\TAU}^{}}
\def\uET{\underline{\EE}_\TAU}
\def\oPT{\hspace*{.3em}\overline{\hspace*{-.3em}\ZJEp\hspace*{-.1em}}_{\TAU}^{}}
\def\oTT{\hspace*{.15em}\overline{\hspace*{-.15em}\TT}_{\TAU}^{}}
\def\oDT{\wb{\DD}_{\!\TAU}}
\def\otT{\wb{\theta}_\TAU}
\def\utT{\underline{\theta}_\TAU}
\def\ouT{\wb{\Ent}_\TAU}
\def\vvk{\vv_\TAU^k}
\def\vvkk{\vv_\TAU^{k-1}}

\def\ALPHA{\alpha}

\noindent{\LARGE\bf Time discretization in convected linearized
\\[.4em]            thermo-visco-elastodynamics at large displacements.}

\bigskip\bigskip

\noindent{\large\sc Tom\'{a}\v{s} Roub\'\i\v{c}ek}\\
{\it Mathematical Institute, Charles University, \\Sokolovsk\'a 83,
CZ--186~75~Praha~8,  Czech Republic
}\\and\\
{\it Institute of Thermomechanics, Czech Academy of Sciences,\\Dolej\v skova~5,
CZ--182~08 Praha 8, Czech Republic
\\{\tt ORCID}: {\rm 0000-0002-0651-5959}, Email: {\tt tomas.roubicek@mff.cuni.cz}}

\bigskip\medskip

\begin{center}\begin{minipage}[t]{14.6cm}

{\small

\noindent{\bfseries Abstract.}
The fully-implicit time discretization (i.e.\ the 
backward Euler formula) is applied to compressible nonlinear
dynamical models of thermo-viscoelastic solids in the
Eulerian description, i.e.\ in
the actual deforming configuration, formulated in terms of rates.
The Kelvin-Voigt rheology or also, in the deviatoric part, the
Jeffreys rheology (which covers creep or plasticity) are considered,
using the additive Green-Naghdi
decomposition of total strain into the elastic and the inelastic strains
formulated in terms of (objective) rates exploiting the Zaremba-Jaumann
time derivative. A linearized convective model at large displacements
is considered, focusing on the case where the internal energy additively
splits the (convex) mechanical and the thermal parts.  A fully
implicit  time-discrete scheme is devised. Considering the
multipolar 2nd-grade viscosity, the numerical stability and convergence
towards weak solutions are proven by exploiting, in particular, the
convexity of the kinetic energy when written in terms of linear momentum
instead of velocity and by estimating the temperature gradient from
the entropy-like inequality. 

\medskip


\medskip


} 
\end{minipage}
\end{center}

\smallskip

\baselineskip=14pt

\section{Introduction}

The models of {\it thermo-visco-elastodynamics} in continuum mechanics
{\it at finite} (also called {\it large}) {\it strains} lead to
strongly  nonlinear systems of evolution partial differential
equations. The basic modelling paradigms are the Lagrangian description
(using a referential configuration) versus the  {\it Eulerian description}
(using the actual deforming configuration). In this paper, we focus on the
latter. The Eulerian approach is particularly justified in situations where 
no referential configuration exists for physical reasons, such as in
fluids or in solids undergoing creep as in geophysical models
on long timescales of millions of years). The description in the actual
configuration reveals the ``actual'' physics more explicitly than
Lagrangian approach.

To avoid truly large-strain Eulerian models, a certain modeling compromise
is to use a linearization that leads to a symmetric small-strain tensor while
allowing {\it large displacements} and incorporating suitable convected
objective time derivatives. Such models are widely used in engineering
and geophysics.

Such models can be approximated numerically using the time discretization
({\it Rothe}'s) {\it method}. After further space discretization, this
method yields computationally implementable schemes.
An analytically rigorous treatment of time discretization has been
launched for isothermal compressible fluids in 
\cite{GaMaNo19EEIM,FHMN17EENM,FeKaPo16MTCV,FLMS19CFVS,Kar13CFEM,Zato12ASCN}.
and for isothermal compressible solids in \cite{Roub25TDVE}.
This paper aims to extend this treatment to extend this treatment to
the anisothermal compressible solids.

The main attributes of this paper are as follows:
\begin{itemize}
\vspace*{-.8em}\item[$\rhd$]the (generalized) entropy inequality is exploited
for a-priori estimation of the temperature gradient,
\vspace*{-.8em}\item[$\rhd$]the convexity of the kinetic energy of the
linear momentum is used for the time discretization as in \cite{Roub25TDVE},
 as well as the convexity of $\varrho\mapsto\sigma:=1/\varrho$ for
estimation of the so-called sparsity $\sigma$, 
\vspace*{-.8em}\item[$\rhd$]
the discrete Gronwall inequality is used for a sufficiently small time
step without regularing the discrete kinematics when estimating the strain
gradient, which also allows the strong convexity assumption in
\cite{Roub25TDVE} to be weakened to mere convexity, and
\vspace*{-.8em}\item[$\rhd$]
the time discretization is applied directly to the original, non-regularized
problem.
\end{itemize}

\vspace*{-.7em}

\noindent
In comparison with the aforementioned paper \cite{Roub25TDVE}, we account
for Maxwellian {\it creep}  in the deviatoric part, which, in combination
with the Kelvin-Voigt rheology, yields the {\it Jeffreys} (also called the
anti-Zener) {\it viscoelastic model}. This  
creep may be nonlinear and may include plasticity,
see Remark~\ref{rem-plasticity} below.  Although we focus on the
basic visco-elastodynamics, the coupling with other phenomena such as damage
or aging (like \cite{Roub23SPTC}) or, e.g., also diffusion in poroelasticity
or magnetoelasticity, is conceptually well possible and would expand the
applicability of this model and approach. Even in the basic scenario
presented here, we improve \cite{Roub23SPTC} where only a semi-compressible
model and a thermally decoupled free energy were considered.

In Section~\ref{sec-linearized} we formulate the linearized
convective model involving a small-strain tensor at large displacements.
We use the objective strain rate due to the Zaremba-Jaumann time derivative,
which is most commonly used in linearized Eulerian models for solid mechanics.
It should be emphasized that this model is formulated
entirely in terms of rates, so that neither the displacement
nor the deformation occurs explicitly in the resulting system.
In Section~\ref{sec-discrtete}, we will devise a fully coupled
implicit time discretization. Finally, in Section~\ref{sec-anal},
we perform the stability and convergence analysis in a special (but
commonly used) case when the heat capacity is temperature dependent only. 

In order to facilitate the rigorous analysis particularly
for models of solids (and to model various dispersion of the velocity of
propagation of elastic waves, as analyzed in \cite{Roub24SGTL}),
some higher-order gradients in the dissipative part of the
models can be considered. These higher (here 2nd order) gradients lead to
the concept of (here 2nd-grade)
{\it nonsimple media}, which has been discussed in literature since
the works by
R.A.\,Toupin \cite{Toup62EMCS} and R.D.\,Mindlin \cite{Mind64MSLE}.
In the dissipative part as used in this paper, 
it was also developed by J.\,Ne\v cas at al.\
\cite{BeBlNe92PBMV,NeNoSi89GSIC,NecRuz92GSIV} as {\it multipolar fluids}
and later e.g.\ in \cite{FriGur06TBBC,PoGiVi13SPBC}.

For readers' convenience, let us summarize the basic notation used in what
follows:
\vspace*{-.9em}
\begin{center}
\fbox{
\begin{minipage}[t]{17em}\small\smallskip
$\yy$ deformation,\\
  $\vv$ velocity,\\
$\varrho$ mass density,\\
$\pp=\varrho\vv$ the linear momentum,\\
$\theta$ (absolute) temperature,\\
$\Ee$ small strain,\\
 $\Ent$ thermal part of the internal energy,\\
$\Ke,\Ge$ elastic bulk and shear moduli,\\
$\Kv,\Gv$ viscosity bulk and shear moduli,\\
$\Gm$ Maxwellian viscosity modulus,\\
$\HYPER$ the hyper-viscosity coefficient,\\
$I=[0,T]$ a time interval, $T>0$,\\
$\psi:\R^{3\times 3}\to\R$ free energy,\\[.1em]
$(\cdot)'$ (partial) derivative of a  mapping,\\[.1em]
$(\cdot)\!\DT{^{\,}}$ convective time derivative,
\end{minipage}
\begin{minipage}[t]{21em}\small\smallskip
$\TT$ Cauchy stress,\\
$\DD$ dissipative stress,\\
$\Lp$ inelastic distortion rate,\\
$\GRAVITY$ gravity acceleration,\\
$\bbI$ the identity matrix,\\
tr$(\cdot)$ trace of a matrix,\\
sph$(\cdot)$ spherical 
part of a matrix, sph\,$E:=$(tr\,$E)\bbI/3$,\\
dev$(\cdot)$ deviatoric part of a matrix, dev$E:=E{-}{\rm sph}E$,\\
$\mathscr{R}:\Rsym{\times}\R\to\Rdev$ the inelastic-strain rate,\\
$\R_{\rm sym}^{3\times3}$ set of symmetric matrices,\\
$\R_{\rm dev}^{3\times3}=\{A\in\R_{\rm sym}^{ 3\times 3};\ {\rm tr}A=0\}$,\\
``$\:\Cdot\:$'', ``$\:\Colon\:$'' scalar products of vectors or matrices,\\ 
``$\,\Vdots\,$'' scalar products 3rd-order tensors,\\ 
$\tau>0$ a time step for discretization,\\
$(\cdot)\!\ZJ{^{\,}}$ Zaremba-Jaumann corotational time derivative.
\smallskip \end{minipage}
}\end{center}
\nopagebreak
\vspace{-.9em}
\nopagebreak
\begin{center}
{\small\sl Table\,1.\ }
{\small
Summary of the basic notation used. 
}
\end{center}

We will use the standard notation for the Lebesgue and the Sobolev
spaces of functions on the Lipschitz bounded domain $\varOmega\subset\R^3$,
namely $L^p(\varOmega;\R^n)$ for Lebesgue measurable $\R^n$-valued functions
$\varOmega\to\R^n$ whose Euclidean norm is integrable with $p$-power, and
$W^{k,p}(\varOmega;\R^n)$ for functions from $L^p(\varOmega;\R^n)$ whose
all derivatives  up to the order $k$ have their Euclidean norm integrable with
$p$-power.  For short,  we also write $H^k=W^{k,2}$. We have the embedding
$H^1(\varOmega)\subset L^6(\varOmega)$. Furthermore, for a Banach space $X$
and for $I=[0,T]$, we will use the notation $L^p(I;X)$ for the Bochner
space of Bochner measurable functions $I\to X$ whose norm is in $L^p(I)$, and
$H^1(I;X)$ for functions $I\to X$ whose distributional derivative is in
$L^2(I;X)$. Occasionally, we will use $L_{\rm w}^p(I;X)$  for  the
space of weakly* measurable mappings $I\to X$ if $X$ has a predual,
i.e.\ there is $X'$ such that $X=(X')^*$ where $(\cdot)^*$ denotes the dual
space. The space of continuous functions on the closure $\barOmega$ of
$\varOmega$ will be denoted by $C(\barOmega)$.

\section{Linearized large-deformation convective model}\label{sec-linearized}

Most materials typically cannot withstand too much large elastic strains
without initiating inelastic processes such as damage, creep, or
plastification. Thus, elastic strains always remain rather small (and the
elastic distortion tensor is close to the identity tensor $\bbI$),
and a common  small-strain linearization is well acceptable and widely used
in many applications.  However, small strains do not preclude large
displacements, which typically occur in fluids but also in solids. In the
latter case large displacements can  occur especially  when the Kelvin-Voigt
model is combined with Maxwellian-type rheology in the deviatoric part. This
suggests the use of the Eulerian small-strain models combined with a properly
designed transport of the strain tensor, in addition to the usual density
transport. Ultimately, any rational model must respect thermodynamic
consistency in terms of mass, momentum, and energy.

This compromise between small elastic strains and large deformations
and displacements requires appropriate formulations in a convected
coordinate system. In particular, it requires the proper choice and treatment
of {\it objective rates}. Here, objectivity means that the time
derivatives do not depend on the evolving reference frame.
Many possibilities are used in the literature for different models.
For applications in solid mechanics, it is reasonable to require that the
tensor time derivatives, i.e.\ the tensor rates (in particular the stress
rate), be identical and {\it corotational}. This means 
respectively  that the stress rate vanishes for all rigid body motions
and, roughly speaking, only takes into account the rotation of the material
element. Importantly, practically all corotational rates enjoy a
chain-rule property for isotropic tensor functions, which is an attribute
not shared by other, merely objective rates, as discussed in
\cite{NHAB25HECE}.  The simplest corotational variant is the
{\it Zaremba-Jaumann time derivative} \cite{Jaum11GSPC,Zare03FPTR},
denoted it by a circle $(\cdot)\!\ZJ{^{\,}}$, which is well justified when
applied to the Cauchy stress tensor, as demonstrated by Biot
\cite[p.494]{Biot65MID}, cf.\ also
\cite{Bruh09EEBI,Fial20OTDR,MorGio22OREM}. This choice is most commonly used
especially in geophysics and more generally in engineering, See, for example,
the monographs \cite[Chap.12]{Gery19INGM} or also \cite[Sec.8.6]{HasYam13IFST}
or \cite[Sec.8.3]{Maug92TPF}. However, in some other applications in which 
cycling regimes are expected, this choice may exhibit undesirable
``ratchetting'' effects \cite{JiaFis17ADRD,MeXiBrMe03ESRC}.

In this paper, we focus on {\it isotropic materials}, for which this
derivative also affects the symmetric small-strain tensor $\Ee$, which is 
considered as an independent variable in the models. For an overview of
objective corotational strain rates see \cite{XiBrMe98SRMS,MeShBr00SCOR}.
Specifically, given the Eulerian velocity $\vv$, we define 
\begin{align}\label{ZJ}
\ZJ\Ee=
\pdt\Ee+(\vv\Cdot\nabla)\Ee-\WW\Ee+\Ee\WW\ \ \ \text{ with }\
\WW={\rm skw}(\nabla\vv)=\frac{\nabla\vv{-}(\nabla\vv)^\top}2\,.
\end{align}
The Eulerian velocity $\vv$ is also used in the convective time derivative
for scalars and, component-wise, for vectors or tensors
\begin{align}
(\bm\cdot)\!\DT{^{}}=\pdt{}(\bm\cdot)+(\vv{\cdot}\nabla)(\bm\cdot)\,.
\end{align}
Then \eq{ZJ} can be written shortly as
$\ZJ\Ee=\DT\Ee-{\rm skw}(\nabla\vv)\Ee+\Ee{\rm skw}(\nabla\vv)$.

\subsection{The thermodynamical system and its energetics}
\label{sec-linearized-system}

We consider the visco-elastodynamics in the {\it Kelvin-Voigt rheology} in
the volumetric part and the {\it Jeffreys} (also called {\it anti-Zener})
{\it rheology} in the deviatoric (isochoric) part. This is a fairly general
model that allows for an isochoric creep or plasticity. For the
displacement $\uu$, we implement the {\it Green-Naghdi additive decomposition}
\cite{GreNag65GTEP} of the total small strain $\strain(\uu)$ into the elastic
and the inelastic strain $\strain(\uu)=\Ee+\Ep$, where $\Ep$ is in the
position of an internal variable. This decomposition will be expressed in
objective rates as $\strain(\vv)=\ZJ\Ee+\ZJ\Ep$, to be further written as
\begin{align}\label{additive-rate}
\strain(\vv)=\ZJ\Ee+\Lp\ \ \text{ with }\ \ \Lp:=\ZJ\Ep\,,
\end{align}
where $\Ep$ denotes the inelastic strain. Thus both $\uu$ and $\bm{P}$ can
be eliminated in the end, although these variables can be reconstructed
a-posteriori if the corresponding initial conditions were prescribed.
Therefore, the model allows for {\it large displacements} while still using
small elastic strains. 

The main ingredient is the Helmholtz' free energy
$\psi:\Rsym{\times}\R^+\to\R$ acting on the small-strain tensor
$\bm{E}\in\Rsym$ and on temperature $\theta$. Here, and in what
follows, we will not notationally distinguish the variable as a placeholder
and a time-space field. Here, it means both $\bm{E}\in\Rsym$ and 
$\Ee:I{\times}\varOmega\to\Rsym$. Standardly, we define the (conservative
part of the) Cauchy stress tensor $\TT$ as
\begin{align}\label{stress-entropy-Ch4}
\TT=\mathscr{T}(\Ee,\theta)\ \ \ \text{ with }\ 
\mathscr{T}(\bm{E},\theta):=\psi'_{\bm{E}}(\bm{E},\theta)+\psi(\bm{E},\theta)\bbI\,.
\end{align}
We note that the ``negative-pressure'' term $\psi(\Ee,\theta)\bbI$ in
\eq{stress-entropy-Ch4} balances the energetics in the convected model,
cf.\ the calculations \eq{Euler-small-divT.v++}--\eq{Euler-small-calc} below. 
From the free energy $\psi$, we also obtain (by Gibbs' relation) the
{\it internal energy} $\ENG$ and the {\it entropy} $\eta$, and also the
{\it  heat part} $\ENT$ of $\ENG$ as
\begin{subequations}\label{def-E-eta-U}\begin{align}\label{def-E-eta}
&\ENG(\bm{E},\theta)=\psi(\bm{E},\theta)+\theta\eta(\bm{E},\theta)
\ \text{ with }\ \eta(\bm{E},\theta)=-\psi_\theta'(\bm{E},\theta)
\ \text{ and }\
\\&\ENT(\bm{E},\theta)=\ENG(\bm{E},\theta)-\ENG(\bm{E},0)\,.
\label{def-U}\end{align}\end{subequations}

Another ingredient is the (pseudo)potential of the dissipative forces
$\zeta:\R^+{\times}\Rsym{\times}\Rdev\to\R$ acting on the rates
$\EE:=\strain(\vv)$ and $\Lp$, and also depending on $\theta$.
This potential governs the dissipative contribution to the
Cauchy stress $\DD=\zeta_\EE'(\theta,\strain(\vv),\Lp)$. Based on the
Kelvin-Voigt rheology summing the stresses $\TT$ and $\DD$, we can
set up the {\it momentum equation}
\begin{align}\label{Euler-small-viscoelastodyn+1}
\varrho\DT\vv={\rm div}\big(\TT{+}\DD\big)+\varrho\GRAVITY
\end{align}
with $\varrho$ the mass density which is governed by the mass {\it continuity
equation}
\begin{align}\label{Euler-small-viscoelastodyn+0}
&\!\!\DT\varrho=-\varrho\,{\rm div}\,\vv\,.
\end{align}
From the potential $\zeta$, we also obtain a ``creep stress'' $\zeta_\Lp'$
which balances the  deviatoric part of $\TT$ and which stands here
in the position of a ``linearized  Mandel stress''.
Thus we compose the flow rule for the creep rate as
\begin{align}\label{def-Pi-small}
\zeta_{\Lp}'(\theta;\strain(\vv),\Lp)=
{\rm dev}\,\TT\,.
\end{align}
Of course, the pressure contribution  $\psi\,\bbI$ of $\TT$ in
\eq{stress-entropy-Ch4}  is purely volumetric and thus has no
effect on \eq{def-Pi-small}.

The entropy, considered in Pa/K\,=\,J/(m$^3$K), is an intensive variable and
is governed by the {\it entropy equation}:
\begin{align}
\pdt\eta+{\rm div}(\eta\vv)=\frac{\xi-{\rm div}\,{\bm j}}\theta\ \ \ \
\text{ with the heat flux }\ \ \jj=-\kappa(\theta)\nabla\theta\,,
\label{entropy-eq-Ch4}\end{align}
where $\xi=\xi(\theta;\strain(\vv),\ZJEp)$ denotes the heat production rate
(=\,the specific power in W/m$^3$) and is considered
equal to the dissipation rate of the mechanical energy
\begin{align}
\xi(\theta;\strain(\vv),\ZJEp)&=
\zeta_\EE'(\theta,\strain(\vv),\Lp)\Colon\strain(\vv)
+\zeta_{\Lp}'(\theta;\strain(\vv),\Lp)\Colon\Lp
=\DD\Colon\strain(\vv)+\TT\Colon\Lp.
\end{align}
The latter equality in \eq{entropy-eq-Ch4} is the {\it Fourier law} which
phenomenologically determines the heat flux $\jj$  to be  proportional
to the negative temperature gradient through the
thermal conductivity coefficient $\kappa=\kappa(\theta)$ considered
independent of $\Ee$; for a general continuous dependence of $\kappa$ also
on $\Ee$ see Remark~\ref{rem-Euler-thermo-kappa} below. Assuming $\xi\ge0$ and
$\kappa\ge0$ and integrating \eq{entropy-eq-Ch4} over the domain $\varOmega$
while imposing the impenetrability of the boundary in the sense that the
normal velocity $\vv{\cdot}\nn$ vanishes over the boundary $\varGamma$ of
$\varOmega$, we obtain the {\it Clausius-Duhem inequality}:
\begin{align}
\frac{\d}{\d t}\!\!\!\!\!\!\lineunder{\int_\varOmega\eta\,\d\xx}{total entropy}\!\!\!\!\!\!
=\int_\varOmega\!\!\!\!\!\!\!\!\!\!\lineunder{\frac\xi\theta+\kappa(\theta)\frac{|\nabla\theta|^2\!}{\theta^2}}{entropy production rate}\!\!\!\!\!\!\!\!\!\!\d\xx
+\!\!\int_\varGamma\!\!\!\!\!\lineunder{\!\Big(\kappa(\theta)\frac{\nabla\theta}{\theta}-\eta\vv\Big)\!}{entropy flux}\!\!\!\!\!\Cdot\nn\,\d S
\ge\!\int_\varGamma\!\kappa(\theta)\frac{\!\nabla\theta\!}{\theta}{\cdot}\nn\,\d S.
\label{entropy-ineq-Ch4}\end{align}
If the system is thermally isolated in the sense that the normal heat flux
$\jj{\cdot}\nn$ vanishes across the boundary $\varGamma$, we recover the
{\it 2nd law of thermodynamics}, i.e., the total entropy in isolated
systems is nondecreasing with time.

Restricting our attention to a case when
\begin{align}\text{$\TT$ and $\Ee$ commute,}
\label{ass-T-E-commute}\end{align}
we perform the calculus for the time derivative of the internal
energy $\Eng=\ENG(\Ee,\theta)$ from \eq{def-U} as an extensive variable:
\begin{align}\nonumber
  \pdt\Eng&+{\rm div}(\Eng\vv)=
 \DT\Eng+\Eng\,{\rm div}\,\vv
 =\ \DT{\overline{\psi(\Ee,\theta)-\theta\psi_\theta'(\Ee,\theta)}}
\,+\big(\psi(\Ee,\theta)-\theta\psi_\theta'(\Ee,\theta)\big)\,{\rm div}\,\vv
\\[-.2em]&\nonumber\ \ =\
\psi'_{E}(\Ee,\theta){:}\DT\Ee+\psi_\theta'(\Ee,\theta)\DT\theta
-\theta\psi_{E\theta}''(\Ee,\theta)\DT\Ee-\theta\psi_{\theta\theta}''(\Ee,\theta)\DT\theta
-\DT\theta\psi_\theta'(\Ee,\theta)
\\[-.2em]&\nonumber\hspace{18.7em}
+\big(\psi(\Ee,\theta)-\theta\psi_\theta'(\Ee,\theta)\big)\,{\rm div}\,\vv
\\[-.6em]&\nonumber\ \ =\ \psi'_{E}(\Ee,\theta){:}\DT\Ee+\theta\,\DT{\overline{\eta(\Ee,\theta)}}
+\big(\psi(\Ee,\theta)-\theta\psi_\theta'(\Ee,\theta)\big)\,{\rm div}\,\vv
\\[-.2em]&\nonumber\stackrel{\scriptsize\eq{entropy-eq-Ch4}}{=}\!\psi'_{E}(\Ee,\theta){:}\DT\Ee
+\xi-{\rm div}\,\jj-\theta\eta(\Ee,\theta){\rm div}\,\vv
+\big(\psi(\Ee,\theta)-\theta\psi_\theta'(\Ee,\theta)\big)\,{\rm div}\,\vv
\\[-.2em]&\nonumber
\stackrel{\scriptsize\eq{additive-rate}}{=}
\!\psi'_{E}(\Ee,\theta){:}\big(\strain(\vv)
+{\rm skw}(\nabla\vv)\Ee+\Ee{\rm skw}(\nabla\vv)
-\ZJEp\big)+\xi-{\rm div}\,\jj+\psi(\Ee,\theta){\rm div}\,\vv
\\[-1em]&
\!\!\UUU{=}{\eq{E:WS=E:SW}}{\eq{stress-entropy-Ch4}}\!\!
\xi-{\rm div}\,\jj
+\TT{:}\big(\strain(\vv){-}\ZJEp\big)
\,,\label{Euler-thermodynam3-Ch4-}\end{align}
where we have used the  matrix  algebra
$A\Colon (BC)=(B^\top\!A)\Colon C=(AC^\top)\Colon B$ so that,
for $W={\rm skw}\,L$, it holds 
\begin{align}\nonumber
S\Colon(WE-EW)&=\frac12S\Colon\big(LE-L^\top E-EL+EL^\top\big)
=\frac12\big(SE^\top\!\!-E^\top S\big)\Colon L
-\frac12\big(LS-SL\big)\Colon E
\\&=\frac12\big(SE^\top\!\!-ES^\top\!\!-E^\top S+S^\top E\big)\Colon L
=\big(\!\!\!\!\lineunder{SE-ES}{$\ =0$}\!\!\!\!\big)\Colon L=0\,.
\label{E:WS=E:SW}
\end{align}
Note that in \eq{Euler-thermodynam3-Ch4-}, we have used \eq{E:WS=E:SW}
for $\psi'_{\bm{E}}(\Ee,\theta)$, $\Ee$, $\Nabla\vv$, and
${\rm skw}(\Nabla\vv)$ in place of $S$, $E$, $L$, and $W$, respectively, so
that, in particular, the last equality in \eq{E:WS=E:SW} exploited
\eq{ass-T-E-commute}. Moreover, it should be emphasized that in
\eq{E:WS=E:SW}, we have to assume that the
initial condition for $\Ee$ is symmetric and to exploit that the
Zaremba-Jaumann corotational derivative in \eq{additive-rate} maintains
the symmetry of $\Ee$ throughout the entire evolution.

In fact, \eq{Euler-thermodynam3-Ch4-} with $\TT=\mathscr{T}(\Ee,\theta)$
forms the {\it internal-energy equation} for
$\Eng=\ENG(\Ee,\theta)$ with $\ENG$ from \eq{def-E-eta}. Similarly, by
making an obvious modification of \eq{Euler-thermodynam3-Ch4-}, we can
obtain the equation for the heat part of the internal energy
$\Ent=\ENT(\Ee,\theta)$:
\begin{align}\nonumber
\!\!\!\!\pdt\Ent+{\rm div}\big(\Ent\vv{+}{\bm j}\big)&=\xi
+\big(\mathscr{T}(\Ee,\theta){-}\mathscr{T}(\Ee,0)\big)\Colon
\big(\strain(\vv){-}\Lp\big)\ \
\text{ with }\ \Ent=\ENT(\Ee,\theta)\ \text{ and}
\\\hspace{5em}
\ \text{ with }\ 
\jj&=-\kappa(\theta)\nabla\theta\,,\ \
\text{ where $\mathscr{T}$ is from \eq{stress-entropy-Ch4} and
$\ENT$ from \eq{def-E-eta-U}}\,.
\label{Euler-thermodynam3-Ch4}\end{align}

\subsection{The thermomechanical system with a multipolar viscosity}

The momentum equation \eq{Euler-small-viscoelastodyn+1}, the mass continuity
equation \eq{Euler-small-viscoelastodyn+0}, the  kinematic  equation
\eq{additive-rate} reflecting the Green-Naghdi decomposition, the creep flow
rule \eq{def-Pi-small}, and the heat internal energy equation
\eq{Euler-thermodynam3-Ch4} altogether form a closed system for
$(\varrho,\vv,\Ee,\ZJEp,\theta)$.  

We restrict our focus to the thermal coupling only in the volumetric
part and on the additive split of the Stokes and the Maxwellian viscosities,
specifically
\begin{align}\label{Euler-small-ass}
{\rm dev}\big(\psi''_{\bm{E}\theta}(\Ee,\theta)\big)\equiv0
\ \ \text{ and }\ \ 
\zeta(\theta,\EE,\Lp)=\tfrac12\bbD\EE\Colon\EE+\zeta_{\rm p}(\theta,\Lp)\,.
\end{align}
The former means that temperature may influence only the volumetric part of
the Cauchy stress but does not influence the deviatoric part. Indeed,
realizing that
\begin{align}\nonumber
\mathscr{T}(\Ee,\theta)=\psi'_{\bm{E}}(\Ee,\theta)+\psi(\Ee,\theta)\bbI
&=\psi'_{\bm{E}}(\Ee,0)+\psi(\Ee,0)\bbI
+\int_0^\theta\!\!\!\psi''_{\bm{E}\theta}(\Ee,\vartheta)
+\psi_\theta'(\Ee,\vartheta)\bbI\,\d\vartheta
\\[-.6em]&=\mathscr{T}(\Ee,0)+\!\!\int_0^\theta\!\!\!\psi''_{\bm{E}\theta}(\Ee,\vartheta)
+\psi_\theta'(\Ee,\vartheta)\bbI\,\d\vartheta
\,,
\end{align}
the  former condition in \eq{Euler-small-ass} implies that
\begin{align}\label{Eular-small-devT0=devT}
{\rm dev}\,\mathscr{T}(\Ee,\theta)={\rm dev}\,\mathscr{T}(\Ee,0)\,.
\end{align}
Notably, \eq{Euler-small-ass} simplifies the adiabatic heat source in
\eq{Euler-thermodynam3-Ch4} as
\begin{align}\nonumber
\big(\mathscr{T}(\Ee,\theta)\,{-}\mathscr{T}(\Ee,0)\big)
\Colon\big(\strain(\vv){-}\ZJEp\big)
&={\rm sph}\big(\mathscr{T}(\Ee,\theta)\,{-}\mathscr{T}(\Ee,0)\big)\Colon
\strain(\vv)
\\&=\tfrac13{\rm tr}\big(\mathscr{T}(\Ee,\theta)\,{-}\mathscr{T}(\Ee,0)\big)\,
{\rm div}\,\vv\,
\label{Eular-small-adiabatic-power}\end{align}
because the inelastic strain (and its rate) naturally does not affect the
volumetric part of the model and is thus modelled as isochoric, i.e.,
${\rm tr}\,\ZJEp=0$. Furthermore, for the latter in \eq{Euler-small-ass}, we
assume (as physically natural) the convexity of $\zeta_{\rm p}(\theta,\cdot)$
and use the standard convex-analysis construction of the convex conjugate
$\zeta_{\rm p}^*(\theta,\cdot):=[\zeta_{\rm p}(\theta,\cdot)]^*$. Realizing that
$\zeta_{\rm p}^*(\theta,\cdot)'=[\zeta_{\rm p}(\theta,\cdot)']^{-1}$, we will use
the abbreviation for the inelastic-strain rate
$\mathscr{R}:\Rsym{\times}\R\to\Rdev$ defined by 
\begin{align}\label{R-abbreviation}
\mathscr{R}:(\bm{E},\theta)\mapsto\big[{\zeta_{\rm p}(\theta,\cdot)^*}\big]'
\big({\rm dev}\mathscr{T}(\bm{E},\theta)\big)\,.
\end{align}

In addition, as mentioned in the Introduction, we incorporate a nonlinear
dissipative higher-order enhancement of the Stokes-type viscosity, also
known as {\it hyper-viscosity}. For simplicity,
we adopt this concept of {\it multipolar continua} equally in
the volumetric and in the shear parts using just one coefficient $\HYPER>0$
and one exponent $p$, leading to the extended dissipation potential
and the extended dissipation rate
\begin{align}\nonumber
&\!\!\!\zeta_{\rm ext}^{}(\theta;\vv,\ZJEp)
=\zeta(\theta,\strain(\vv),\Lp)+\tfrac\HYPER p|\nabla^2\vv|^p
\!\stackrel{\eq{Euler-small-ass}}{=}\!
\tfrac12\bbD\strain(\vv)\Colon\strain(\vv)+
\zeta_{\rm p}(\theta,\Lp)+\tfrac\HYPER p|\nabla^2\vv|^p\ 
\text{ and}
\!\\&
\!\!\!\xi_{\rm ext}^{}(\theta;\vv,\ZJEp)=\xi(\theta,\strain(\vv),\Lp)+\HYPER|\nabla^2\vv|^p
\!\stackrel{\eq{Euler-small-ass}}{=}\!\bbD\strain(\vv)\Colon\strain(\vv)
+[\zeta_{\rm p}]'_{\Lp}(\theta,\Lp)\Colon\Lp+\HYPER|\nabla^2\vv|^p\,,
\label{xi-ext}\end{align}
respectively. 
Mechanically, such higher gradient only amplifies the effect of normal
dispersion and attenuation of elastic waves as exhibited already for simple
Stokes viscosity in the Kelvin-Vogt rheology. In other words, waves with
ultra-high frequencies are more attenuated and do not propagate at all,
while waves with very low frequencies can propagate without substantial
attenuation, cf.\ \cite[Sect.3.1]{Roub24SGTL}. When $\zeta_{\rm p}\not\equiv0$,
the contribution from the Maxwellian-type viscosity combines with the
anomalous dispersion, so that also ultra low frequency waves are also
attenuated and eventually cannot propagate at all, cf.\ again
\cite{Roub24SGTL}.

We express the system in terms of the linear momentum $\pp=\varrho\vv$ and,
later, exploit the convexity of the kinetic energy $(\pp,\varrho)\mapsto
\frac12|\pp|^2/\varrho$, in contrast to the  nonconvex (although equivalent)
form $(\vv,\varrho)\mapsto\frac12\varrho|\vv|^2$, cf.\
\cite{Roub25TDVE}. The explicit use of the convexity of the kinetic energy
expressed in terms of the momentum $\pp$ is also in \cite{FLMS21NACF}. 
Exploiting that
$\varrho\DT\vv=\pdt{}(\varrho\vv)+{\rm div}(\varrho\vv{\otimes}\vv)$ granted
by \eq{Euler-small-viscoelastodyn+0}, we will deal with the system for
$(\varrho,\vv,\Ee,\theta)$
and thus also for $\pp$ and $\Ent$:
\begin{subequations}\label{Euler-small-therm-ED-anal}
\begin{align}\label{Euler-small-therm-ED-anal0}
&
\!\!\pdt\varrho=-{\rm div}\,\pp\ \ \ \text{ with $\ \ \pp=\varrho\vv$}\,,
\\&\!\!\nonumber
\pdt\pp={\rm div}\big(\mathscr{T}(\Ee,\theta){+}\DD
{-}\pp{\otimes}\vv\big)+\varrho\GRAVITY
\ \ \text{ with }\ \mathscr{T}
\ \text{ from \eq{stress-entropy-Ch4}}
\\[-.3em]&\hspace{9em}\text{and }\
\DD=\bbD\strain(\vv)-{\rm div}\mathfrak{H}\ \text{ with }\ 
\mathfrak{H}=\HYPER\big|\nabla^2\vv\big|^{p-2}\nabla^2\vv\,,
\label{Euler-small-therm-ED-anal1}
\\[-.4em]&\!\!
\pdt\Ee=\strain(\vv)-\ZJEp-\bm B_\text{\sc zj}^{}(\vv,\bm E)
\ \ \text{with }\ \ZJEp=\mathscr{R}(\Ee,\theta)\ \text{ where $\mathscr{R}$ is
from \eq{R-abbreviation}}\,,
\label{Euler-small-therm-ED-anal2}
\\&\nonumber
\pdt\Ent={\rm div}\big(\kappa(\theta)\nabla\theta-\Ent\vv\big)
+\xi_{\rm ext}^{}(\theta;\vv,\ZJEp)
+\tfrac13{\rm tr}\big(\mathscr{T}(\Ee,\theta){-}\mathscr{T}(\Ee,0)\big)
\,{\rm div}\,\vv\ \ \text{with }\ \
\\[-.3em]&\hspace{2.6em}
\xi_{\rm ext}^{}(\theta;\vv,\ZJEp)
\ \text{ from \eq{xi-ext}}
\ \text{ and}
\ \ \Ent=\ENT(\Ee,\theta)
:=\psi(\Ee,\theta)-\theta\psi_\theta'(\Ee,\theta)-\psi(\Ee,0)\,,
\label{Euler-small-therm-ED-anal3}
\end{align}\end{subequations}
where, in \eq{Euler-small-therm-ED-anal2}, we have used a shorthand notation
for the bi-linear operator involved in the Zaremba-Jaumann
derivative:
 \begin{align}
 \bm B_\text{\sc zj}^{}(\vv,\bm E)=(\vv{\cdot}\nabla)\bm E-
      {\rm skew}(\nabla\vv)\Ee+\Ee\,{\rm skew}(\nabla\vv)\,.
\label{ZJ-def}\end{align}
The equations \eq{Euler-small-therm-ED-anal1} and
\eq{Euler-small-therm-ED-anal3} are to be completed by the
boundary conditions. For velocity, we allow for free slip
but, as is common in the Eulerian formulation, we prescribe zero
normal velocity (so that the shape of
$\varOmega$ does not evolve). This is here adapted for the multipolar
4th-order enhancement. For the temperature, we consider
the Fourier-type boundary condition. Altogether,
\begin{align}\nonumber
&\vv\Cdot\nn=0,\ \ \ \ [(\mathscr{T}(\Ee,\theta){+}\DD)\nn+\divS(\mathfrak{H}\nn)]_\text{\sc t}^{}
=\bm0\,,\ \ \ \ 
\mathfrak{H}\Colon(\nn{\otimes}\nn)={\bm0}\,,\ \ \ \ 
\text{ and }
\\&
\kappa(\theta)\nabla\theta\Cdot\nn+h(\theta)=h_{\rm ext}\ \ \ \text{ on }\ \ \varGamma\,
\label{Euler-small-therm-ED-BC-hyper}\end{align}
with $h_{\rm ext}$ a prescribed external heat flux and
with $h(\cdot)$ a temperature-dependent boundary heat  out-flux.
We consider the initial-value problem for
\eq{Euler-small-therm-ED-anal} with the initial conditions
\begin{align}\label{Euler-Ch4-thermodynam-IC}
\varrho|_{t=0}^{}=\varrho_0\,,
\ \ \ \ \ \vv|_{t=0}^{}=\vv_0\,,\ \ \ \ \ 
\Ee|_{t=0}^{}=\Ee_0\,,\ \,\text{ and }\,\ \theta|_{t=0}^{}=\theta_0\,.
\end{align}

Now, let us reveal the energy balance behind this system.
To this aim, we test \eq{Euler-small-therm-ED-anal0}
by $|\pp|^2/(2\varrho^2)$ and \eq{Euler-small-therm-ED-anal1}
by $\pp/\varrho$, which yields
\begin{align}
  \pdt{}\bigg(\frac{|\pp|^2}{2\varrho}\bigg)
  &=\frac{\pp}{\varrho}\Cdot\pdt\pp
  -\frac{|\pp|^2}{2\varrho^2}\pdt\varrho
=\frac{\pp}{\varrho}\Cdot
\Big({\rm div}\big(\TT{+}\DD-\vv{\otimes}\pp\big)+\varrho\GRAVITY\Big)
+\frac{|\pp|^2}{2\varrho^2}{\rm div}\,\pp
\,.\label{rate-of-kinetic-modif}\end{align}
Then, we use the calculus
\begin{align}\label{calculus-for-kinetic}
\!\!\int_\varOmega\frac{\pp}{\varrho}\Cdot{\rm div}\big(\pp{\otimes}\vv\big)
-\frac{|\pp|^2\!}{2\varrho^2\!}\,{\rm div}\,\pp\,\d\xx
=\int_\varOmega\!\vv\Cdot{\rm div}\big(\varrho\vv{\otimes}\vv\big)
-\frac{|\vv|^2\!}{2}\,{\rm div}(\varrho\vv)\,\d\xx=0\,,
\end{align}
based on the Green formula with $\vv\Cdot\nn=0$ for the calculus 
\begin{align}\nonumber
\!\!\!\!\!\int_\varOmega\varrho(\vv\Cdot\nabla)\vv\Cdot\vv\,\d\xx
&=\int_\varGamma\varrho|\vv|^2\!\!\!\!\lineunder{\!\!\!\!\vv\Cdot\nn\!\!\!\!}{$=0$}\!\!\!\!\d S
-\!\!\int_\varOmega\!\vv\Cdot{\rm div}(\varrho\vv{\otimes}\vv)\,\d\xx
\\[-.7em]&
=-\!\!\int_\varOmega\varrho(\vv\Cdot\nabla)\vv\Cdot\vv+{\rm div}(\varrho\vv)|\vv|^2\,\d\xx
=-\int_\varOmega\!{\rm div}(\varrho\vv)\frac{|\vv|^2}2\,\d\xx\,.
\label{calculus-convective}\end{align}

\def\SS{{\bm{S}}}

Let us denote the mere {\it stored energy} part of $\psi$ as
$\varphi(\bm{E}):=\psi(\bm{E},0)$ and the corresponding part
of the Cauchy stress tensor $\TT$ as
$\TT_0=\mathscr{T}(\Ee,0)=\varphi'(\Ee)+\varphi(\Ee)\bbI$.
For this part of the Cauchy stress, we use the calculus
\begin{align}\nonumber
\int_\varOmega&\!({\rm div}\,\TT_0)\Cdot\vv\,\d\xx
=\!\int_\varGamma\!\!\vv\Cdot\TT_0\nn\,\d S
-\!\int_\varOmega\!\TT_0\Colon\strain(\vv)\,\d\xx
\!\stackrel{\eq{Euler-small-viscoelastodyn+1}}{=}\!\!
\int_\varGamma\!\vv\Cdot\TT_0\nn\,\d S
-\!\int_\varOmega\!
\big(\varphi'(\Ee){+}\varphi(\Ee)\bbI\big)\Colon\strain(\vv)\,\d\xx
\\&\nonumber\!\stackrel{\eq{ZJ}}{=}\!\int_\varGamma\!\!\vv\Cdot\TT_0\nn\,\d S
-\!\!\int_\varOmega\!\varphi'(\Ee)\Colon\Big(\pdt{\Ee\!}
+(\vv\Cdot\Nabla)\Ee-\WW\Ee+\Ee\WW\Big)
-{\rm dev}\TT_0\Colon\mathscr{R}(\Ee,\theta)
+\varphi(\Ee){\rm div}\,\vv\,\d\xx
\\[-1.1em]&\!\!\!\!\UUU{=}{\eq{Euler-small-calc}}{\eq{E:WS=E:SW}}\!\!\!\!
\int_\varGamma\!\!\vv\Cdot\TT_0\nn\,\d S
-\frac{\d}{\d t}\int_\varOmega\varphi(\Ee)\,\d\xx
-\int_\varOmega 
[\zeta_{\rm p}]'_{\Lp}(\theta,\Lp)\Colon\Lp\,\d\xx\,.
\label{Euler-small-divT.v++}\end{align}
Here, in addition to \eq{Eular-small-adiabatic-power}, we have also used that
\begin{align}\nonumber
{\rm dev}\,\TT_0\Colon\mathscr{R}(\Ee,\theta)\!
&\stackrel{\eq{Eular-small-devT0=devT}}{=}\!
{\rm dev}\mathscr{T}(\Ee,\theta)\Colon\mathscr{R}(\Ee,\theta)
\\&\ \ =
{\rm dev}\mathscr{T}(\Ee,\theta)\Colon[\zeta_{\rm p}(\theta,\cdot)^*]'({\rm dev}\mathscr{T}(\Ee,\theta))=[\zeta_{\rm p}]'_{\Lp}(\theta,\Lp)\Colon\Lp
\end{align}
and the calculus
\begin{align}\label{Euler-small-calc}
\int_\varOmega\!\varphi'(\Ee)\Colon(\vv\Cdot\Nabla)\Ee+\varphi(\Ee){\rm div}\,\vv\,\d\xx=\!\!\int_\varOmega\!\nabla\varphi(\Ee)\Cdot\vv+\varphi(\Ee){\rm div}\,\vv\,\d\xx=\!\!\int_\varGamma\!\varphi(\Ee)\!\!\lineunder{\!\!(\vv\Cdot\nn)\!\!}{$=0$}\!\!\d S\,.
\end{align}
Also, taking into account the form of the spin
$\WW=\frac12\Nabla\vv-\frac12(\Nabla\vv)^\top$, we have used 
$\varphi'(\Ee)\Colon(\WW\Ee-\Ee\WW)=\bm0$, cf.\ \eq{E:WS=E:SW}.

Relying on the boundary conditions \eq{Euler-small-therm-ED-BC-hyper},
we finally obtain the {\it energy-dissipation balance}
\begin{align}\nonumber
\!\!\!\!\!\frac{\d}{\d t}\int_\varOmega\!\!\!\!
\linesunder{\frac\varrho2|\vv|^2}{kinetic}{energy}
\!\!\!\!\!+\!\!\!\!\!\linesunder{\varphi(\Ee)}{stored}{energy}\!\!\!\!\,\d\xx
+\!\int_\varOmega
\hspace{-.7em}\linesunder{\bbD\strain(\vv)\Colon\strain(\vv)}{disipation rate due to}{the  Stokes  viscosity}\hspace{-.7em}
+\hspace{-.7em}\linesunder{\HYPER|\nabla^2\vv|^p}{disipation rate due}{to the  hyper-viscosity}\hspace{-1.5em}
+\hspace{-.2em}\linesunder{[\zeta_{\rm p}]'_{\Lp}(\theta,\Lp)\Colon\Lp}{disipation rate due to}{Maxwellian viscosity}\hspace{-1em}\d\xx\
\\[-.4em]
=\int_\varOmega\hspace{-1.4em}\linesunder{\varrho\,\GRAVITY\Cdot\vv}{\ power of}{gravity field\ }\hspace{-1.4em}-\hspace{-.9em}\lineunder{
\tfrac13{\rm tr}\big(\mathscr{T}(\Ee,\theta){-}\mathscr{T}(\Ee,0)\big)
\,{\rm div}\,\vv}{\ power of adiabatic effects}\hspace{-.9em} \,\d\xx\,;\ 
\label{Euler-small-energy-balance-stress}\end{align}
cf.\ also \cite[Formulas (2.10)-(2.13)]{Roub25TDVE} for detailed handling
of the boundary conditions by using Green's formula over $\varOmega$ twice
and the Green surface formula over $\varGamma$.

By adding \eq{Euler-small-therm-ED-anal3} tested by 1, we obtain the
{\it total-energy balance}:
\begin{align}
  &\!\!\!\!\!\frac{\d}{\d t}
  \int_\varOmega\!\!\!\!
  \linesunder{\frac\varrho2|\vv|^2}{kinetic}{energy}\!\!\!\!\!\!+\!\!\!\!
  \linesunder{
\ENG(\Ee,\theta)}{internal}{energy}\!\!\!\d\xx
+\!\int_\varGamma\!\!\!\!\!\threelinesunder{h(\theta)}{power of}{boundary}{flux}\!\!\!\!\!\!\!\d S
=\int_\varOmega\!\!\!\!\!\!\linesunder{\varrho\,\GRAVITY{\cdot}\vv}{power of}{gravity field}\!\!\!\!\!\!\!
\d\xx
+\!\int_\varGamma
\!\!\!\!\threelinesunder{h_{\rm ext}}{external\!\!}{heat}{flux}\!\!\d S\,,
\label{Euler-engr-finite}\end{align}
expressing the {\it 1st law of thermodynamics}: in particular, in isolated
systems, the total energy (here as a sum of the kinetic and the internal
energies) is conserved.

Furthermore, for the (here rather formal)  test of
\eq{Euler-thermodynam3-Ch4} by the so-called {\it coldness} $1/\theta$, we need
to assume positivity of $\theta$ and then use the calculus
\begin{align}\nonumber
\!\!\!\!\!\!\!\!\!\!\frac1\theta\Big(\pdt\Eng+{\rm div}(\Eng\vv)\Big)
&=\frac1\theta\,\DT{{\overline{\ENG(\Ee,\theta)}}}
+\frac{\ENG(\Ee,\theta)\!}{\theta}\,{\rm div}\,\vv
\\&\nonumber=\frac{\ENG_\theta'(\Ee,\theta)}{\theta}\DT\theta
+\frac{\ENG_\Ee'(\Ee,\theta)}{\theta}\Colon\DT\Ee
+\frac{\ENG(\Ee,\theta)\!}{\theta}\,{\rm div}\,\vv
\\&\nonumber=\DT{{\overline{\eta(\Ee,\theta)}}}
+\Big(\frac{\ENG_\Ee'(\Ee,\theta)}{\theta}{-}\eta_\Ee'(\Ee,\theta)
\Big)\Colon\DT\Ee+\frac{\ENG(\Ee,\theta)\!}{\theta}\,{\rm div}\,\vv
\\&\nonumber=\DT{{\overline{\eta(\Ee,\theta)}}}
+\eta(\Ee,\theta){\rm div}\,\vv
+\frac{\psi_\Ee'(\Ee,\theta)}{\theta}\Colon\DT\Ee
+\frac{\psi(\Ee,\theta)\!}{\theta}\,{\rm div}\,\vv
\\&=\pdt{}\eta(\Ee,\theta)+{\rm div}(\eta(\Ee,\theta)\vv)
+\!\!\!\lineunder{\!\!\frac{\psi_\Ee'(\Ee,\theta)}{\theta}\Colon\DT\Ee
+\frac{\psi(\Ee,\theta)\!}{\theta}\,{\rm div}\,\vv}{$=\TT(\Ee,\theta){:}(\strain(\vv){-}\ZJEp)/\theta$}\!\!\!\!\,.
\nonumber\\[-1.2em]\label{calculus-to-entropy++}\\[-2.2em]\nonumber\end{align}
Here, in addition to \eq{E:WS=E:SW}, we also used
\begin{align}
\eta_\theta'(\Ee,\theta)=\frac{\ENG_\theta'(\Ee,\theta)}\theta\ \ \ \text{ and }\
\ \ \frac{\psi(\Ee,\theta)}\theta-\frac{\ENG(\Ee,\theta)}\theta=
\psi_\theta'(\Ee,\theta)=-\eta(\Ee,\theta)\,.
\label{calculus-to-entropy+++}\end{align}
 The positivity of $\theta$ mentioned above does not seem easy to prove,
cf.\ also \cite[Remark 3.4]{Roub25TEFS}, and we leave it open here. In any
case,  by the mentioned  (formal)  test of
\eq{Euler-thermodynam3-Ch4-}, we obtain the {\it entropy balance}:
\begin{align}
  &\frac{\d}{\d t}\!\!\!\!
  \linesunder{\int_\varOmega\!\eta(\Ee,\theta)\,\d\xx\!\!}{total}{entropy}\!\!\!\!
  =\int_\varOmega\!\!\!
  \threelinesunder{\!\!
\frac{\!\xi_{\rm ext}^{}(\theta;\vv,\ZJEp)\!}\theta
}{entropy}{production due}{to viscosity}\!\!\!\!
+\!\!\!\!\!\threelinesunder{\kappa(\theta)
  \frac{|\nabla\theta|^2\!\!}{\theta^2}}{entropy
  produ-}{ ction due to}{heat transfer}\!\!\!\!\!\!\d\xx
+\int_\varGamma\!\!\!
\threelinesunder{\!\!\!\frac{h_{\rm ext}-h(\theta)}\theta\!\!\!}{entropy flux}{through}{boundary}\!\d S
  \label{Euler-entropy-finite++}\end{align}
with the dissipation rate $\xi_{\rm ext}^{}(\theta;\vv,\ZJEp)$ from
\eq{Euler-small-therm-ED-anal3}.
The identity \eq{Euler-entropy-finite++} expresses the {\it 2nd law of
thermodynamics}:\index{2nd law of thermodynamics}
in isolated systems, the total entropy is not decaying.

Using \eq{Eular-small-adiabatic-power}, the modification of
\eq{calculus-to-entropy++} for $\Ent=\ENT(\Ee,\theta)$ instead of
$\Eng=\ENG(\Ee,\theta)$ gives
\begin{align}
\!\frac1\theta\Big(\pdt\Ent+{\rm div}(\Ent\vv)\Big)
&=\pdt{}\eta(\Ee,\theta)+{\rm div}(\eta(\Ee,\theta)\vv)
+
\frac1{3\theta}{\rm tr}\big(\TT(\Ee,\theta){-}\TT(\Ee,0)\big)\,{\rm div}\,\vv\,.
\label{Euler-entropy-finite-u}\end{align}

The weak formulation of \eq{Euler-small-therm-ED-anal1} and
\eq{Euler-small-therm-ED-anal3}
uses the by-part integration in
time and the Green formula in $\varOmega$, in the former equation even twice
and still combined with a surface Green formula over $\varGamma$:

\begin{definition}[Weak formulation of
\eq{Euler-small-therm-ED-anal}]\label{def-thermo-Ch4}
The quadruple $(\varrho,\vv,\Ee,\theta)$ with
$\varrho\in L^\infty(I{\times}\varOmega)\cap W^{1,1}(I{\times}\varOmega)$,
$\vv\in L^p(I;W^{2,p}(\varOmega;\R^3))$,
$\Ee\in W^{1,1}(I{\times}\varOmega;\Rsym)$, and
$\theta\in L^1(I{\times}\varOmega)$ such that
${\rm tr}(\mathscr{T}(\Ee,\theta){-}\mathscr{T}(\Ee,0))\in L^{p'}(I;L^1(\varOmega))$, 
$\ENT(\Ee,\theta)\in L^{p'}(I;L^1(\varOmega))$, and
$\intkappa(\theta)\in L^1(I{\times}\varOmega)$ is called a weak solution to
the system \eq{Euler-small-therm-ED-anal} with the boundary
conditions \eq{Euler-small-therm-ED-BC-hyper} and the initial
conditions \eq{Euler-Ch4-thermodynam-IC} if
\begin{subequations}\label{def-thermo-Ch4-momentum}\begin{align}
&\nonumber
\!\!\!\int_0^T\!\!\!\!\int_\varOmega\bigg(\Big(\psi'_{\bm{E}}(\Ee,\theta)+
\bbD\strain(\vv)-\varrho\vv{\otimes}\vv\Big){:}\strain(\widetilde\vv)
+\psi(\Ee,\theta)\,{\rm div}\,\widetilde\vv
-\varrho\vv{\cdot}\pdt{\widetilde\vv}
\\[-.6em]&\hspace{3em}
+\HYPER|\nabla^2\vv|^{p-2}\nabla^2\vv\Vdots\nabla^2\widetilde\vv\bigg)\,\d\xx\d t
=\!\int_0^T\!\!\!\!\int_\varOmega\varrho\GRAVITY{\cdot}\widetilde\vv\,\d\xx\d t
+\!\int_\varOmega\!\varrho_0\vv_0{\cdot}\widetilde\vv(0)\,\d\xx
\label{def-thermo-Ch4-momentum1}
\intertext{holds for any $\widetilde\vv$ with $\widetilde\vv{\cdot}\nn={\bm0}$
on $I{\times}\varGamma$ and $\widetilde\vv(T)=0$, and}
\nonumber
&\!\!\!\!\int_0^T\!\!\!\!\int_\varOmega\!\bigg(\ENT(\Ee,\theta)\pdt{\wt\theta}+
\ENT(\Ee,\theta)\vv{\cdot}\nabla\wt\theta-\intkappa(\theta)\Delta\wt\theta
\\[-.6em]&\hspace{3em}\nonumber
+\Big(
\tfrac13{\rm tr}\big(\mathscr{T}(\Ee,\theta){-}\mathscr{T}(\Ee,0)\big)\,{\rm div}\,\vv
+\xi_{\rm ext}^{}(\theta;\vv,\ZJEp)\Big)\,\wt\theta\bigg)\,\d\xx\d t
\\[-.6em]&\hspace{7em}
+\!\!\int_0^T\!\!\!\!\int_\varGamma h(\theta)\wt\theta\,\d S\d t
+\!\!\int_\varOmega\!\!\ENT(\Ee_0,\theta_0)\,\wt\theta(0)\,\d\xx
=\!\int_0^T\!\!\!\!\int_\varGamma h_{\rm ext}\wt\theta\,\d S\d t
\label{def-thermo-Ch4-momentum3}\end{align}\end{subequations}
with $\intkappa(\theta):=\int_0^\theta\kappa(\vartheta)\,\d\vartheta$
and with $\xi_{\rm ext}^{}(\theta;\vv,\ZJEp)$ from
\eq{Euler-small-therm-ED-anal3} and $\ZJEp=\mathscr{R}(\Ee,\theta)$
holds for any $\wt\theta$ smooth with $\wt\theta(T)=0$ and
$\nn{\cdot}\nabla\wt\theta=0$ on $I{\times}\varGamma$, and if
\eq{Euler-small-therm-ED-anal0} and \eq{Euler-small-therm-ED-anal2} hold
a.e.\ on $I{\times}\varOmega$ together with the respective initial conditions
for $\varrho$ and $\Ee$ in \eq{Euler-Ch4-thermodynam-IC}.
\end{definition}

Notably, due to the general embedding $W^{1,1}(I;X)\subset C(I;X)$, the initial
conditions mentioned in the above definition for $\varrho$ and $\Ee$ have
a good sense in $L^1(\varOmega)$ and $L^1(\varOmega;\R_{\rm sym}^{3\times3})$,
respectively. In fact, it will be improved further in $W^{1,r}(\varOmega)$ and
$W^{1,s}(\varOmega;\R_{\rm sym}^{3\times3})$ later on, relying on the
estimates \eq{est-of-rho-disc} and \eq{Euler-small-E-Wr} below.

\begin{remark}[{\sl The condition $\vv\Cdot\nn=0$ in \eq{Euler-small-therm-ED-BC-hyper}}]\upshape
A common disadvantage of the Eulerian approach is that we have assumed the
shape of the considered bounded domain $\varOmega$ fixed by assuming zero
normal velocity on the boundary. This simplification is commonly used in
both theoretical and computational studies.  In applications where the
shape of the domain intended  to evolve, this disadvantage can be overcome by
embedding $\varOmega$ into an artificial soft medium whose boundary is fixed,
which is known as the sticky-air approach (in geophysics) or the
fictitious-domain approach or the immersed-boundary method (in engineering).
\end{remark}

\section{ Fully implicit  time discretization of \eq{Euler-small-therm-ED-anal}}
\label{sec-discrtete}

For simplicity, the viscosity coefficients $\bbD$ and $\HYPER$ in
\eq{Euler-small-therm-ED-anal} are considered constant here, although their
continuous dependence on $\theta$ and on $\Ee$ could easily be considered,
too. Similarly, a continuous dependence of $\kappa$ also on $\Ee$
 could  be considered for a stable scheme (i.e.\ for a-priori
estimates), although convergence  would  then  be 
analytically troublesome, cf.\ Remark~\ref{rem-Euler-thermo-kappa}.

\def\NOINDEX{}
\def\TAU{\tau}

We use the so-called {\it Rothe method},
i.e.\ the fully implicit time discretization with an equidistant partition of
the time interval $I$ with the time step $\tau>0$ such that $T/\tau$
 is an integer. 
We denote by $\varrho_\TAU^k$, $\vv_\TAU^k$, $\Ee_\TAU^k$, $\theta_\TAU^k$, ...
the approximate values of $\varrho$, $\vv$, $\Ee$, $\theta$, ...
at time instants $t=k\tau$ with $k=1,2,...,T/\tau$. Using $\mathscr{T}$ from
\eq{stress-entropy-Ch4} and
$\mathscr{R}$ from \eq{R-abbreviation}, we will then use the following
recursive regularized time-discrete scheme 
\begin{subequations}\label{Euler-small-therm-ED+discr}
\begin{align}\label{Euler-small-therm-ED+0discr}
&\!\!\frac{\varrho_\TAU^k{-}\varrho_\TAU^{k-1}}\tau=
-{\rm div}\,\pp_\TAU^k\ \ \ \ \text{ with }\ \
\pp_\TAU^k=\varrho_\TAU^k\vv_\TAU^k\,,\!\!
\\[-.2em]&\nonumber
\!\!\frac{\pp_\TAU^k{-}\pp_\TAU^{k-1}\!\!}\tau\,=
{\rm div}\big(
\mathscr{T}(\Ee_\TAU^k,\theta_\TAU^k){+}\DD_\TAU^k{-}
\pp_\TAU^k{\otimes}\vv_\TAU^k\big)+\varrho_\TAU^k\GRAVITY_{\tau}^k\,,\ 
\\[-.3em]&
\hspace*{6em}\text{ where }\ 
\DD_\TAU^k\!=\bbD\strain(\vv_\TAU^k)-{\rm div}\mathfrak{H}_\TAU^k\ 
\text{ with }\ \mathfrak{H}_\TAU^k\!=\HYPER\big|\nabla^2\vv_\TAU^k\big|^{p-2}
\nabla^2\vv_\TAU^k\,,\!\!
\label{Euler-small-therm-ED+1disc}
\\[-.5em]&\!\!\frac{\Ee_\TAU^k{-}\Ee_\TAU^{k-1}\!\!}\tau\,
=\strain(\vv_\TAU^k)-\ZJEp_\TAU^k
-\bm B_\text{\sc zj}^{}(\vv_\TAU^k,\Ee_\TAU^k)
\ \ \text{ with }\ \ZJEp_\TAU^k=\mathscr{R}(\Ee_\TAU^k,\theta_\TAU^k)
\,,\label{Euler-small-therm-ED+2disc}
\\[-.1em]&\nonumber\!\!\frac{\Ent_\TAU^k{-}\Ent_\TAU^{k-1}\!\!}\tau
={\rm div}\big(\kappa(\theta_\TAU^k)\nabla\theta_\TAU^k
+\Ent_\TAU^k\vv_\TAU^k\big)
+\xi_{\rm ext}^{}(\theta_\TAU^k;\vv_\TAU^k,\ZJEp_\TAU^k)
\\[-.2em]&\hspace*{5em}
+\tfrac13{\rm tr}\big(\mathscr{T}(\Ee_\TAU^k,\theta_\TAU^k)
-\mathscr{T}(\Ee_\TAU^k,0)\big)\,{\rm div}\,\vv_\TAU^k
\ \ \ \text{ with $\ \ \Ent_\TAU^k\!=\ENT(\Ee_\TAU^k,\theta_\TAU^k)$}\,,
\label{Euler-small-therm-ED+3discr}
\end{align}\end{subequations}
where $\xi_{\rm ext}^{}$ is from \eq{Euler-small-therm-ED-anal3} and
$\ENT(E,\theta)=\psi(E,\theta)-\theta\psi_\theta'(E,\theta)-\psi(E,0)$ is from
\eq{def-E-eta-U}.

The corresponding boundary conditions \eq{Euler-small-therm-ED-BC-hyper}
lead to
\begin{subequations}\label{BC-disc}\begin{align}\label{BC-disc-a}
&\!\vv_\TAU^k{\cdot}\nn=0\,,\ \ \,
\big[\big(\mathscr{T}(\Ee_\TAU^k,\theta_\TAU^k){+}\DD_\TAU^k\big)\nn
-\divS\big(\mathfrak{H}_\TAU^k\Cdot\nn\big)\big]_\text{\sc t}^{}\!\!=\bm0\,,\ \ \,
\\&
\nabla^2\vvk{:}(\nn{\otimes}\nn)=0\,,\ \ \ \text{ and }\ \ \
\kappa(\theta_\TAU^k)\nabla\theta_\TAU^k\Cdot\nn+h(\theta_\TAU^k)
=h_{{\rm ext},\tau}^k\,,
\label{BC-disc-b}\end{align}\end{subequations}
where $h_{{\rm ext},\tau}^k(\xx):=\int_{(k-1)\tau}^{k\tau}h_{\rm ext}(t,\xx)\,\d t$.
Similarly, in \eq{Euler-small-therm-ED+1disc}, we used
$\GRAVITY_{\tau}^k(\xx):=\int_{(k-1)\tau}^{k\tau}\GRAVITY(t,\xx)\,\d t$.

This system of boundary-value problems
\eq{Euler-small-therm-ED+discr}--\eq{BC-disc} for the quadruple
$(\varrho_\TAU^k,\pp_\TAU^k,\Ee_\TAU^k,\theta_\TAU^k)$ and thus also for
$\Ent_\TAU^k$ is to be solved recursively for $k=1,2,...,T/\tau$, starting for
$k=1$ with the initial conditions analogous to \eq{Euler-Ch4-thermodynam-IC},
i.e.\ 
\begin{align}\label{IC-disc}
\varrho_\TAU^0=\varrho_0\,,\ \ \ \ \ \ \ \ \pp_\TAU^0=\varrho_0\vv_{0}\,,\  
 \ \ \ \ \ \ {\Ee}_\TAU^0={\Ee}_0\,,\  
 \ \text{ and }\ \ \Ent_\TAU^0=\ENT({\Ee}_0,\theta_0)\,.
\end{align}
Thus, from \eq{Euler-small-therm-ED+0discr}, we obtain also
$\vvk=\pp_\TAU^k/\varrho_\TAU^k$  provided that  $\varrho_\TAU^k>0$,
as  will  indeed  be  proved in  Step~2 below 
in the next section~\ref{sec-proof}. In that section, we will also
prove the existence of a
solution to \eq{Euler-small-therm-ED+discr} with \eq{BC-disc} in Step~5. 

We will use the ``compact'' notation that exploits the interpolants.
Specifically, using the values $(\vvk)_{k=0}^{T/\tau}$, we define the piecewise
constant and the piecewise affine interpolants respectively as
\begin{align}
&\ovT(t)\!:=\vvk\ \ \text{ and }\ \ 
\vv_\TAU(t)\!:=\Big(\frac t\tau\,{-}\,k{+}1\Big)\vvk
\!+\Big(k\,{-}\,\frac t\tau\Big)\vvkk
\ \ \text{ for }\ (k{-}1)\tau<t\le k\tau
\label{def-of-interpolants}
\end{align}
for $k=0,1,...,T/\tau$. Analogously, we define also $\orT$,
$\varrho_\TAU^{}$, $\opT$, $\pp_\TAU^{}$, $\oEeT$,
$\Ee_\TAU^{}$, etc. Written ``compactly'' in terms of these interpolants, the
recursive system \eq{Euler-small-therm-ED+discr} writes as
\begin{subequations}\label{Euler-small-therm-ED+dis}
\begin{align}\label{Euler-small-therm-ED+0dis}
&\!\!\pdt{\varrho_\TAU}=-{\rm div}\,\opetau\ \ \
\text{ with $\ \ \opT=\orT\ovT$}\,,
\\&\nonumber
\!\!\pdt{\pp_\TAU}=
{\rm div}\big(\mathscr{T}(\oEeT,\otT){+}\oDT{-}\opT{\otimes}\ovT\Big)
+\orT\overline\GRAVITY_{\tau}\,,\ \ \ 
\\[-.1em]&
\hspace*{4.em}\text{ where }\ \oDetau\!=\bbD\strain(\ovT)
-{\rm div}\overline{\mathfrak{H}}_\TAU\,\ 
\text{ with }\ \overline{\mathfrak{H}}_\TAU\!=
\HYPER\big|\nabla^2\ovT\big|^{p-2}\nabla^2\ovT\,,
\label{Euler-small-therm-ED+1dis}
\\[-.3em]&
\!\!\pdt{\Ee_\TAU}=\strain(\ovT)-\oPT-\bm B_\text{\sc zj}^{}(\otT,\oEeT)
\ \ \text{ with }\ \oPT=\mathscr{R}(\oEeT,\otT)\,,
\label{Euler-small-therm-ED+2dis}
\\[.1em]&\nonumber\!\!\pdt{\Ent_\TAU}
={\rm div}\big(\kappa(\otT)\nabla\otT{-}\ouT\ovT\big)+
\xi_{\rm ext}^{}(\otT;\ovT,\oPT)\\[-.2em]&
\hspace*{3em}
+\tfrac13{\rm tr}\big(\mathscr{T}(\oEetau,\otT)
-\mathscr{T}(\oEeT,0)\big)\,{\rm div}\,\ovT
\ \ \ \text{ with $\ \ \ouT\!=\ENT(\oEeT,\otT)$}\,,
\label{Euler-small-therm-ED+3dis}
\end{align}\end{subequations}
where the dissipation rate $\xi_{\rm ext}^{}(\theta;\vv,\ZJEp)$ is from
\eq{Euler-small-therm-ED-anal3}. Of course, the boundary conditions
\eq{BC-disc} are ``translated'' accordingly, i.e.
\begin{subequations}\label{BC-disc+}\begin{align}
&\ovT\Cdot\nn=0,\ \ \ \
\big[(\mathscr{T}(\oEeT,\otT){+}\oDT)\nn
+\divS(\overline{\mathfrak{H}}_\TAU\nn)\big]_\text{\sc t}^{}=\bm0\,,\ \ \ \ 
\\&\overline{\mathfrak{H}}_\TAU\Colon(\nn{\otimes}\nn)={\bm0}\,,
\ \ \text{ and }\ \ \kappa(\otT)\nabla\otT\Cdot\nn+h(\otT)
 =\bar h_{{\rm ext},\tau}\ \text{ on }\ \varGamma\,,
\end{align}\end{subequations}
and with the corresponding initial conditions for
$(\varrho_\TAU,\pp_\TAU,\Ee_\TAU,\Ent_\TAU)$ arising for $t=0$ from \eq{IC-disc}.

For the mechanical part (\ref{Euler-small-therm-ED-anal}a-c) alone, such a
time discretization was used for the compressible (isothermal) Navier-Stokes
equations in \cite{Kar13CFEM,Zato12ASCN} or \cite[Ch.7]{FeKaPo16MTCV}.
The explicit use of the convexity of the kinetic energy expressed in terms of
the momentum $\pp$ is in \cite{FLMS21NACF}.

\section{Analysis of the time-discrete scheme
\eq{Euler-small-therm-ED+discr} in a special case}\label{sec-anal}

Proving the stability of the discrete scheme \eq{Euler-small-therm-ED+discr}
is quite technical and it is worth simplifying  by making the heat
capacity independent of the strain variable $\bm{E}$.

\subsection{A special ansatz of partly linearized thermomechanical coupling}
Without significant loss of applicability, we consider a commonly used
special form
\begin{align}\label{ansatz-for-entropy-base-estimate}
\psi(\bm{E},\theta)=\varphi(\bm{E})+\theta\,\phi({\rm tr}\,\bm{E})+\gamma(\theta)
\end{align}
with $\gamma(0)=0$, which will allow for a relatively simple estimation strategy
by exploting the entropy balance. It is physically natural to consider also
$\gamma'(0)=0$. The ansatz \eq{ansatz-for-entropy-base-estimate} implies
that ${\rm dev}(\psi''_{\bm{E}\theta}(\bm{E},\theta))
={\rm dev}(\phi'({\rm tr}\,\bm{E})\bbI)\equiv0$
so, in particular, it complies with \eq{Euler-small-ass}. Also, the adiabatic
power \eq{Eular-small-adiabatic-power} takes a specific form as
\begin{align}
\tfrac13
{\rm tr}\big(\mathscr{T}(\bm{E},\theta)\,{-}\mathscr{T}(\bm{E},0)\big)\,
{\rm div}\,\vv=\big(\theta\,\phi'({\rm tr}\,\bm{E})+\theta\,\phi({\rm tr}\,\bm{E})
+\gamma(\theta)\big)\,{\rm div}\,\vv\,.
\label{Eular-small-adiabatic-power+}\end{align}

The ansatz \eq{ansatz-for-entropy-base-estimate} has the effect of separating
the mechanical and the thermal variables additively in the internal energy
$\ENG$ and the entropy $\eta$ defined in \eq{def-E-eta-U}. Here it results in:
\begin{align}\label{Ch4-additive-split-energy}
\ENG(\bm{E},\theta)
=\varphi(\bm{E})+\!\!\!\!\!\lineunder{\gamma(\theta)
-\theta\gamma'(\theta)}{$\qquad=:\ENT(\theta)$}
\ \ \ \text{ and }\ \ \ \ 
\eta(\bm{E},\theta)
=-\,\phi({\rm tr}\,\bm{E})\!\!\!\!\!\!\lineunder{-\,\gamma'(\theta)}{$\qquad=:\eta_1(\theta)$}\!\!\!\!\!\!\,;
\end{align}
note that the thermal part of the internal energy
$\ENT(\theta)=\ENG(\bm{E},\theta)-\ENG(\bm{E},0)
=\ENG(\bm{E},\theta)-\varphi(\bm{E})$ here depends only on
temperature, i.e.\ $\ENT'_{\!\bm{E}}=0$. Likewise, the heat capacity
$\ENG_{\!\theta}'(\bm{E},\theta)=\ENT'(\theta)
=-\theta \gamma''(\theta)=:c(\theta)$ depends only on temperature. Also, the
thermal part of the entropy
$\eta_1(\theta)=\eta(\bm{E},\theta)-\eta(\bm{E},0)$
depends only on temperature. This simplifies many calculations but it slightly
corrupts the thermodynamical consistency (cf.~Remark~\ref{rem-3rd-law-partly}
below) no matter how often such models are used in the literature.

Given this special ansatz \eq{ansatz-for-entropy-base-estimate}, we
have the following calculus
\begin{align}\nonumber
\frac1{\theta}\Big(\pdt\Ent+{\rm div}(\Ent\vv)\Big)
&=\frac1{\theta}\Big(\pdt{}\ENT(\theta)+{\rm div}(\ENT(\theta)\vv)\Big)
\\\nonumber
&=\frac1{\theta}\pdt{}\big(\gamma(\theta){-}\theta\gamma'(\theta)\big)
+\frac1{\theta}{\rm div}\Big(\big(\gamma(\theta){-}\theta\gamma'(\theta)\big)\vv\Big)
\\
&=-\pdt{}\gamma'(\theta)-{\rm div}\big(\gamma'(\theta)\vv\big)
+\frac{\gamma(\theta)}{\theta}{\rm div}\,\vv\,.
\label{calculus-to-entropy-simpler}\end{align}
Reminding the thermal part of the entropy
$\eta_1=\eta_1(\theta)=-\gamma'(\theta)$ from \eq{Ch4-additive-split-energy}
and testing \eq{Euler-small-therm-ED-anal3} by the coldness $1/\theta$, we
obtain the balance of $\eta_1$ as 
\begin{align}\nonumber
\!\!\!
\pdt{}\eta_1^{}(\theta)=\frac{{\rm div}\big(\kappa(\theta)\nabla\theta\big)\!}\theta
-{\rm div}
\big(\eta_1^{}(\theta)\vv\big)+\frac{\xi_{\rm ext}^{}(\theta;\vv,\ZJEp
)\!}{\theta}
+\!\!\!\!\lineunder{\!\Big({\rm tr}\frac{\mathscr{T}(\Ee,\theta){-}\mathscr{T}(\Ee,0)}
{3\theta}{-}\frac{\gamma(\theta)}\theta\Big)\!}{$=\phi'({\rm tr}\,\Ee)+\phi({\rm tr}\,\Ee)$}\!\!\!\!
{\rm div}\,\vv\,
\\[-1.7em]\label{thermal-entropy-balance}\\[-2.5em]\nonumber\end{align}
with the heat power $\xi_{\rm ext}^{}=\xi_{\rm ext}^{}(\theta;\vv,\ZJEp)$ from
\eq{Euler-small-therm-ED-anal3}.

\begin{remark}[{\sl 3rd law of thermodynamics}]\label{rem-3rd-law-partly}
\upshape
A physically relevant requirement is that the entropy at zero temperature
is bounded from below and independent of the mechanical state. This is
called the 3rd law of thermodynamics. By default, the entropy at zero
temperature $\eta(\cdot,0)$ is thus
calibrated to zero. It should be openly pointed out that the ansatz
\eq{ansatz-for-entropy-base-estimate} satisfies this 3rd law of thermodynamic
only partially. Namely, the entropy
$\eta(\Ee,\theta)=\eta_1(\theta)-\phi({\rm tr}\,\Ee)$ is bounded from below
at zero absolute temperature when proving that the small-strain
field $\Ee$ and hence also $\phi$ is bounded, as we will indeed proved later;
recall that $\gamma'(0)$ is assumed to be bounded (even zero). However, if
$\phi(\cdot)$ is not constant, then $\eta(\cdot,0)$ is not constant either, and
thus  does not fully comply with the third law.
\end{remark}

\def\ALPHA{\alpha}
\def\ZETA{\gamma}
\def\TWO{2}
\def\TWOprime{2}
\def\EXP{\mu}

\def\othetatauEXP{\hspace*{.1em}\overline{\hspace*{-.1em}\theta}_{\etau}^{\,\EXP}}
\def\othetatauExp{\hspace*{.1em}\overline{\hspace*{-.1em}\theta}_{\etau}^{\,{\lambda}}}
\def\othetatauEXp{\hspace*{.1em}\overline{\hspace*{-.1em}\theta}_{\etau}^{\,_{1+\lambda}}}
\def\othetatauexp{\hspace*{.1em}\overline{\hspace*{-.1em}\theta}_{\etau}^{\,{1-\lambda}}}

\subsection{Stability and convergence of the time-discrete scheme
\eq{Euler-small-therm-ED+discr}}\label{sec-proof}

We will now perform the analysis for the
special ansatz \eq{ansatz-for-entropy-base-estimate}.

Recalling the notation $\ENG$ and $\ENT$ from \eq{def-E-eta-U} and
$\mathscr{R}$ from \eq{R-abbreviation}, let us summarize the assumptions used
in what follows. For suitable exponents $\alpha$ and $\beta$, which
will be specified later in \eq{Euler-small-ass-alpha-beta} and which refer,
respectively, to the poly\-nomial growth of the heat capacity and heat
conductivity, we assume 
\begin{subequations}\label{Euler-thermo-small-ass}
\begin{align}\label{Euler-thermo-small-ass-phi}
&\!\!\psi\in C^2(\Rsym{\times}\R)\,,\ 
\ \ \ \kappa\in C(\R)\,,\ \ \ h\in C(\R)\,,\ \ \ 
\phi(\cdot,0)\ \text{ convex}\,,
\\[-.1em]&\nonumber\hspace{-.2em}
\exists\, 0< c_0^{}\le C_1^{}<+\infty\ \
\ \ \forall(\bm{E},\theta)\in\Rsym{\times}\R{:}
\\[-.1em]&\nonumber
\hspace{5em}
\ENG(\bm{E},\theta)\ge c_0^{}(|\bm{E}|^2+(\theta^+)^{1+\ALPHA})\ \,\text{ and }\ \
\ENG_\theta'(\bm{E},\cdot)>0\ \text{ on $(0,+\infty)$\,,}
\\[-.1em]&
\hspace{5em}
|\ENG'_{\bm{E}}(\bm{E},\theta)|\le C_1^{}\big(1{+}|\bm{E}|
\big)\,\ \text{ and }\ \,\theta^+\ENG_\theta'(\bm{E},\theta)
\le C_1^{}\big(1{+}(\theta^+)^{1+\ALPHA}\big)\,,
\label{Euler-small-ass-therm-psi-1}
\\[-.1em]&\!\!\label{Euler-thermo-small-ass-entropy}
\max\big(|\eta(\bm{E},0)|\,,\,|\eta'_{\bm{E}}(\bm{E},0)|\big)\le
C_1^{}\big(1{+}|\bm{E}|^{C_1}\big)\,,
\\[-.1em]&\!\!\label{Euler-thermo-small-ass-Tth}
\big|\mathscr{T}(\bm{E},\theta)\big|
\le C_1\big(1+\ENG(\bm{E},\theta)\big)\,,
\\&\!\!\!\max\big(|\mathscr{R}(\bm{E},\theta)|\,,\,
|\mathscr{R}_\theta'(\bm{E},\theta)|\big)\le C_1^{}(1{+}|\bm{E}|)
\ \text{ and }\mathscr{R}'_{\bm{E}}\ \text{ positive semi-definite},
\label{Euler-thermo-small-ass-zeta+}
\\[-.1em]&\!\!\label{Euler-thermo-small-ass-D}
\bbD:\Rsym\to\Rsym\ \text{ linear symmetric positive semi-definite}\,,\ \
\ \ \HYPER>0\,,
\\[-.1em]&\!\!
\exists\,
\,\kappa_{\max}^{}\ge\kappa_{\min}^{}>0\ \ \forall\,\theta\in\R:\ \ \
\kappa_{\min}^{}(\theta^+)^{\beta^+}
\le\kappa(\theta)\le\kappa_{\max}^{}\big(1{+}(\theta^+)^\beta\big)\,,
\label{Euler-ass-therm-kappa}
\\&\nonumber
\zeta_{\rm p}\in C(\R{\times}\Rdev)\,,\ \
\theta\in\R:\ \ \zeta_{\rm p}(\theta,\cdot):\Rdev\to\R
\ \text{ convex},
\\&\hspace{4em}
0<\!\!\inf_{\theta\in\R^+,\ \Lp\in\Rdev\setminus\{0\}}\!\!\!\!\!
\frac{\zeta_{\rm p}(\theta,\Lp)}{|\Lp|^2}\ \ \ \ \text{ and }\ 
\sup_{\theta\in\R^+,\ \Lp\in \Rdev\setminus\{0\}}\!\!\!
\frac{\zeta_{\rm p}(\theta,\Lp)}{1{+}|\Lp|^2}\!<+\infty\,,
\label{Euler-thermo-small-ass-zeta}
\\[-.1em]&h:\R^+\!\!\to\R^+\text{ increasing},\ \ h(0)=0\,,\ \ 
h_{\rm ext}\,{\in}\,L^1(I{\times}\varGamma)\,,\ \
h_{\rm ext}\ge0\text{ a.e.\ on }I{\times}\varGamma,\!
\label{Euler-thermo-small-ass-BC}
\\[-.1em]&
\varrho_0\in W^{1,r}(\varOmega)\,,\ \ \ {\rm min}_{\barOmega}^{}\varrho_0>0\,,\ \ \
\vv_0\in L^2(\varOmega;\R^3)\,,\ \ \
\GRAVITY\in L^1(I;L^\infty(\varOmega;\R^3))\,,
\label{Euler-small-ass-IC}
\\[-.1em]&\Ee_0\in W^{1,\infty}(\varOmega;\Rsym)\,,\ \ \
\theta_0\in L^\alpha(\varOmega),\ \ \ \theta_0>0\ 
\text{ on }\ \varOmega\,,\ \ \ \ENG(\Ee_0,\theta_0)\in L^1(\varOmega)\,.
\label{Euler-thermo-small-ass-IC}
\end{align}\end{subequations}
An example for $h$ complying with \eq{Euler-thermo-small-ass-BC}, occurred in
\eq{Euler-small-therm-ED-BC-hyper}, is
$h(\theta)=\alpha_1\theta+\alpha_1\theta^4$ with $\alpha_1\ge0$, $\alpha_2\ge0$,
and $\alpha_1+\alpha_2>0$ should hold for $h$ to be increasing. Yet, in fact,
the an analysis in \eq{Euler-entropy-finite-disc} allows for $h\equiv0$, too.

\def\EXS{s}

\begin{proposition}[Stability of the discrete scheme and solutions to
\eq{Euler-small-therm-ED+discr}]\label{prop-thermo-Ch4-existence}
For the ansatz \eq{ansatz-for-entropy-base-estimate} leading, for some
$0<\lambda<2$, to the convex function
$1/\ENT(\cdot)^\lambda$, let the assumptions \eq{ass-T-E-commute},
(\ref{Euler-small-ass}a,b,d,e), and \eq{Euler-thermo-small-ass} hold with
$p> r>3$ and the exponents $\alpha$ and $\beta$ in
(\ref{Euler-thermo-small-ass}b,e) satisfy
\begin{align}\label{Euler-small-ass-alpha-beta}
1+\lambda>\beta^+\!\ge\frac23\alpha+\lambda-\frac13
\ \ \ \text{ and }\ \ \ \alpha\ge\Big(\frac32\lambda-1\Big)^+\,.
\end{align}
Moreover, let $h_{\rm ext}/h^{-1}(h_{\rm ext})^\lambda\in L^1(I{\times}\varGamma)$.
Then:\\
\Item{(i)}{The~time-discrete~scheme~\eq{Euler-small-therm-ED+dis}~possesses
a solution $(\varrho_\TAU,\vv_\TAU,\Ee_\TAU,\theta_\TAU)$ and is stable (in the
spaces specified later within the proof) with respect to $\tau>0$ provided
$\tau$ is sufficiently small.}
\Item{(ii)}{For $\tau\to0$, $(\varrho_\TAU,\vv_\TAU,\Ee_\TAU,\theta_\TAU)$ 
converges weakly* (in terms of subsequences) in the topologies specified in
\eq{Euler-small-converge}  and \eq{Euler-thermo-Ch4-conv2}  below and
every such a limit $(\varrho,\vv,\Ee,\theta)$ solves (in a weak sense  of
Definition~\ref{def-thermo-Ch4})
the initial-boundary-value problem in the sense.}
\Item{(iii)}{In particular, \eq{Euler-small-therm-ED-anal} has at least one
weak solution in the sense of Definition~\ref{def-thermo-Ch4} such that also
$\varrho\in L^\infty(I;W^{1,r}(\varOmega))\,\cap\,C(I{\times}\barOmega)$ with
$\min_{I{\times}\barOmega}\varrho>0$, $\vv\in L^\infty(I;L^2(\varOmega;\R^3))$,
$\Ee\in L^\infty(I;W^{1,\EXS}(\varOmega;\Rsym))$ with any $1\le\EXS<\EXP$, and
$\theta\in L^\infty(I;L^{1+\ALPHA}(\varOmega))\,\cap\,L^\EXP(I;W^{1,\EXP}(\varOmega))$
with}
\begin{align}\nonumber\\[-2.5em]
\label{Euler-entropy-nabla-theta-mu}
\ \ \EXP\le\frac{5+2\alpha+3\beta^+\!-3\lambda}{4+\alpha}\,.
\end{align}
\Item{(iv)}{M{}oreover, the energy-dissipation balance
\eq{Euler-small-energy-balance-stress} and the total-energy balance
\eq{Euler-engr-finite} integrated over the time interval $[0,t]$ holds for
any $t\in I$.}
\end{proposition}

\begin{proof}
We expand  and modify  the arguments from \cite[Sect.2]{Roub25TDVE},
where the isothermal variant with $\zeta_{\rm p}'=\zeta_{\rm p}'(\Lp)$ linear was
handled. Compared to \cite{Roub25TDVE}, the discrete scheme is 
 substantially  simpler here without any gradient regularization of the
equation \eq{Euler-small-therm-ED+2disc} for $\Ee$. Here we will need
A careful usage of the discrete Gronwall inequality, for which we refer e.g.\ to
\cite{QuaVal94NAPD,Roub13NPDE,Thom97GFEM}, is to be used here.
 Moreover,  unlike \cite{Roub25TDVE},  we avoid any regularization
of (\ref{Euler-small-therm-ED+discr}a,b).  For the heat part, which is
completely absent in \cite{Roub25TDVE}, we use the test of a
generalized-entropy type, similarly as in \cite{AbBuLe24ESGN} in the
incompressible case with $\alpha=0=\beta$.

For the sake of clarity, we divide the following argumentation into seven steps,
and  partly  use  some  calculations for the mechanical part
in \cite[Sect.2]{Roub25TDVE}.

\medskip{\it Step 1:  Basic stability of \eq{Euler-small-therm-ED+discr}
and first a-priori estimates  from the total energy}. 
Using the convexity of the kinetic energy as the functional
$(\varrho,\pp)\mapsto\int_\varOmega\frac12|\pp|^2/\varrho\,\d\xx$, cf.\
\cite[Sect.2.4]{Roub25TDVE}, and the assumed convexity
\eq{Euler-thermo-small-ass-phi} of the stored energy $\varphi$,
we now test \eq{Euler-small-therm-ED+1disc} by
$\vv_\TAU^k=\pp_\TAU^k/\varrho_\TAU^k$ while using \eq{Euler-small-therm-ED+0discr}
tested by $\frac12|\pp_\TAU^k|/(\varrho_\TAU^k)^2$ and
\eq{Euler-small-therm-ED+2disc} tested by $\varphi'(\Ee_\TAU^k)$, and derive
the discrete analog of \eq{Euler-small-energy-balance-stress} an inequality:
\begin{align}\nonumber
&\!\!\!\int_\varOmega\!\frac{|\pp_\TAU^k|^2\!\!}{2\varrho_\TAU^k\!}
  +\varphi(\Ee_\TAU^k)\,\d\xx
+\tau\!\sum_{m=1}^k\int_\varOmega\!\xi_{\rm ext}^{}(\theta_\TAU^m;\vv_\TAU^m,\ZJEp_\TAU^m)
  \,\d\xx
\le\int_\varOmega\!\frac{|\pp_\TAU^0|^2\!\!}{2\varrho_\TAU^0\!}+\varphi(\Ee_\TAU^0)\,\d\xx
\\[-.8em]&\qquad\qquad
+\tau\!\sum_{m=1}^k
\int_\varOmega\!\bigg(\varrho_\TAU^m\GRAVITY_{\EPS\DELTA}^m{\cdot}\vv_\TAU^m
-\Big(\theta_\TAU^m\phi'({\rm tr}\,\bm{E}_\TAU^m)+
\theta_\TAU^m\phi({\rm tr}\,\bm{E}_\TAU^m)+\gamma(\theta_\TAU^m)\Big)
\,{\rm div}\,\vv_\TAU^m\bigg)\,\d\xx\,.
\label{Euler-engr-finite-disc}\end{align}
In addition,  we  also perform the test of
\eq{Euler-small-therm-ED+3discr} by 1. This leads to the discrete analog of
the total-energy balance \eq{Euler-engr-finite}, namely 
\begin{align}\nonumber
\!\!\!\!\!\!\!\!&\!\!\!\!\int_\varOmega\!\frac{|\pp_\TAU^k|^2\!\!}{2\varrho_\TAU^k\!}
+\!\!\!\!\!\!\lineunder{\varphi(\Ee_\TAU^k)+
\ENT(\theta_\TAU^k)}{$\qquad\qquad=\ENG(\Ee_\TAU^k,\theta_\TAU^k)$}\!\!\!\!\d\xx
+\tau\!\sum_{m=1}^k\!\int_\varGamma\!h(\theta_\TAU^k)\,\d S
\\[-.8em]&\hspace{1.5em}
\le\int_\varOmega\!\frac{|\pp_0|^2\!\!}{2\varrho_0\!\!}
+\ENG(\Ee_0,\theta_0)\,\d\xx
+\tau\!\sum_{m=1}^k\bigg(\int_\varOmega\!\varrho_\TAU^m\GRAVITY_{\tau}^m\Cdot\vv_\TAU^m\,\d\xx+\!\!
\int_\varGamma\!h_{{\rm ext},\tau}^k\,\d S\bigg)
\,,\!\!
\label{Euler-thermo-basic-energy-balance-disc}\end{align}
where $\ENG(\bm{E},\theta):=\varphi(\bm{E})+\ENT(\theta)$ as in
\eq{Ch4-additive-split-energy}. To obtain some a-priori estimates, we handle
the gravity-loading term $\varrho_\TAU^m\GRAVITY_{\tau}^m$ by using the
continuity equation \eq{Euler-small-therm-ED-anal0} and the impenetrability
of the boundary, which guarantees 
\begin{align}\label{Euler-small-est-1}
\int_\varOmega\varrho_\TAU^k\,\d\xx=\!\int_\varOmega\varrho_0\,\d\xx=:M\ \ \
\text{ for any $k=1,...,T/\tau$}\,.
\end{align}

 Let us imagine, for a moment, that
\eq{Euler-small-therm-ED+0discr} is replaced by
$\varrho_\TAU^k{-}\varrho_\TAU^{k-1}=\tau{\rm div}((\varrho_\TAU^k)^+\vv_\TAU^k)$
with $(\cdot)^+$  denoting the non-negative part. By testing it by the
non-positive part $(\varrho_\TAU^k)^-$ of $\varrho_\TAU^k$ and realizing that
${\rm div}((\varrho_\TAU^k)^+\vv_\TAU^k)(\varrho_\TAU^k)^-=0$ a.e.,
we can see that $\varrho_\TAU^k\ge0$ if also $\varrho_\TAU^{k-1}\ge0$.
Therefore, $(\varrho_\TAU^k)^+=\varrho_\TAU^k$ so that, in fact, we arrive
at \eq{Euler-small-therm-ED+0discr}. For another argumentation more
directly for the original equation \eq{Euler-small-therm-ED+0discr},
we refer to \eq{calulus-nonlin-rho-r++} below.
Then, by exploiting that $\varrho_\TAU^k\ge0$, we can estimate 
\begin{align}\nonumber
\!\!\int_\varOmega&\varrho_\TAU^m\GRAVITY_\TAU^m\Cdot\vv_\TAU^m\,\d\xx
=\!\int_\varOmega\GRAVITY_\TAU^m\Cdot\pp_\TAU^m\,\d\xx
\le\big\|\!\sqrt{\varrho_\TAU^m}\GRAVITY_\TAU^m\big\|_{L^2(\varOmega;\R^3)}^{}
\Bigg\|\frac{\pp_\TAU^m}{\sqrt{\varrho_\TAU^m}}\Bigg\|_{L^2(\varOmega;\R^3)}^{}
\!\!\!
\\[-.4em]&\nonumber
\le\|\varrho_\TAU^m\|_{L^1(\varOmega)}^{1/2}\|\GRAVITY_\TAU^m\|_{L^\infty(\varOmega;\R^3)}^{}
\Bigg\|\frac{\pp_\TAU^m}{\sqrt{\varrho_\TAU^m}}\Bigg\|_{L^2(\varOmega;\R^3)}^{}\!\!\!
\le M^{1/2}\|\GRAVITY_\TAU^m\|_{L^\infty(\varOmega;\R^3)}^{}
\bigg(\!1+\frac14\bigg\|\frac{\pp_\TAU^m}{\sqrt{\varrho_\TAU^m}}\bigg\|_{L^2(\varOmega;\R^3)}^2\,\bigg)
\\[-.1em]&\hspace{18.2em}
=M^{1/2}\|\GRAVITY_\TAU^m\|_{L^\infty(\varOmega;\R^3)}^{}\bigg(\!1+\frac14\int_\varOmega\!\varrho_\TAU^m|\vv_\TAU^m|^2\,\d\xx\bigg).
\label{Euler-small-est-Gronwall-}\end{align}
 Furthermore, we need to know that $\theta_\TAU^k\ge0$, which can
be seen by testing \eq{Euler-small-therm-ED+3discr} with $(\theta_\TAU^k)^-$ and
by relying on that $\theta_\TAU^{k-1}\ge0$ so that also $\ENT(\theta_\TAU^{k-1})\ge0$.
Then, 
relying on $L^\infty(I{\times}\varOmega;\R^3)$-bound of $\GRAVITY$, we can use
the discrete Gronwall inequality to obtain the a-priori estimates
\begin{subequations}\label{Euler-small-est}\begin{align}
\label{Euler-small-est1}
&\bigg\|\frac{\opT}{\sqrt{\orT}}\bigg\|_{L^\infty(I;L^2(\varOmega;\R^3))}^{}\le C\,,
\intertext{and a similar $L^\infty(I;L^1(\varOmega))$-estimate of
$\ENG(\oEeT,\otT)$.  By using 
the coercivity assumption in \eq{Euler-small-ass-therm-psi-1}, we also obtain}
\label{Euler-small-therm-ED-est1}
&\|\oEeT\|_{L^\infty(I;L^2(\varOmega;\R^{3\times3}))}^{}\le C
\ \ \ \ \text{ and }\ \ \ \ \|\otT\|_{L^\infty(I;L^{1+\alpha}(\varOmega))}^{}\le C\,.
\intertext{Thus, by exploiting \eq{Euler-thermo-small-ass-Tth}, we obtain 
an estimate for $\oTT\!=\mathscr{T}(\oEeT,\otT)$:}
&\|\oTT\|_{L^\infty(I;L^1(\varOmega;\R^{3\times3}))}^{}\le C\,.
\label{Euler-small-therm-ED-est1+}
\intertext{ Moreover, from the first estimate in \eq{Euler-small-est1},
we have also }
& \|\opT\|_{L^\infty(I;L^1(\varOmega;\R^3))}^{}\le\!\!\!\!
\lineunder{\big\|\!\sqrt{\orT}\big\|_{L^\infty(I;L^2(\varOmega))}^{}}
{$=M^{1/2}$}\!\!\!
\bigg\|\frac{\opT}{\sqrt{\orT}}\bigg\|_{L^\infty(I;L^2(\varOmega;\R^3))}^{}\!\!
\le C\,.
\label{Euler-small-est3+++}
\end{align}\end{subequations}
It should also be emphasized that all the estimates \eq{Euler-small-est} hold
only for sufficiently small $\tau>0$ for which the discrete Gronwall inequality
can be applied. Specifically they hold for all
$\tau\le 1/(M^{1/2}\|\GRAVITY\|_{L^\infty(I{\times}\varOmega;\R^3)})$.

\medskip{\it Step 2: Further estimates from the mechanical energy dissipation}.
We will now exploit \eq{Euler-engr-finite-disc} itself. Due to
\eq{Euler-thermo-small-ass-Tth} and the already obtained bounds, we have that
$\mathscr{T}(\Ee_\TAU^m,\theta_\TAU^m)$ is bounded in $L^1(\varOmega;\R^{3\times3})$
uniformly in $m$. Moreover, the growth of $\ENG(\cdot,0)$ in
\eq{Euler-small-ass-therm-psi-1} with \eq{Euler-small-therm-ED-est1} ensures
also that $\mathscr{T}(\Ee_\TAU^m,0)$ is bounded in $L^1(\varOmega;\Rsym)$.
Altogether, 
$\mathcal{T}_\TAU^m:=
{\rm tr}(\mathscr{T}(\Ee_\TAU^m,\theta_\TAU^m){-}\mathscr{T}(\Ee_\TAU^m,0))$
is bounded in $L^1(\varOmega;\Rsym)$ uniformly in $m$.
Thus we can estimate the power of the adiabatic effects as
\begin{align}\nonumber
\!\!-\!\int_\varOmega&\mathcal{T}_\TAU^m\,{\rm div}\,\vv_\TAU^m\,\d\xx
\le
\|\mathcal{T}_\TAU^m\|_{L^1(\varOmega)}^{}\|{\rm div}\,\vv_\TAU^m\|_{L^\infty(\varOmega)}^{}
\\[-.5em]&\nonumber\
\le C\|\mathcal{T}_\TAU^m\|_{L^1(\varOmega;\R^{3\times3})}^{}
\Big(\|\strain(\vv_\TAU^m)\|_{L^2(\varOmega;\R^{3\times3})}^{}
+\|\nabla^2\vv_\TAU^m\|_{L^p(\varOmega;\R^{3\times3\times3})}^{}\Big)
\\[-.1em]&\
\le C_{\delta}^{}\max_{n=1,...,T/\tau}\|\mathcal{T}_\TAU^n\|_{L^1(\varOmega;\R^{3\times3})}^{}\!
+\delta\|\strain(\vv_\TAU^m)\|_{L^2(\varOmega;\R^{3\times3})}^2\!
+\delta\|\nabla^2\vv_\TAU^m\|_{L^p(\varOmega;\R^{3\times3\times3})}^p
\label{Euler-small-est2-}\end{align}
with $\delta>0$ sufficiently small with respect to $\inf_{|E|=1}\bbD E\Colon E>0$
and to $\HYPER>0$. The power of the gravity in \eq{Euler-engr-finite-disc} can
be treated by emploing the information about the total mass as in
\eq{Euler-small-est-Gronwall-}. Thus, from \eq{Euler-engr-finite-disc}, we
obtain the estimate for the dissipation rate. From this, we can read the
estimates
\begin{align}
&\label{Euler-small-est2}
\|\strain(\ovT)\|_{L^2(I{\times}\varOmega;\R^{3\times3})}^{}\le C
\ \ \ \,\text{ and }\ \ \ \,
\|\nabla^2\ovT\|_{L^p(I{\times}\varOmega;\R^{3\times3\times3})}^{}\le C\,.
\end{align}

Then \eq{Euler-small-therm-ED+0discr} can be tested by 
$|\varrho_\TAU^k|^{s-2}\varrho_\TAU^k=(\varrho_\TAU^k)^{s-1}$ with some $s>1$.
Using the Green formula with the boundary condition $\nn\Cdot\vv_\TAU^k=0$,
the convective term can be handled as
\begin{align}\nonumber
\!\!\int_\varOmega{\rm div}(\varrho_\TAU^k\vv_\TAU^k)&
(|\varrho_\TAU^k|^{s-2}\varrho_\TAU^k)\,\d\xx
=\!\!\int_\varOmega({\rm div}\,\vv_\TAU^k)|\varrho_\TAU^k|^s
+(\vv_\TAU^k\Cdot\nabla\varrho_\TAU^k)
(|\varrho_\TAU^k|^{s-2}\varrho_\TAU^k)\,\d\xx
\\[-.3em]&
=\!\int_\varOmega\!({\rm div}\,\vv_\TAU^k)|\varrho_\TAU^k|^s
-\varrho_\TAU^k\,{\rm div}\big(|\varrho_\TAU^k|^{s-2}\varrho_\TAU^k\vv_\TAU^k\big)\,\d\xx
=\Big(1{-}\frac1s\Big)\int_\varOmega({\rm div}\,\vv_\TAU^k)|\varrho_\TAU^k|^s\d\xx\,,
\label{transport-rho}\end{align}
so that this test gives the inequality
\begin{align}
  &\int_\varOmega\!\frac{|\varrho_\TAU^k|^s\!-|\varrho_\tau^{k-1}|^s\!}\tau
  \,\d\xx\,\le\!\frac1{s'}\int_\varOmega\!({\rm div}\,\vv_\tau^k)|\varrho_\TAU^k|^s
\,\d\xx
\le\frac1{s'}\|{\rm div}\,\vv_\tau^k\|_{L^\infty(\varOmega)}^{}\int_\varOmega\!|\varrho_\TAU^k|^s\,\d\xx\,.
\label{transport-rho+}\end{align}
Due to \eq{Euler-small-est1} with $p>3$, we have
$\tau\sum_{k=1}^{T/\tau}\|{\rm div}\,\vv_\tau^k\|_{L^\infty(\varOmega)}^p\le C$,
from which we can see that
$\max_{1\le k\le T/\tau}\|{\rm div}\,\vv_\tau^k\|_{L^\infty(\varOmega)}^{}\le(C/\tau)^{1/p}$.
This allows us to use the discrete Gronwall inequality which, for a sufficiently
small $\tau$, say for $\tau\le \frac12C^{1/(p+1)}$, gives the estimate
\begin{align}\label{Euler-small-est-rho-s}
\|\orT\|_{L^\infty(I;L^s(\varOmega;\R^3))}^{}\le C\,.
\end{align}
Furthermore, we test \eq{Euler-small-therm-ED+0discr} by 
${\rm div}(|\nabla\varrho_\TAU^k|^{r-2}\nabla\varrho_\TAU^k)$ with some $r>1$.
Using the Green formula with the boundary condition
$\nn\Cdot\vv_\TAU^k=0$, the convective term can be handled as
\begin{align}\nonumber
&\!\!\int_\varOmega\!\nabla\big(\vv_\TAU^k\Cdot\nabla\varrho_\TAU^k
  \big)\Cdot\big(|\nabla\varrho_\TAU^k|^{r-2}\nabla\varrho_\TAU^k\big)\,\d\xx
\\[-.5em]&\hspace{.1em}\nonumber
=\!\int_\varOmega|\nabla\varrho_\TAU^k|^{r-2}(\nabla\varrho_\TAU^k{\otimes}\nabla\varrho_\TAU^k)\Colon\strain(\vv_\TAU^k)
+(\vv_\TAU^k\Cdot\nabla)\nabla\varrho_\TAU^k
\Cdot\big(|\nabla\varrho_\TAU^k|^{r-2}\nabla\varrho_\TAU^k\big)\,\d\xx
\\[-.1em]&\hspace{.1em}\nonumber
=\!\int_\varGamma|\nabla\varrho_\TAU^k|^r\,\vv_\TAU^k\Cdot\nn\,\d S
+\!\int_\varOmega\bigg(|\nabla\varrho_\TAU^k|^{r-2}(\nabla\varrho_\TAU^k{\otimes}\nabla\varrho_\TAU^k)\Colon\strain(\vv_\TAU^k)
\\[-.9em]&\hspace{12em}\nonumber
-({\rm div}\,\vv_\TAU^k)|\nabla\varrho_\TAU^k|^r-(r{-}1)|\nabla\varrho_\TAU^k|^{r-2}|\nabla\varrho_\TAU^k|^{r-2}\nabla\varrho_\TAU^k\Cdot
(\vv_\TAU^k\Cdot\nabla)\nabla\varrho_\TAU^k\bigg)\,\d\xx
\\[-.5em]&\hspace{.1em}
=\!\int_\varGamma\frac1r|\nabla\varrho_\TAU^k|^r\!\!\!\!
\lineunder{\vv_\TAU^k\Cdot\nn}{$=0$}\!\!\!\!\d S
+\!\int_\varOmega|\nabla\varrho_\TAU^k|^{r-2}(\nabla\varrho_\TAU^k{\otimes}\nabla\varrho_\TAU^k)\Colon\strain(\vv_\TAU^k)-\frac1r({\rm div}\,\vv_\TAU^k)|\nabla\varrho_\TAU^k|^r\,\d\xx\,.\!\!
\label{calulus-nonlin-rho-r}\end{align}
Thus, counting also with $\nabla(\varrho_\TAU^k\,{\rm div}\,\vv_\TAU^k)\Cdot
(|\nabla\varrho_\TAU^k|^{r-2}\nabla\varrho_\TAU^k)=
({\rm div}\,\vv_\TAU^k)|\nabla\varrho_\TAU^k|^r
+\varrho_\TAU^k|\nabla\varrho_\TAU^k|^{r-2}\nabla\varrho_\TAU^k
\Cdot\nabla({\rm div}\,\vv_\TAU^k)$,
we obtain the inequality
\begin{align}\nonumber
&\frac1r
\int_\varOmega\frac{|\nabla\varrho_\TAU^k|^r{-}|\nabla\varrho_\TAU^{k-1}|^r
\!\!}\tau\,\d\xx
\le\!\int_\varOmega\bigg(\Big(\frac1r{-}1\Big)
({\rm div}\,\vv_\TAU^k)|\nabla\varrho_\TAU^k|^r\!
-|\nabla\varrho_\TAU^k|^{r-2}(\nabla\varrho_\TAU^k{\otimes}
\nabla\varrho_\TAU^k)\Colon\strain(\vv_\TAU^k)
\\[-.7em]&\nonumber\hspace{24em}
-\varrho_\TAU^k|\nabla\varrho_\TAU^k|^{r-2}\nabla\varrho_\TAU^k
\Cdot\nabla({\rm div}\,\vv_\TAU^k)\bigg)\,\d\xx
\\[-.5em]&
\le C\big\|\strain(\vv_\TAU^k)\big\|_{L^\infty(\Omega;\R^{3\times3})}^{}
\int_\Omega|\nabla\varrho_\TAU^k|^r\,\d\xx
+\|\varrho_\tau^k\|_{L^{pr/(p-r)}}^{}\|\nabla\varrho_\tau^k\|_{L^{r}(\varOmega;\R^3)}^{r-1}
\|\nabla({\rm div}\,\vv_\tau^k)\|_{L^p(\Omega;\R^3)}^{}
\,.\!\!
\label{calulus-nonlin-rho-r+}\end{align}
Here, we needed $p>r$ to be able then to choose $s=pr/(p{-}r)<+\infty$ for
\eq{Euler-small-est-rho-s}. Using the discrete Gronwall inequality again,
relying on \eq{Euler-small-est1}, we obtain the estimate for sufficiently
small time steps $\tau$ for the mass density:
\begin{align}\label{est-of-rho-disc}
\|\orT\|_{L^\infty(I;W^{1,r}(\varOmega))}^{}\le C\,.
\end{align}

Furthermore, by testing \eq{Euler-small-therm-ED+0discr} with the nonpositive
part $(\varrho_\TAU^k)^-=\min(0,\varrho_\TAU^k)\in W^{1,r}(\varOmega)$
of $\varrho_\TAU^k$, we obtain
\begin{align}\nonumber
\int_\varOmega((\varrho_\TAU^k)^-)^2\,\d\xx
=\int_\varOmega\!\!\!\!\!\lineunder{\varrho_\TAU^{k-1}(\varrho_\TAU^k)^-}{$\le0$ if
$\varrho_\TAU^{k-1}\ge0$}\!\!\!\!\!\!
-&\tau{\rm div}(\varrho_\TAU^k\vv_\TAU^k)(\varrho_\TAU^k)^-\,\d\xx
\le-\tau\!\int_\varOmega
\varrho_\TAU^k\vv_\TAU^k\Cdot\nabla(\varrho_\TAU^k)^-\,\d\xx
\\[-1.1em]&=-\tau\!
\int_\varOmega(\varrho_\TAU^k)^-\vv_\TAU^k\Cdot\nabla(\varrho_\TAU^k)^-\,\d\xx
=\frac\tau2\!\int_\varOmega((\varrho_\TAU^k)^-)^2\,{\rm div}\,\vv_\TAU^k\,\d\xx\,,
\label{calulus-nonlin-rho-r++}\end{align}
where we also used the calculus like \eq{transport-rho} for $s=2$.
From this, for each $\tau>0$ sufficiently small, namely for
$\tau<2/\|{\rm div}\,\vv_\TAU^k\|_{L^\infty(\varOmega)}$,
we can see that $\varrho_\TAU^k\ge0$ provided $\varrho_\TAU^{k-1}\ge0$
on $\Omega$. Even, we can prove positivity of $\varrho_\TAU^k$. This can be
seen via a contradiction argument: Assuming that $\varrho_\TAU^{k-1}>0$
on $\barOmega$ and that the minimum of a (momentarily smooth) solution
$\varrho_\TAU^k$ is attained at some $\xx\in\varOmega$ and
$\varrho_\TAU^k(\xx)=0$, we have that
$\varrho_\TAU^k(\xx)-\varrho_\TAU^{k-1}(\xx)<0$ and
$\nabla\varrho_\TAU^k(\xx)=\bm0$. This yields that
$\varrho_\TAU^k(\xx)-\varrho_\TAU^{k-1}(\xx)
=-\tau\vv_\TAU^k(\xx)\Cdot\nabla\varrho_\TAU^k(\xx)-\tau\varrho_\TAU^k(\xx){\rm div}\vv_\TAU^k(\xx)<0$ but, in view of
\eq{Euler-small-therm-ED+0discr} written at $\xx$, it should be equal to 0. Thus
we obtain the contradiction, showing that $\varrho_\TAU^k(\xx)>0$.

This a.e.\ positivity allows for a test with $\sigma_\TAU^k:=1/\varrho_\TAU^k$.
Using the convexity of the function $\varrho\mapsto1/\varrho$ on $(0,+\infty)$
for the convective time differences and \eq{Euler-small-therm-ED+0discr}, we
obtain the inequality a.e.\ on $\varOmega$:
\begin{align}\nonumber
&\frac{\sigma_\TAU^k-\sigma_{\tau}^{k-1}\!\!\!}\tau
+\vv_\TAU^k\Cdot\nabla\sigma_\TAU^k
\ =\frac1\tau\Big(\frac1{\varrho_\TAU^k}-\frac1{\varrho_{\tau}^{k-1}}\Big)
+\vv_\TAU^k\Cdot\nabla\frac1{\varrho_\TAU^k}
\\&\qquad
\le-
\frac1{(\varrho_\TAU^k)^2\!}\Big(\frac{\varrho_\TAU^k-\varrho_{\tau}^{k-1}\!\!\!}\tau
+\vv_\TAU^k\Cdot\nabla\varrho_\TAU^k\Big)
=\frac1{\varrho_\TAU^k\!}{\rm div}\,\vv_\TAU^k
=({\rm div}\,\vv_\TAU^k)\,\sigma_\TAU^k\,,
\label{ineq-for-sparsity}\end{align}
where $\sigma_{\tau}^{k-1}:=1/\varrho_\tau^{k-1}$. Testing  \eq{ineq-for-sparsity}
by $|\sigma_\TAU^k|^{s-2}\sigma_\TAU^k$ for some $s>1$ and using the slightly
modified calculus \eq{transport-rho}, this test gives the following inequality
\begin{align}\nonumber
\frac1s\int_\varOmega\frac{(\sigma_\TAU^k)^s
-(\sigma_{\tau}^{k-1})^s}\tau\,\d\xx
&\le\int_\varOmega\Big(1{+}\frac1s\Big)({\rm div}\,\vv_\TAU^k)
(\sigma_\TAU^k)^s\,\d\xx
\\&
\le\Big(1{+}\frac1s\Big)\|{\rm div}\,\vv_\TAU^k\|_{L^\infty(\varOmega)}^{}
\int_\varOmega(\sigma_\TAU^k)^{s}\,\d\xx\,.
\label{est-of-sparsity}\end{align}
For a sufficiently small $\tau>0$, say for
$\tau\le1/((1{+}s)\max_{1\le k\le T/\tau}
\|{\rm div}\,\vv_\TAU^k\|_{L^\infty(\varOmega)}^{})$, we obtain a uniform
bound for $\|\sigma_\TAU^k\|_{L^s(\varOmega)}^{}$ by using the discrete Gronwall
inequality. Here note that we have $\max_{1\le k\le T/\tau}
\|{\rm div}\,\vv_\TAU^k\|_{L^\infty(\varOmega)}^{}\le\sqrt{NC/\tau}$ with
$C$ from the second estimate in 
\eq{Euler-small-est2} and with $N$ the norm of the embedding
$L^2(I;L^2(\varOmega))\,\cap\, L^p(I;W^{1,p}((\varOmega))\subset
L^2(I;L^\infty((\varOmega))$. Thus, choosing $\tau\le 1/((1{+}s)^2NC)$ is sufficient.
Having $\|\sigma_\TAU^k\|_{L^s(\varOmega)}^{}$ bounded and using
\eq{Euler-small-est1} while realizing that
$\vv_\TAU^k=\sqrt{\sigma_\TAU^k}\,\pp_\TAU^k/\sqrt{\varrho_\TAU^k}$,
we obtain the bound 
\begin{align}
\|\vv_\TAU^k\|_{L^{2s/(s+1)}(\varOmega;\R^3)}^{}
\le C_s\big\|\sqrt{\sigma_\TAU^k}\big\|_{L^{2s}(\varOmega)}^{}
\Big\|\frac{\pp_\TAU^k}{\sqrt{\varrho_\TAU^k}}\Big\|_{L^2(\varOmega;\R^3)}\,.
\label{Euler-large-est-v-2-eps}\end{align}
From \eq{Euler-large-est-v-2-eps} with an arbitrarily large $s>1$, we have obtained
an estimate
\begin{align}
\|\ovT\|_{L^\infty(I;L^a(\varOmega;\R^3))}^{}\le C_a\ \ \text{ with
an arbitrary }\ 1\le a<2\,.
\label{Euler-large-est-v-2-eps+}\end{align}
Together with \eq{Euler-small-est1}, we thus obtain the bound for
$\ovT$ in $L^a(I;W^{2,p}(\varOmega;\R^3))$ with $1\le a<2$. In particular,
\begin{align}
\|\ovT\|_{L^{p}(I;W^{1,\infty}(\varOmega;\R^3))}^{}\le C\,.
\label{Euler-large-est-v-2-eps++}\end{align}
Moreover, \eq{Euler-large-est-v-2-eps+} together with
\eq{Euler-small-est1} and \eq{est-of-rho-disc} also allow us to augment 
\eq{Euler-small-est3+++} by  an estimate for $\nabla\opT=
\nabla(\orT\ovT)=\orT\nabla\ovT+\nabla\orT{\otimes}\ovT$, namely
\begin{align}
\|\opT\|_{L^{p}(I;W^{1,r}(\varOmega;\R^3))}^{}\le C\,.
\label{Euler-large-est-p+}\end{align}

By the comparison from \eq{Euler-small-therm-ED+0dis}, using
\eq{Euler-large-est-p+}, we obtain
\begin{subequations}\label{Euler-small-est++}\begin{align}
&
\Big\|\pdt{\varrho_\TAU}\Big\|_{ L^p(I;L^r(\varOmega))}\!\le C\,. 
\label{Euler-small-est6}
\intertext{Moreover, by the comparison from
\eq{Euler-small-therm-ED+1disc} when using  \eq{Euler-small-therm-ED-est1},
\eq{est-of-rho-disc}, \eq{Euler-large-est-v-2-eps++}, and
\eq{Euler-large-est-p+},  we obtain} 
\nonumber
&\Big\|\pdt{\pp_\TAU}\Big\|_{L^{p'}(I;W^{2,p}(\varOmega;\R^3)^*)}^{}\!
=\!\!\sup_{\|\wt\vv\|_{L^p(I;W^{2,p}(\varOmega;\R^3))}\le1}\int_0^T\!\!\!\int_\varOmega
\Big(\oTT{+}\bbD\strain(\ovT){-}
\opT{\otimes}\ovT\Big)\Colon\strain(\wt\vv)\!\!\!
\\[-.8em]&\hspace{22em}
+\overline{\mathfrak{H}}_\TAU\Vdots\nabla^2\wt\vv
-\orT\overline\GRAVITY_{\tau}\Cdot\wt\vv
\,\d\xx\d t\le C\,;\label{Euler-small-est6+}
\end{align}\end{subequations}
note that, as $p\ge3$ is assumed,
$
\orT{\otimes}\ovT\in L^{p/2}(I;L^\infty(\varOmega;\Rsym))$ is indeed in duality
with $\strain(\wt\vv)\in L^p(I;L^\infty(\varOmega;\Rsym))$.

\medskip{\it Step 3: A test towards a generalized entropy}.
The next test is for \eq{Euler-small-therm-ED+3discr} by
$(1/\theta_\TAU^k)^\lambda$ for some $0<\lambda$. Note that, for $\lambda=1$,
such test leads to the standard entropy balance. To this aim,
 reminding that $\theta_\TAU^k\ge0$, 
we can also show that $\theta_\TAU^k$ is positive
a.e.\ on $\varOmega$. More specifically, let us assume that $\theta_\TAU^{k-1}>0$
on $\barOmega$ and that the minimum of a (momentarily smooth) solution
$\theta_\tau^k$ is attained at some $\xx\in\varOmega$ and $\theta_\TAU^k(\xx)=0$.
Thus $\theta_\TAU^k(\xx)-\theta_\tau^{k-1}(\xx)<0$ and
$\nabla\theta_\tau^k(\xx)=\bf0$ and also
${\rm div}(\kappa(\theta_\tau^k(\xx))\nabla\theta_\tau^k(\xx))\ge0$. The heat
equation \eq{Euler-small-therm-ED+3discr} written as $\Ent_\TAU^k{-}\Ent_\TAU^{k-1}
=\tau({\rm div}\,\vv_\TAU^k)\Ent_\TAU^k+\tau\vv_\TAU^k\Cdot\nabla\Ent_\TAU^k
+\tau{\rm div}(\kappa(\theta_\TAU^k)\nabla\theta_\TAU^k)
+\frac13\tau
{\rm tr}(\mathscr{T}(\Ee_\TAU^k,\theta_\TAU^k){-}\mathscr{T}(\Ee_\TAU^k,0))
\,{\rm div}\,\vv_\TAU^k+\tau\xi_{\rm ext}^{}(\theta_\TAU^k;\vv_\TAU^k,\ZJEp_\TAU^k)\ge0$
at this point $\xx$. Since $\ENT$ is increasing on $[0,+\infty)$,
 we have  $\Ent_\TAU^k(\xx){-}\Ent_\TAU^{k-1}(\xx)<0$
and we thus obtain the contradiction, showing  that  $\theta_\TAU^k(\xx)>0$.

We now use the assumption that the function $1/[\ENT^{-1}](\cdot)^\lambda$ is convex,
so that we have the inequality
\begin{align}
\frac1{(\theta_\TAU^k)^\lambda}(\Ent_\TAU^k\!-\Ent_\TAU^{k-1})=
\frac1{[\ENT^{-1}](\Ent_\TAU^k)^\lambda}(\Ent_\TAU^k\!-\Ent_\TAU^{k-1})
\le\wt\eta_\lambda^{}(\Ent_\TAU^k)\,-\,\wt\eta_\lambda^{}(\Ent_\TAU^{k-1})\,,
\label{calculus-to-entropy-disc-}\end{align}
where $\wt\eta_\lambda^{}$ is a primitive function of $1/[\ENT^{-1}](\cdot)^\lambda$.
Note that, up to a possible additive constant which can be calibrated as 0,
it holds that $\wt\eta_\lambda^{}\circ\ENT=\eta_\lambda$, where
\begin{align}\label{calculus-to-entropy-lambda}
\eta_\lambda(\theta)
:=-\int_0^\theta\!\!\!\vartheta^{1-\lambda}\gamma''(\vartheta)\,\d\vartheta\,;
\end{align}
note that we have $[\wt\eta_\lambda^{}{\circ}\ENT]'(\theta)
=((1/\ENT^{-1})^\lambda\circ\ENT(\theta))\ENT'(\theta)
=(1/\theta^\lambda)\ENT'(\theta)=
-\theta^{1-\lambda}\gamma''(\theta)=\eta_\lambda'(\theta)$.
Further, we will exploit the discrete version of the calculus
\begin{align}\nonumber
\frac1{\theta^\lambda}\Big(\pdt\Ent+{\rm div}(\Ent\vv)\Big)
&=\frac1{\theta^\lambda}\Big(\pdt{}\ENT(\theta)+{\rm div}(\ENT(\theta)\vv)\Big)
\\\nonumber
&=\frac1{\theta^\lambda}\pdt{}\big(\gamma(\theta){-}\theta\gamma'(\theta)\big)
+\frac1{\theta^\lambda}{\rm div}\Big(\big(\gamma(\theta){-}\theta\gamma'(\theta)\big)\vv\Big)
\\
&=\pdt{}\eta_\lambda(\theta)+{\rm div}\big(\eta_\lambda(\theta)\vv\big)
+\bigg(\frac{\gamma(\theta){-}\theta\gamma'(\theta)}{\theta^\lambda}-\eta_\lambda(\theta)\bigg)\,{\rm div}\,\vv\,;
\label{calculus-to-entropy-simpler-modified}\end{align}
note that, for $\lambda=1$, it reduces to the calculus
\eq{calculus-to-entropy-simpler} for $\eta_1$ as defined in
\eq{Ch4-additive-split-energy}. Thus, using
\eq{calculus-to-entropy-disc-}, we have now the inequality
\begin{align}\nonumber
\frac1{(\theta_\TAU^k)^\lambda}\bigg(\frac{\Ent_\TAU^k\!{-}\Ent_\TAU^{k-1}}\tau
+{\rm div}(\Ent_\TAU^k\vv_\TAU^k)\bigg)
&\le\frac{\eta_\lambda(\theta_\TAU^k)\,{-}\,\eta_\lambda(\theta_\TAU^{k-1})}\tau
\\[-.7em]&\hspace{0em}
+{\rm div}\Big(\eta_\lambda(\theta_\TAU^k)\vv_\TAU^k\Big)+
\bigg(\frac{\gamma(\theta_\TAU^k){-}\theta_\TAU^k\gamma'(\theta_\TAU^k)}{(\theta_\TAU^k)^\lambda}
-\eta_\lambda(\theta_\TAU^k)\bigg)\,{\rm div}\,\vv_\TAU^k\,.
\label{calculus-to-entropy-simpler-disc}\end{align}
Multiplying \eq{Euler-small-therm-ED+3discr} by $(\theta_\TAU^k)^{-\lambda}$ and 
integrating it over $\varOmega$ while using the Green formula with the boundary
condition for the normal velocity, we obtain
\begin{align}\nonumber
\!\!\int_\varOmega\!\frac{\xi_{\rm ext}^{}(\theta_\TAU^k;\vv_\TAU^k,\ZJEp_\TAU^k)}
{(\theta_\TAU^k)^\lambda}
&+\kappa(\theta_\TAU^k)\frac{|\nabla\theta_\TAU^k|^2\!\!\!}{(\theta_\TAU^k)^{1+\lambda}}\,\d\xx
+\!\int_\varGamma\frac{\kappa(\theta_\TAU^k)\nabla\theta_\TAU^k\!}{(\theta_\TAU^k)^\lambda}
\Cdot\nn\,\d S\
\\[-.3em]
&=\!\int_\varOmega\!
\frac1{(\theta_\TAU^k)^\lambda}\bigg(\frac{\Ent_\TAU^k\!{-}\Ent_\TAU^{k-1}\!\!}\tau
+{\rm div}(\Ent_\TAU^k\vv_\TAU^k)
-\tfrac13\mathcal{T}_\TAU^m{\rm div}\,\vv_\TAU^m\bigg)
\,\d\xx\,.
\label{calculus-to-entropy-simpler-disc+}\end{align}
The handling of the boundary term in \eq{calculus-to-entropy-simpler-disc+} is quite
tricky. Since $h(\cdot)$ is strictly monotone, cf.\ \eq{Euler-thermo-small-ass-BC},
we can identify an external temperature
\begin{align}\nonumber
\theta_{{\rm ext},\tau}^k:=h^{-1}(h_{{\rm ext},\tau}^k)\,.
\end{align}
Together with the non-negativity of $h$, this yields the following estimate
\begin{align}\nonumber
&\frac{\kappa(\theta_\TAU^k)\nabla\theta_\TAU^k\!}{(\theta_\TAU^k)^\lambda}\Cdot\nn
=\frac{h_{{\rm ext},\tau}^k{-}h(\theta_\TAU^k)\!}{(\theta_\TAU^k)^\lambda}
=\frac{h(\theta_{{\rm ext},\tau}^k){-}h(\theta_\TAU^k)\!}{(\theta_\TAU^k)^\lambda}
\pm\frac{h(\theta_\TAU^k){-}h(\theta_{{\rm ext},\tau}^k)\!}{(\theta_{{\rm ext},\tau}^k)^\lambda}
\\&\qquad=\frac{h(\theta_\TAU^k)-h(\theta_{{\rm ext},\tau}^k)\!}
{(\theta_{{\rm ext},\tau}^k)^\lambda}
+\frac{\big(h(\theta_\TAU^k)-h(\theta_{{\rm ext},\tau}^k)\big)\big((\theta_\TAU^k)^\lambda
{-}(\theta_{{\rm ext},\tau}^k)^\lambda\big)\!}
{(\theta_\TAU^k)^\lambda(\theta_{{\rm ext},\tau}^k)^\lambda}
\ge-\frac{h(\theta_{{\rm ext},\tau}^k)\!}{(\theta_{{\rm ext},\tau}^k)^\lambda}\,.
\nonumber\end{align}
By merging it with \eq{calculus-to-entropy-simpler-disc}, realizing
\eq{Eular-small-adiabatic-power+}, and summing it over the time levels
$k=1,...,T/\tau$, we obtain the inequality
\begin{align}&\nonumber
  \!\sum_{k=1}^{T/\tau}\int_\varOmega\!\frac{\xi_{\rm ext}^{}(\theta_\TAU^k;\vv_\TAU^k,\ZJEp_\TAU^k)\!}{(\theta_\TAU^k)^\lambda}
+\kappa(\theta_\TAU^k)\frac{|\nabla\theta_\TAU^k|^2\!\!\!}{(\theta_\TAU^k)^{1+\lambda}}\,\d\xx
\le\!\!\int_\varOmega\!\frac{\eta_\lambda(\theta_\TAU^{T/\tau})-\eta_\lambda(\theta_0)\!}\tau\,\d\xx
+\!\sum_{k=1}^{T/\tau}\!\int_\varGamma\frac{h(\theta_{{\rm ext},\tau}^k)\!}
{(\theta_{{\rm ext},\tau}^k)^\lambda}\,\d S
\\[-.4em]&\hspace{6em}
-\!\sum_{k=1}^{T/\tau}\int_\varOmega\!\bigg(
\theta_\TAU^k\frac{\phi'({\rm tr}\Ee_\TAU^k){+}\phi({\rm tr}\Ee_\TAU^k)\!}
{(\theta_\TAU^k)^{\lambda}}
+\frac{\theta_\TAU^k\gamma'(\theta_\TAU^k){+}\gamma(\theta_\TAU^k)\!\!}
{(\theta_\TAU^k)^\lambda}
+\eta_\lambda(\theta_\TAU^k)\bigg)\,{\rm div}\,\vv_\TAU^k\,\d\xx
\label{Euler-entropy-finite-disc-}\end{align}
with $\xi_{\rm ext}^{}=\xi_{\rm ext}^{}(\theta;\vv,\ZJEp)$ from
\eq{Euler-small-therm-ED-anal3}. This inequality written in terms
of the interpolants yields the variant of the thermal entropy balance
\eq{thermal-entropy-balance} integrated over the time integral $I$ and over
$\varOmega$ as an inequality, i.e.
\begin{align}&\nonumber
\!\!\!\int_0^T\!\!\!\!\int_\varOmega\!
\frac{\xi_{\rm ext}^{}(\otT;\ovT,\oPT)}{\otT^\lambda}
+\kappa(\otT)\frac{|\nabla\otT|^2\!\!\!}{\ \otT^{1+\lambda}}\ \d\xx\d t
\le\!\int_\varOmega\!\!
\eta_\lambda\big(\theta_\TAU(T)\big)-\eta_\lambda(\theta_0)\,\d\xx  
+\!\int_0^T\!\!\!\!\int_\varGamma\!\frac{h(\bar\theta_{{\rm ext},\tau}^{})}
{\bar\theta_{{\rm ext},\tau}^\lambda}\,\d S\d t
\\[-.5em]&\hspace{.5em}
-\!\int_0^T\!\!\!\!\int_\varOmega\!
\otT^{1-\lambda}\big(\phi'({\rm tr}\oEeT){+}\phi({\rm tr}\oEeT)\big)
{\rm div}\,\ovT+\Big(\otT^{1-\lambda}\gamma'(\otT)+\otT^\lambda\gamma(\otT)
+\eta_\lambda(\otT)\Big)\,{\rm div}\,\ovT
\,\d\xx\d t\,.\!\!\!
\label{Euler-entropy-finite-disc}\end{align}

The last integral in \eq{Euler-entropy-finite-disc} is bounded due to
\eq{Euler-small-est2} and \eq{Euler-small-therm-ED-est1}. In more detail,
we first improve the former estimate in \eq{Euler-small-therm-ED-est1} by
testing \eq{Euler-small-therm-ED+2disc} with $|\Ee_\TAU^k|^{r-2}\Ee_\TAU^k$,
where $r$ is arbitrarily large. Reminding the abbreviation $\mathscr{R}$
from \eq{R-abbreviation} and the assumption
\eq{Euler-thermo-small-ass-zeta+}, we can estimate 
\begin{align}\nonumber
\frac1r\int_\varOmega\frac{|\Ee_\TAU^k|^r\!-|\Ee_\TAU^{k-1}|^r\!\!}{\tau}\,\d\xx
&\le\!\!\int_\varOmega\!|\Ee_\TAU^k|^{r-2}\Ee_\TAU^k\Colon\Big(\strain(\vv_\TAU^k)
-\bm B_\text{\sc zj}^{}(\vv_\TAU^k,\Ee_\TAU^k)
-\mathscr{R}(\Ee_\TAU^k,\theta_\TAU^k)\Big)\,\d\xx
\\[-.1em]
&
\le\Big(1+
\|\nabla\vv_\TAU^k\|_{L^\infty(\varOmega;\R^{3\times3})}^{}\!+\|\Ee_\TAU^k\|_{L^r(\varOmega;\R^{3\times3})}^{}\Big)
\|\Ee_\TAU^k\|_{L^r(\varOmega;\R^{3\times3})}^{r-1}\,;
\label{transport-}\end{align}
here, for the term $(\vv_\TAU^k{\cdot}\nabla\Ee_\TAU^k)$ involved in
$\bm B_\text{\sc zj}^{}(\vv_\TAU^k,\Ee_\TAU^k)$, we have used the Green formula with
the boundary condition $\vv_\TAU^k{\cdot}\nn=0$ on $\varGamma$ as
\begin{align}\nonumber
&\!\!\!\int_\varOmega(\vv_\TAU^k{\cdot}\nabla\Ee_\TAU^k)\Colon\big(|\Ee_\TAU^k|^{ r-2}\Ee_\TAU^k\big)\,\d\xx
=\int_\varGamma\!|\Ee_\TAU^k|^r(\vv_\TAU^k{\cdot}\nn)\,\d S
\\[-.4em]&\hspace{1em}
-\!\int_\varOmega\!
(r{-}1)|\Ee_\TAU^k|^{r-2}\Ee_\TAU^k\Colon\big(\vv_\TAU^k{\cdot}\nabla\Ee_\TAU^k\big)
+({\rm div}\,\vv_\TAU^k)|\Ee_\TAU^k|^r\d\xx
=-\!\int_\varOmega\!\!\frac{{\rm div}\,\vv_\TAU^k\!}r\,|\Ee_\TAU^k|^r\d\xx.\!
\label{transport}\end{align}
By the assumed regularity \eq{Euler-thermo-small-ass-IC} of the initial
condition $\Ee_0$ with the growth assumption \eq{Euler-thermo-small-ass-zeta+}
on $\mathscr{R}$, by the discrete Gronwall inequality, for sufficiently
small $\tau>0$, we obtain
\begin{align}\label{Euler-small-E-Lr}
&\|\oEeT\|_{L^\infty(I;L^r(\varOmega;\R^{3\times3}))}^{}\le C_r\,.
\end{align}
Here, we have used that
$\tau\sum_{k=1}^{T/\tau}\|\nabla\vv_\TAU^k\|_{L^\infty(\varOmega;\R^{3\times3})}^{2}$
is bounded  due to \eq{Euler-small-est2} with Korn's inequality, so that
$\tau\max_{k=1,...,T/\tau}\|\nabla\vv_\TAU^k\|_{L^\infty(\varOmega;\R^{3\times3})}=\mathscr{O}(\sqrt{\tau})$
is sufficiently small for $\tau$ small enough, which allows us to use of the
mentioned  discrete Gronwall inequality.

This allows estimation of the terms on the right-hand side of
\eq{Euler-entropy-finite-disc}. Using \eq{Euler-small-E-Lr} with $r$
sufficiently large and the growth assumption
\eq{Euler-thermo-small-ass-entropy} when realizing that
$\eta_E'(E,0)=\phi'({\rm tr}E)\bbI$ due to the ansatz
\eq{ansatz-for-entropy-base-estimate}, the term
$\otT^{1-\lambda}\phi'({\rm tr}\oEeT)\,{\rm div}\,\ovT$ can be estimated as
\begin{align}\nonumber
&-\int_\varOmega\!\otT^{1-\lambda}\phi'({\rm tr}\oEeT)\,{\rm div}\,\ovT\,\d\xx
\le\|\otT\|_{L^{(1+\alpha)/(1-\lambda)}(\varOmega)}^{1-\lambda}
\|\phi'({\rm tr}\oEeT)\|_{L^{r/C_1}(\varOmega)}^{}
\|\strain(\ovT)\|_{L^\infty(\varOmega;\R^{3\times3})}^{}\,.
\nonumber\end{align}
The term $\otT^{1-\lambda}\phi({\rm tr}\oEeT)\,{\rm div}\,\ovT$ can be
estimated similarly. To estimate the term
$\otT^{1-\lambda}\gamma'(\otT)\,{\rm div}\,\ovT$, we use the growth
$|\gamma'(\theta)|=|\int_0^\theta\gamma''(\vartheta)\,\d\vartheta|=\mathscr{O}(\theta^\ALPHA)$. Here note the growth
$|\gamma''(\theta)|=|\ENT'(\theta)/\theta|=|\ENG_{\!\theta}'(E,\theta)/\theta|=\mathscr{O}(\theta^{\ALPHA-1})$
due to the assumption \eq{Euler-small-ass-therm-psi-1}.
Finally, in view of \eq{calculus-to-entropy-lambda}, when using the ansatz
\eq{ansatz-for-entropy-base-estimate} and the estimate
\eq{Euler-small-therm-ED-est1}, we can see that $\eta_\lambda(\otT)$ with
$\eta_\lambda(\theta)=\int_0^\theta\ENT'(\vartheta)/\vartheta^\lambda\,\d\vartheta$
is bounded in $L^\infty(I;L^{1/(1-\lambda)}(\varOmega))$.
Therefore, $\eta_\lambda(\theta_\TAU(T))$ is bounded in $L^{1/(1-\lambda)}(\varOmega)$ and
$\eta_\lambda(\otT)\,{\rm div}\,\ovT$ is bounded in
$L^p(I;L^{1/(1-\lambda)}(\varOmega))$. Eventually, the boundary term
$h(\bar\theta_{{\rm ext},\tau}^{})/\bar\theta_{{\rm ext},\tau}^\lambda
=\bar h_{{\rm ext},\tau}^{}/h^{-1}(\bar h_{{\rm ext},\tau}^{})^\lambda$ in
\eq{Euler-entropy-finite-disc} is uniformly bounded in
$L^1(I{\times}\varGamma)$ under the assumption
$h_{\rm ext}^{}/h^{-1}(h_{\rm ext}^{})^\lambda\in L^1(I{\times}\varGamma)$.

With the left-hand side of \eq{Euler-entropy-finite-disc} estimated, we
therefore obtain a bound for $\kappa(
\otT)|\nabla\otT|^2/\otT^{1+\lambda}$ in
$L^1(I{\times}\varOmega)$, which can be used to estimate
$\nabla\othetatau$ in $L^\EXP(I{\times}\varOmega;\R^3)$ for some 
$1\le \EXP<2$. To this end, we apply H\"older's inequality to 
$|\nabla\theta_\EPS|^\EXP$ written as the product of
${\otT}^{\hspace*{-.3em}\EXP(1+\lambda)/2}\!/\kappa(\otT)^{\EXP/2}$ and 
$\kappa(\otT)^{\EXP/2}|\nabla\otT|^\EXP/{\otT}^{\hspace*{-.3em}\EXP(1+\lambda)/2}$:
\begin{align}
&\!\!\!\!\int_0^T\!\!\!\!\int_\varOmega|\nabla\otT|^\EXP\,\d\xx\d t
\le C_{\EXP,\lambda}\Bigg(\int_0^T\!\!\!\!\int_\varOmega
\frac{
\overline{\hspace*{-.1em}\theta}_\TAU^{\:\EXP(1+\lambda)/(2-\EXP)}}
{\kappa(\otT)^{\EXP/(2-\EXP)}\!\!}\,\d\xx\,\d t\Bigg)^{\!\!1-\EXP/\TWO}
\Bigg(\int_0^T\!\!\!\!\int_\varOmega
\frac{\kappa(\otT)|\nabla\otT|^2\!\!\!}{\!\!\otT^{1+\lambda}}
\,\d\xx\,\d t\Bigg)^{\EXP/\TWO}\!\!\!\!\!\!
\label{Euler-entropy-nabla-theta}
\end{align}
with a constant $C_{\EXP,\lambda}$ dependent on $\EXP$ and $\lambda$. The last
integral in \eq{Euler-entropy-nabla-theta} is estimated from
\eq{Euler-entropy-finite-disc}, while the penultimate integral is to be
estimated using the latter estimate in \eq{Euler-small-therm-ED-est1}.
The first integral on the right-hand side is to be estimated by interpolation with
the latter estimate in \eq{Euler-small-therm-ED-est1}, specifically
\begin{align}\nonumber
\int_0^T\!\!\!\int_\varOmega\frac{\otT^{\:\EXP(1+\lambda)/(2-\EXP)}}
{\kappa(\otT)^{\EXP/(2-\EXP)}\!\!}\,\d\xx\,\d t
=\int_0^T\!\!\!\int_\varOmega\!\bigg(\frac{\otT^{\:1+\lambda}}
{\kappa(\otT)\!\!}\,\bigg)^{\EXP/(2-\EXP)}\!\!\!\d\xx\,\d t
&\!\!\stackrel{\eq{Euler-ass-therm-kappa}}{\le}\!\!
C\!\!\int_0^T\!\!\!\int_\varOmega\!\!\Big(1{+}\,
\otT^{\;\EXP(1+\lambda-\beta^+)/(2-\EXP)}\Big)\,\d\xx\,\d t
\\[-.3em]&\ \le C'
\bigg(1+\int_0^T\!\!\!\!\int_\varOmega|\nabla\otT|^\EXP\,\d\xx\d t\bigg)\,.
\label{Euler-entropy-nabla-theta+}\end{align}
The last inequality relies on the latter estimate in
\eq{Euler-small-therm-ED-est1} via the Gagliardo-Nirenberg inequality.
Specifically, for each time instant $t$ (not explicitly denoted), we have
\begin{align}\nonumber
&\|1{+}\otT\|_{L^{\EXP(1+\lambda-\beta^+)/(2-\EXP)}(\varOmega)}^{}
\le C\|1{+}\otT\|_{L^{1+\alpha}(\varOmega)}^{}
\Big(\|1{+}\otT\|_{L^{1+\alpha}(\varOmega)}^{}\!+
\|\nabla\otT\|_{L^{\EXP}(\varOmega;\R^3)}^{}\Big)
\\
&\hspace{12em}
\text{ for }\ \ \frac{2{-}\EXP}{\EXP(
1{+}\lambda{-}\beta^+)}
\ge a\Big(\frac1\EXP-\frac13\Big)+\frac{1{-}a}{1{+}\ALPHA}
\ \text{ with }\ 0<a\le1\,.
\label{8-cond}
\end{align}
We raise it to the power $\EXP(1{+}\lambda{-}\beta^+)/(2{-}\EXP)$ and choose
$a$ to obtain the desired exponent $a\EXP(1{+}\lambda{-}\beta^+)/(2{-}\EXP)=\EXP$,
i.e.\ $a=(2{-}\EXP)/(1{+}\lambda{-}\beta^+)$, which allows for the 
last inequality in \eq{Euler-entropy-nabla-theta+}. Since the exponent
$1-\EXP/2$ is less than 1, merging \eq{Euler-entropy-nabla-theta} and
\eq{Euler-entropy-nabla-theta+} yields the bound
\begin{align}\label{Euler-entropy-nabla-theta++}
\|\nabla\otT\|_{L^\EXP(I{\times}\varOmega;\R^3)}^{}\le C\,.
\end{align}
After some algebra, substituting the mentioned choice
$a=(2{-}\EXP)/(1{+}\lambda{-}\beta^+)$ into \eq{8-cond} yields the bound
\eq{Euler-entropy-nabla-theta-mu}. In more detail, note that, for this choice
of $a$, the inequality in \eq{8-cond} simply becomes $a/3\ge(1{-}a)/(1{+}\ALPHA)$,
i.e. $a\ge3/(4{+}\ALPHA)$. Realizing $a=(2{-}\EXP)/(1{+}\lambda{-}\beta^+)$,
we eventually obtain \eq{Euler-entropy-nabla-theta-mu}.

According to the Sobolev embedding theorem , \eq{Euler-entropy-nabla-theta++} yields
a bound for $\otT$ in $L^\EXP(I;L^{\EXP^*}(\varOmega))$ and, by interpolation with
the bound in $L^\infty(I;L^{1+\alpha}(\varOmega))$ with the weights
$3/(4{+}\alpha)$ and $(1{+}\alpha)/(4{+}\alpha)$, also the bound
\begin{align}\label{Euler-entropy-nabla-theta+++}
\|\otT\|_{L^{(4{+}\alpha)\,\EXP/3}(I{\times}\varOmega)}^{}\le C\,.
\end{align}
From this, realizing the bounds $c(\theta)=\mathscr{O}(\theta^\alpha)$
and $\ENT(\theta)=\mathscr{O}(\theta^{1+\alpha})$ and 
that $\nabla\ouT=\nabla\ENT(\otT)
=c(\otT)\nabla\otT$, we can read also the estimates
\begin{align}\label{Euler-entropy-nabla-theta++++}
\|\ouT\|_{L^{(4+\alpha)\,\EXP/(3+3\alpha)}(I{\times}\varOmega)}^{}\le C
\ \ \text{ and }\ \ 
\|\nabla\ouT\|_{L^{(4+\alpha)\,\EXP/(4+4\alpha)}(I{\times}\varOmega;\R^3)}^{}\le C\,.
\end{align}
Recalling $\intkappa$ as a primitive function of $\kappa$ as used in
\eq{def-thermo-Ch4-momentum3}, we realize
that $\intkappa(\theta)=\mathscr{O}(\theta^{1+\beta})$ and thus, from
\eq{Euler-entropy-nabla-theta+++}, we obtain the estimate
\begin{align}\label{Euler-hat-kappa-est}
\|\intkappa(\otT)\|_{L^{(4{+}\alpha)\,\EXP/(3+3\beta)}(I{\times}\varOmega)}^{}\le C
\end{align}

Let us specify the feasible exponents $\alpha$ and $\beta$ depending on
$\lambda$. First, we assume that $\alpha\ge0$ in order to work with the
conventional $L^{1+\alpha}$-Lebesgue space.
Recalling \eq{Euler-entropy-nabla-theta-mu}, the above-specified bounds
$1\le\EXP<2$ needed for \eq{Euler-entropy-nabla-theta} lead respectively
to the restrictions $1+\alpha+3\beta^+\ge3\lambda$ and $\beta^+<1+\lambda$.
The latter restriction is involved in \eq{Euler-small-ass-alpha-beta}
while the former restriction is automatically satisfied
if $\alpha\ge0$ and if the exponents in \eq{Euler-entropy-nabla-theta++++}
are greater than (or equal to) 1. This last requirement means the
restriction $1+3\beta^+\ge 2\alpha+3\lambda$, which is also included in
\eq{Euler-small-ass-alpha-beta}. Furthermore, also the exponent
in \eq{Euler-entropy-nabla-theta+++} should be greater than (or
equal to) 1, which means another restriction $2+2\alpha+3\beta^+\ge2\lambda$.
This is automatically satisfied by the previous restrictions.
Eventually, also the exponent in \eq{Euler-hat-kappa-est} should be greater
than (or equal to) 1, which needs 
$2+2\alpha\ge3\lambda+3(\beta{-}\beta^+)$. This restricts $\alpha\ge0$
further as $\alpha\ge\frac32\lambda-1$ when $\lambda>\frac23$.
Altogether, we obtain the conditions in \eq{Euler-small-ass-alpha-beta}.
 Cf.\ also Remark~\ref{rem-alpha-beta} below. 

\medskip{\it Step 4: An estimate of $\nabla\Ee_\TAU$}.
For a mere Kelvin-Voigt rheology, i.e., when the Maxwellian viscosity
vanishes by considering  $\ZJEp\equiv0$, the estimate of
$\nabla\Ee_\TAU$ would follow from the obtained regularity of the
velocity field \eq{Euler-small-est2} together with the assumed regularity
of the initial condition $\Ee_0$. For the general temperature-dependent
Jeffreys rheology, the estimation of $\nabla\Ee_\TAU$ is more involved
and we also rely on the previously obtained estimate for $\nabla\theta_\TAU$.

To obtain an $L^\EXS$-estimate $\nabla\Ee_\TAU$ with some $\EXS>1$,
we can test \eq{Euler-small-therm-ED+2disc} by
${\rm div}(|\nabla\Ee_\TAU^k|^{\EXS-2}\nabla\Ee_\TAU^k)$ and choose $r$ in
\eq{Euler-small-E-Lr} sufficiently large so that $1/r+1/\EXP+1/\EXS'\le1$.
Therefore, we need $1/r+1/\EXP\le1/\EXS$, which requires that $\EXS<\EXP$.
For this, we can choose $r\ge\EXS\EXP/(\EXP{-}\EXS)$.
Additionally, we also need $1/p+1/\EXS'\le1$, which requires that $\EXS\le p$,
a condition that is always met. This yields
\begin{align}&\nonumber
\frac1\EXS\!\int_\varOmega\frac{|\nabla\Ee_\TAU^k|^\EXS-|\nabla\Ee_\TAU^{k-1}|^\EXS\!\!}{\tau}\,\d\xx
\le\!\int_\varOmega|\nabla\Ee_\TAU^k|^{\EXS-2}\nabla\Ee_\TAU^k\Vdots\,
\Big(\nabla\strain(\vv_\TAU^k)-\nabla\mathscr{R}(\Ee_\TAU^k,\theta_\TAU^k)
-\nabla\bm B_\text{\sc zj}^{}(\vv_\TAU^k,\Ee_\TAU^k)\Big)\,\d\xx
\\[-.3em]&\nonumber\ \ \ \le
\int_\varOmega|\nabla\Ee_\TAU^k|^{\EXS-2}\nabla\Ee_\TAU^k\Vdots\,
\Big(\nabla\strain(\vv_\TAU^k)
-\mathscr{R}_{\Ee}'(\Ee_\TAU^k,\theta_\TAU^k){\cdot}\nabla\Ee_\TAU^k\!
-\mathscr{R}_\theta'(\Ee_\TAU^k,\theta_\TAU^k){\otimes}\nabla\theta_\TAU^k\!
-\nabla\bm B_\text{\sc zj}^{}(\vv_\TAU^k,\Ee_\TAU^k)\Big)\,\d\xx
\\&
\ \ \ \le C\Big(1{+}\|\nabla\vv_\TAU^k\|_{W^{1,p}(\varOmega;\R^{3\times3})}^{}
{+}\|\Ee_\TAU^k\|_{L^r(\varOmega;\R^{3\times3})}^{}
\|\nabla\theta_\TAU^k\|_{L^\EXP(\varOmega;\R^3)}^{}\Big)
\|\nabla\Ee_\TAU^k\|_{L^\EXS(\varOmega;\R^{3\times3\times3})}^{\EXS-1}
\label{nabla-E-calulus}\end{align}
with $C$ depending here on $(p,r,\EXS,\EXP)$. Here we have also used the
assumption \eq{Euler-thermo-small-ass-zeta+}, so that,
in particular, we have $|\nabla\Ee_\TAU^k|^{\EXS-2}\nabla\Ee_\TAU^k\Vdots
(\mathscr{R}_{\Ee}'(\Ee_\TAU^k,\theta_\TAU^k){\cdot}\nabla\Ee_\TAU^k)\ge0$, 
and for the term $\nabla((\vv_\TAU^k\Cdot\nabla)\Ee_\TAU^k)$ contained in
$\nabla\bm B_\text{\sc zj}^{}(\vv_\TAU^k,\Ee_\TAU^k)$, we have used the Green formula:
\begin{align}\nonumber
&\!\!\int_\varOmega\!\nabla\big((\vv_\TAU^k\Cdot\nabla)\Ee_\TAU^k
  \big)\Vdots|\nabla\Ee_\TAU^k|^{\EXS-2}\nabla\Ee_\TAU^k\,\d\xx
\\[-.5em]&\hspace{.1em}\nonumber
=\!\int_\varOmega\!|\nabla\Ee_\TAU^k|^{\EXS-2}(\nabla\Ee_\TAU^k{\boxtimes}\nabla\Ee_\TAU^k)\Colon\strain(\vv_\TAU^k)
+(\vv_\TAU^k\Cdot\nabla)\nabla\Ee_\TAU^k\Vdots|\nabla\Ee_\TAU^k|^{\EXS-2}\nabla\Ee_\TAU^k\,\d\xx
\\[-.1em]&\hspace{.1em}\nonumber
=\!\int_\varGamma|\nabla\Ee_\TAU^k|^\EXS\!\!\!\lineunder{\!\vv_\TAU^k\Cdot\nn\!}{$=0$}\!\!\!\d S
+\int_\varOmega\Big(|\nabla\Ee_\TAU^k|^{\EXS-2}(\nabla\Ee_\TAU^k{\boxtimes}\nabla\Ee_\TAU^k)\Colon\strain(\vv_\TAU^k)
\\[-1.1em]&\hspace{12em}\nonumber
-({\rm div}\,\vv_\TAU^k)|\nabla\Ee_\TAU^k|^\EXS-(\EXS{-}1)|\nabla\Ee_\TAU^k|^{\EXS-2}\nabla\Ee_\TAU^k\Vdots
(\vv_\TAU^k\Cdot\nabla)\nabla\Ee_\TAU^k\Big)\,\d\xx
\\[-.2em]&\hspace{.1em}
=
\!\int_\varOmega\!|\nabla\Ee_\TAU^k|^{\EXS-2}(\nabla\Ee_\TAU^k{\boxtimes}\nabla\Ee_\TAU^k)\Colon\strain(\vv_\TAU^k)-\frac1\EXS({\rm div}\,\vv_\TAU^k)|\nabla\Ee_\TAU^k|^\EXS\,\d\xx\,,\!\!
\label{calulus-nonlin-hyper}\end{align}
where the product $\boxtimes$ of the 3rd-order tensors is defined as
$[\bm G\boxtimes\bm G]_{ij}=\sum_{k,l=1}^3G_{ikl}G_{jkl}$. Relying on the already
obtained estimate \eq{Euler-small-E-Lr} and using the discrete Gronwall
inequality for sufficiently small $\tau>0$, from \eq{nabla-E-calulus}
we obtain
\begin{align}\label{Euler-small-E-Wr}
&\|\nabla\oEeT\|_{L^\infty(I;L^\EXS(\varOmega;\R^{3\times3\times3}))}^{}\le C_{p,\EXS}
\ \ \text{ with any }\ 1\le \EXS<\EXP
\,;
\end{align}
here we have used the assumption $\Ee_0\in W^{1,\EXS}(\varOmega;\Rsym)$, cf.\
\eq{Euler-thermo-small-ass-IC}. Here, for the aforementioned discrete Gronwall
inequality, we~have~used~that~both
$\tau\sum_{k=1}^{T/\tau}\|\nabla\vv_\TAU^k\|_{L^\infty(\varOmega;\R^{3\times3})}^p$ is
bounded for $p>1$ and
$\tau\sum_{k=1}^{T/\tau}\|\nabla\theta_\TAU^k\|_{L^\EXP(\varOmega;\R^3)}^\EXP$
is bounded for $\EXP>1$. Thus both
$\tau\max_{k=1,...,T/\tau}\|\nabla\vv_\TAU^k\|_{L^\infty(\varOmega;\R^{3\times3})}
=\mathscr{O}(\tau^{1-1/p})$ and
$\tau\max_{k=1,...,T/\tau}\|\nabla\theta_\TAU^k\|_{L^\EXP(\varOmega;\R^3)}^{}
=\mathscr{O}(\tau^{1-1/\EXP})$ are sufficiently small for $\tau$ small enough.

Moreover, by comparison, from
$\pdt{}\Ee_\TAU=\strain(\ovT)-\bm B_\text{\sc zj}^{}(\ovT,\oEeT)-\mathscr{R}(\oEeT,\otT)$,
we also obtain
\begin{align}\label{Euler-small-therm-ED-est3}
\Big\|\pdt{\Ee_\TAU}\Big\|_{L^p(I;L^\EXS(\varOmega;\R^{3\times3}))}^{}\le C_{r,\EXS}\,.
\end{align}

\medskip{\it Step 5: Existence of a solution to \eq{Euler-small-therm-ED+discr}}.
The rigorous proof of the existence of a solution
$(\varrho_\tau^k,\pp_\tau^k,\Ee_\tau^k,\theta_\tau^k)$ and thus also $\vv_\tau^k$
of the system \eq{Euler-small-therm-ED+discr} is slightly delicate.
Convincing arguments may rely on a suitable regularization towards
a quasilinear elliptic system for which standard theory can be used, and
then making a limit passage.

Given $(\varrho_{\tau}^{k-1},\pp_{\tau}^{k-1},\Ee_{\tau}^{k-1},\theta_{\tau}^{k-1})\in
W^{1,r}(\varOmega)\times W^{1,r}(\varOmega;\R^3)\times W^{1,s}(\varOmega;\Rsym)
\times L^{1+\alpha}(\varOmega)$ with $\varrho_{\tau}^{k-1}>0$ and
$\theta_{\tau}^{k-1}>0$ a.e., we seek a solution
$(\varrho_{\EPS\DELTA\tau}^k,\vv_{\EPS\DELTA\tau}^k,\Ee_{\EPS\DELTA\tau}^k,\theta_{\EPS\DELTA\tau}^k)$
and thus also $\pp_{\EPS\DELTA\tau}^k$ and $\Ent_{\EPS\DELTA\tau}^k$
of the regularized quasilinear elliptic system 
\begin{subequations}\label{Euler-small-therm-reg-discr}
\begin{align}\label{Euler-small-therm-reg-0discr}
&\!\!\varrho_{\EPS\DELTA\tau}^k+
{\rm div}\big(\tau\pp_{\EPS\DELTA\tau}^k{-}\DELTA|\nabla\varrho_{\EPS\DELTA\tau}^k|^{r-2}\nabla\varrho_{\EPS\DELTA\tau}^k\big)=\varrho_\TAU^{k-1}
\ \ \ \ \text{ with }\ \ \pp_{\EPS\DELTA\tau}^k=\varrho_{\EPS\DELTA\tau}^k\vv_{\EPS\DELTA\tau}^k\,,\!\!
\\[.1em]&\nonumber
\!\!\pp_{\EPS\DELTA\tau}^k+
\tau{\rm div}\big(\pp_{\EPS\DELTA\tau}^k{\otimes}\vv_{\EPS\DELTA\tau}^k
{-}\DD_{\EPS\DELTA\tau}^k\big)
\,=\pp_\TAU^{k-1}+\tau\varrho_{\EPS\DELTA\tau}^k\GRAVITY_{\tau}^k
+\tau{\rm div}\,\mathscr{T}(\Ee_{\EPS\DELTA\tau}^k,\theta_{\EPS\DELTA\tau}^k)
\\[-.1em]&\nonumber\hspace*{17.4em}
-\EPS|\vv_{\EPS\DELTA\tau}^k|^{p-2}\vv_{\EPS\DELTA\tau}^k-\DELTA|\nabla\varrho_{\EPS\DELTA\tau}^k|^{r-2}(\nabla\vv_{\EPS\DELTA\tau}^k)\nabla\varrho_{\EPS\DELTA\tau}^k\,,\ 
\\[-.1em]&
\hspace*{6.5em}\text{ where }\ 
\DD_{\EPS\DELTA\tau}^k\!=
\bbD\strain(\vv_{\EPS\DELTA\tau}^k)-{\rm div}\mathfrak{H}_{\EPS\DELTA\tau}^k\ 
\text{ with }\ \mathfrak{H}_{\EPS\DELTA\tau}^k\!
=\HYPER\big|\nabla^2\vv_{\EPS\DELTA\tau}^k\big|^{p-2}\nabla^2\vv_{\EPS\DELTA\tau}^k\,,\!\!
\label{Euler-small-therm-reg-1discr}
\\[.1em]&\nonumber
\!\!\Ee_{\EPS\DELTA\tau}^k=\Ee_\TAU^{k-1}\!+
\tau\strain(\vv_{\EPS\DELTA\tau}^k)-\tau\ZJEp_{\EPS\DELTA\tau}^k\!
-\tau\bm B_\text{\sc zj}^{}(\vv_{\EPS\DELTA\tau}^k,\Ee_{\EPS\DELTA\tau}^k)
+{\rm div}\big(\EPS|\nabla\Ee_{\EPS\DELTA\tau}^k|^{s-2}\nabla\Ee_{\EPS\DELTA\tau}^k\big)
\\[-.1em]&\hspace*{24em}
\ \text{ with }\,\ZJEp_{\EPS\DELTA\tau}^k
=\mathscr{R}(\Ee_{\EPS\DELTA\tau}^k,\theta_{\EPS\DELTA\tau}^k)\,,
\label{Euler-small-therm-reg-2discr}
\\[.1em]&\nonumber\!\!\Ent_{\EPS\DELTA\tau}^k+\tau
{\rm div}\big(\Ent_{\EPS\DELTA\tau}^k\vv_{\EPS\DELTA\tau}^k{-}\kappa(\theta_{\EPS\DELTA\tau}^k)\nabla\theta_{\EPS\DELTA\tau}^k\!\big)=\Ent_\TAU^{k-1}+\tau\,
\xi_{\rm ext}^{}(\theta_{\EPS\DELTA\tau}^k;\vv_{\EPS\DELTA\tau}^k,\ZJEp_{\EPS\DELTA\tau}^k)
\\[-.1em]&\hspace*{7.5em}
+\tfrac\tau3{\rm tr}\big(\mathscr{T}(\Ee_{\EPS\DELTA\tau}^k,\theta_{\EPS\DELTA\tau}^k)
-\mathscr{T}(\Ee_{\EPS\DELTA\tau}^k,0)\big)\,{\rm div}\,\vv_{\EPS\DELTA\tau}^k
\ \ \ \text{ with $\ \ \Ent_{\EPS\DELTA\tau}^k\!=\ENT(\theta_{\EPS\DELTA\tau}^k)$}\,,
\label{Euler-small-therm-reg-3discr}
\end{align}\end{subequations}
considered with $r>3$ and $s>3$. For the $\DELTA$-terms in
(\ref{Euler-small-therm-reg-discr}a,b) in the case $r=2$
(not used here) we refer to \cite{Zato12ASCN}.
Of course, we consider the corresponding boundary conditions \eq{BC-disc}
now written for the $(\EPS,\DELTA)$-regularization and 
augmented also by $\nn\Cdot\nabla\varrho_{\EPS\DELTA\tau}^k=0$ and
by $(\nn\Cdot\nabla)\Ee_{\EPS\DELTA\tau}^k=\bm0$ on $\varGamma$.

Like \eq{Euler-thermo-basic-energy-balance-disc}, we obtain the inequality
\begin{align}\nonumber
\!\!\!&\int_\varOmega\!\bigg(\frac{|\pp_{\EPS\DELTA\tau}^k|^2\!\!\!}{2\varrho_{\EPS\DELTA\tau}^k\!\!\!}
+\varphi(\Ee_{\EPS\DELTA\tau}^k)+\ENT(\theta_{\EPS\DELTA\tau}^k)
+\tau[\zeta_{\rm p}]'_{\Lp}(\theta_{\EPS\DELTA\tau}^k,\Lp_{\EPS\DELTA\tau}^k)\Colon\Lp_{\EPS\DELTA\tau}^k
\\[-.5em]\!\!\!&\hspace{5em}\nonumber
+\EPS\varphi''(\Ee_{\EPS\DELTA\tau}^k)\nabla\Ee_{\EPS\DELTA\tau}^k\Vdots\big(|\nabla\Ee_{\EPS\DELTA\tau}^k|^{s-2}\nabla\Ee_{\EPS\DELTA\tau}^k\big)
+\EPS|\vv_{\EPS\DELTA\tau}^k|^p\bigg)\,\d\xx
+\tau\!\!\int_\varGamma\!h(\theta_{\EPS\DELTA\tau}^k)\,\d S
\\[-.5em]\!\!\!&\hspace{8em}
\le\int_\varOmega\!\frac{|\pp_{\tau}^{k-1}|^2\!\!}{2\varrho_{\tau}^{k-1}\!\!}
+\varphi(\Ee_{\tau}^{k-1})+\ENT(\theta_{\tau}^{k-1})
+\tau\varrho_{\EPS\DELTA\tau}^k\GRAVITY_{\tau}^k\Cdot\vv_{\EPS\DELTA\tau}^k
\,\d\xx+\!\!\int_\varGamma\!h_{{\rm ext},\tau}^k\,\d S\,.\!\!
\label{Euler-large-basic-engr-balance-disc+}\end{align}
Here we used the cancellation of the two $\DELTA$-regularizing terms as in
\cite{Roub25TDVE,Zato12ASCN}. From \eq{Euler-large-basic-engr-balance-disc+},
treating the term $\tau\varrho_{\EPS\DELTA\tau}^k\GRAVITY_{\tau}^k\Cdot\vv_{\EPS\DELTA\tau}^k$
as in \eq{Euler-small-est-Gronwall-}, we can directly read
the a-priori estimates for $\vv_{\EPS\DELTA\tau}^k\in L^p(\varOmega;\R^3)$,
$\Ee_{\EPS\DELTA\tau}^k\in L^2(\varOmega;\Rsym)$,  and
$\theta_{\EPS\DELTA\tau}^k\in L^{1+\alpha}(\varOmega)$. Here we also used that
$\theta_{\EPS\DELTA\tau}^k\ge0$, which can be seen by testing
\eq{Euler-small-therm-reg-3discr} by $(\theta_{\EPS\DELTA\tau}^k)^-$.

Then, test of \eq{Euler-small-therm-reg-1discr} by $\vv_{\EPS\DELTA\tau}^k$
while using \eq{Euler-small-therm-reg-0discr} multiplied by
$|\vv_{\EPS\DELTA\tau}^k|^2/2$ and \eq{Euler-small-therm-reg-2discr} multiplied by
$\varphi(\Ee_{\EPS\DELTA\tau}^k)$ gives
\begin{align}\nonumber
\!\!\!\!\!&\int_\varOmega\!\bigg(\frac{|\pp_{\EPS\DELTA\tau}^k|^2\!\!\!}{2\varrho_{\EPS\DELTA\tau}^k\!\!\!}
+\varphi(\Ee_{\EPS\DELTA\tau}^k)
+\EPS|\vv_{\EPS\DELTA\tau}^k|^p+\tau\bbD\strain(\vv_{\EPS\DELTA\tau}^k)\Colon\strain(\vv_{\EPS\DELTA\tau}^k)
+\tau\HYPER|\nabla^2\vv_{\EPS\DELTA\tau}^k|^p\bigg)\,\d\xx
\\[-.1em]&\hspace{1em}
\le\int_\varOmega\!\frac{|\pp_{\tau}^{k-1}|^2\!\!}{2\varrho_{\tau}^{k-1}\!\!}
+\varphi(\Ee_{\tau}^{k-1})
+\tau\varrho_{\EPS\DELTA\tau}^k\GRAVITY_{\tau}^k\Cdot\vv_{\EPS\DELTA\tau}^k
+\frac13{\rm tr}\big(\mathscr{T}(\Ee_{\EPS\DELTA\tau}^k,\theta_{\EPS\DELTA\tau}^k)
{-}\mathscr{T}(\Ee_{\EPS\DELTA\tau}^k,0)\big)\,{\rm div}\,\vv_{\EPS\DELTA\tau}^k
\,\d\xx\,.\!\!
\label{Euler-large-basic-engr-balance-disc++}\end{align}
Here the terms $[\zeta_{\rm p}]'_{\Lp}(\theta_{\EPS\DELTA\tau}^k,\Lp_{\EPS\DELTA\tau}^k)\Colon\Lp_{\EPS\DELTA\tau}^k\ge0$
and $\EPS\varphi''(\Ee_{\EPS\DELTA\tau}^k)\nabla\Ee_{\EPS\DELTA\tau}^k\Vdots((|\nabla\Ee_{\EPS\DELTA\tau}^k|^{s-2}\nabla\Ee_{\EPS\DELTA\tau}^k)\ge0$
in \eq{Euler-large-basic-engr-balance-disc+} are omitted because they would not
yield a useful information. Treating the right-hand side as in
\eq{Euler-small-est2-}, this gives the further a-priori estimates
$\nabla\vv_{\EPS\DELTA\tau}^k\in W^{1,p}(\varOmega;\R^{3\times3})$. Then, testing
\eq{Euler-small-therm-reg-0discr} by $\varrho_{\EPS\DELTA\tau}^k$ leads to an
estimate \eq{transport-rho+} augmented by the left-hand-side term
$\DELTA|\nabla\varrho_{\EPS\DELTA\tau}^k|^r$ and thus an estimate of
$\varrho_{\EPS\DELTA\tau}^k\in W^{1,r}(\varOmega)$ and, from \eq{calulus-nonlin-rho-r+}
augmented correspondingly, we still obtain an estimate for
${\rm div}(\DELTA|\nabla\varrho_{\EPS\DELTA\tau}^k|^{r-2}\nabla\varrho_{\EPS\DELTA\tau}^k)$
in $L^2(\varOmega)$.

Moreover, by the same arguments as in Step~3, we can show that
$\theta_{\EPS\DELTA\tau}^k>0$ a.e.\ on $\varOmega$. Then, we can test
\eq{Euler-small-therm-reg-3discr} by $(\theta_{\EPS\DELTA\tau}^k)^{-\lambda}$, and
thus we arrive at a direct analogue of \eq{Euler-entropy-finite-disc-}, namely
\begin{align}&\nonumber
  \!\int_\varOmega\!\frac{\xi_{\rm ext}^{}(\theta_{\EPS\DELTA\tau}^k;\vv_{\EPS\DELTA\tau}^k,\ZJEp_{\EPS\DELTA\tau}^k
)\!}{(\theta_\TAU^k)^\lambda}
+\kappa(\theta_\TAU^k)\frac{|\nabla\theta_{\EPS\DELTA\tau}^k|^2\!\!\!}{(\theta_{\EPS\DELTA\tau}^k)^{1+\lambda}}\,\d\xx
\le\!\!\int_\varOmega\!\frac{\eta_\lambda(\theta_{\EPS\DELTA\tau}^k)-\eta_\lambda(\theta_\TAU^{k-1})\!}\tau\,\d\xx
+\!\int_\varGamma\frac{h(\theta_{{\rm ext},\tau}^k)\!}
{(\theta_{{\rm ext},\tau}^k)^\lambda}\,\d S
\\[-.4em]&\hspace{4em}
-\!\int_\varOmega\!\bigg(
\theta_{\EPS\DELTA\tau}^k\frac{\phi'({\rm tr}\Ee_{\EPS\DELTA\tau}^k){+}\phi({\rm tr}\Ee_{\EPS\DELTA\tau}^k)\!}
{(\theta_\TAU^k)^{\lambda}}
+\frac{\theta_{\EPS\DELTA\tau}^k\gamma'(\theta_{\EPS\DELTA\tau}^k){+}\gamma(\theta_{\EPS\DELTA\tau}^k)\!\!}
{(\theta_{\EPS\DELTA\tau}^k)^\lambda}
+\eta_\lambda(\theta_{\EPS\DELTA\tau}^k)\bigg)\,{\rm div}\,\vv_{\EPS\DELTA\tau}^k\,\d\xx
\label{Euler-entropy-finite-disc-k}\end{align}

Then, by testing \eq{Euler-small-therm-reg-2discr} with
$|\Ee_{\EPS\DELTA\tau}^k|^{r-2}\Ee_{\EPS\DELTA\tau}^k$, by the calculus
\eq{transport-}--\eq{transport} for $\tau>0$ sufficiently small, we obtain
a bound of $\Ee_{\EPS\DELTA\tau}^k$ in $L^r(\varOmega;\Rsym)$ for $r$ arbitrarily
large. Note that $\int_\varOmega{\rm div}(\EPS|\nabla\Ee_{\EPS\DELTA\tau}^k|^{s-2}\nabla\Ee_{\EPS\DELTA\tau}^k)\Colon
(|\Ee_{\EPS\DELTA\tau}^k|^{r-2}\Ee_{\EPS\DELTA\tau}^k)\,\d\xx=
-\int_\varOmega\EPS|\nabla\Ee_{\EPS\DELTA\tau}^k|^{s-2}\nabla\Ee_{\EPS\DELTA\tau}^k\Vdots
\nabla(|\Ee_{\EPS\DELTA\tau}^k|^{r-2}\Ee_{\EPS\DELTA\tau}^k)\,\d\xx=
(1{-}r)\int_\varOmega\EPS|\Ee_{\EPS\DELTA\tau}^k|^{r-2}|\nabla\Ee_{\EPS\DELTA\tau}^k|^s
\,\d\xx\le0$. Then, by interpolation calculus as in
\eq{Euler-entropy-nabla-theta}--\eq{8-cond}, we obtain a bound of
$\nabla\theta_{\EPS\DELTA\tau}^k$ as in \eq{Euler-entropy-nabla-theta++}, i.e., here
in $L^\EXP(\varOmega;\R^3)$ with $\EXP$ from \eq{Euler-entropy-nabla-theta-mu}.

Also, testing \eq{Euler-small-therm-reg-2discr} by
${\rm div}(|\nabla\Ee_{\EPS\DELTA\tau}^k|^{s-2}\nabla\Ee_{\EPS\DELTA\tau}^k)$
gives an additional estimate for
$\nabla\Ee_{\EPS\DELTA\tau}^k\in L^s(\varOmega;\R^{3\times3\times3})$ and for
${\rm div}(|\nabla\Ee_{\EPS\DELTA\tau}^k|^{s-2}\nabla\Ee_{\EPS\DELTA\tau}^k)$ itself.
More specifically,
\begin{subequations}\label{Euler-small-est-regul}
\begin{align}\label{Euler-small-est-regul1}
&\|\varrho_{\EPS\DELTA\tau}^k\|_{W^{1,r}(\varOmega)}^{}\!\le C\ \  
\ \ \text{ and }\ \ 
\big\|{\rm div}(|\nabla\varrho_{\EPS\DELTA\tau}^k|^{r-2}\nabla\varrho_{\EPS\DELTA\tau}^k)\big\|_{L^2(\varOmega)}^{}\!\le C/\sqrt\DELTA\,,
\\[-.3em]&\label{Euler-small-est-regul2}
\|\nabla\vv_{\EPS\DELTA\tau}^k\|_{W^{1,p}(\varOmega;\R^{3\times3})}^{}\le C\,,\ \ \
\|\vv_{\EPS\DELTA\tau}^k\|_{L^p(\varOmega;\R^3)}^{}\le C/\sqrt[p]\EPS
\,,\ \ \text{ and }\ \
\bigg\|\frac{\pp_{\EPS\DELTA\tau}^k}{\sqrt\varrho_{\EPS\DELTA\tau}^k}
\bigg\|_{L^2(\varOmega;\R^3)}^{}\!\!\le C\,,
\\[-.2em]&\label{Euler-small-est-regul3}
\|\Ee_{\EPS\DELTA\tau}^k\|_{W^{1,s}(\varOmega;\R^{3\times3})}^{}\le C
\ \ \text{ and }\ \ \big\|{\rm div}(|\nabla\Ee_{\EPS\DELTA\tau}^k|^{s-2}\nabla\Ee_{\EPS\DELTA\tau}^k)\big\|_{L^2(\varOmega;\R^{3\times3\times3})}^{}\le C/\sqrt\EPS
\,,
\ \ \text{ and}
\\&\|\theta_{\EPS\DELTA\tau}^k\|_{L^{1+\alpha}(\varOmega)\,\cap\,W^{1,\EXP}(\varOmega)}^{}\le C\,.
\end{align}\end{subequations}

This allows for the existence of weak solutions to \eq{Euler-small-therm-reg-discr}
by rather standard methods for quasilinear elliptic problems when realizing that
the highest-order terms in each equations in (\ref{Euler-small-therm-reg-discr}a--c)
are monotone and when realizing that the equation \eq{Euler-small-therm-reg-3discr}
is semilinear, albeit with the $L^1$-right-hand side. Here, the compactness of
$\Ee$'s and $\theta$'s is important together with the continuity of the
nonlinearities $\mathscr{T}$ and $\mathscr{R}$ in the lower-order terms.

Then we pass to the limit with $\DELTA\to0$ by choosing a subsequence such that
\begin{subequations}\label{Euler-small-converge-k}
\begin{align}
&&&\!\!\varrho_{\EPS\DELTA\tau}^k\to\varrho_{\EPS\tau}^k&&\hspace{-3em}\text{weakly in $\ W^{1,r}(\varOmega)$\,,}
\\&&&\label{Euler-small-converge-bar-p-k}
\!\!\pp_{\EPS\DELTA\tau}^k\to\pp_{\EPS\tau}^k&&\hspace{-3em}\text{weakly in $\
W^{1,r}(\varOmega;\R^3)$\,,}
\\&&&\label{Euler-small-converge-bar-v-k}
\!\!\vv_{\EPS\DELTA\tau}^k\to\vv_{\EPS\tau}^k&&\hspace{-3em}\text{weakly in $\ W^{2,p}(\varOmega;\R^3)$\,,}
\\\label{Euler-small-converge-E-k}
&&&\!\!\Ee_{\EPS\DELTA\tau}^k\to\Ee_{\EPS\tau}^k\!\!\!&&\hspace{-3em}\text{weakly in $\ 
W^{1,\EXS}(\varOmega;\Rsym)$\,,\ and}\!\!&&
\\[-.2em]\label{Euler-thermo-Ch4-conv1-k}
&&&\!\!\theta_{\EPS\DELTA\tau}^k\to\theta_{\EPS\tau}^k&&\hspace{-3em}\text{weakly in }\
W^{1,\EXP}(\varOmega)\,.
&&
\end{align}\end{subequations}
According to the latter estimate in \eq{Euler-small-est-regul1}, the regularizing
term
${\rm div}(\DELTA|\nabla\varrho_{\EPS\DELTA\tau}^k|^{r-2}\nabla\varrho_{\EPS\DELTA\tau}^k)$
in \eq{Euler-small-therm-reg-0discr} is $\mathscr{O}(\sqrt\DELTA)$ in
$L^2(\varOmega)$ and thus vanishes in the limit.
Also, the compensating force $\DELTA|\nabla\varrho_{\EPS\DELTA\tau}^k|^{r-2}
(\nabla\vv_{\EPS\DELTA\tau}^k)\nabla\varrho_{\EPS\DELTA\tau}^k$ in
\eq{Euler-small-therm-reg-1discr} is $\mathscr{O}(\DELTA)$ in
$L^{r'}(\varOmega;\R^3)$ and thus vanishes in the limit. The
strong convergence (in terms of subsequences) of $\vv_{\EPS\DELTA\tau}^k$ in
$W^{2,p}(\varOmega;\R^3)$ can be proved due to the strong monotonicity of the
operator $\vv\mapsto
{\rm div}({\rm div}(\HYPER|\nabla^2\vv|^{p-2}\nabla^2\vv)-\bbD\strain(\vv))$.
In more detail, using \eq{Euler-small-therm-reg-1discr}, we obtain
\begin{align}\nonumber
\big(&\inf_{|E|=1}{\bbD}E\Colon E\big)
\|\strain(\vv_{\EPS\DELTA\tau}^k{-}\wt\vv)\|_{L^2(\varOmega;\R^{3\times3})}^2
+\mu c_p\|\nabla^2(\vv_{\EPS\DELTA\tau}^k{-}\wt\vv)\|_{L^p(\varOmega;\R^{3\times3\times3})}^p
\\[-.4em]&\nonumber\le
\int_\varOmega\bigg(
\Big(\varrho_{\EPS\DELTA\tau}^k\GRAVITY_{\tau}^k-\frac{\pp_{\EPS\DELTA\tau}^k{-}
\pp_\tau^{k-1}\!\!\!}\tau
-\tau\EPS|\vv_{\EPS\DELTA\tau}^k|^{p-2}\vv_{\EPS\DELTA\tau}^k
-\tau\DELTA|\nabla\varrho_{\EPS\DELTA\tau}^k|^{r-2}
(\nabla\vv_{\EPS\DELTA\tau}^k)\nabla\varrho_{\EPS\DELTA\tau}^k\Big)
\Cdot(\vv_{\EPS\DELTA\tau}^k{-}\wt\vv)
\\[-.5em]&\nonumber\hspace{1em}
-\Big(
\bbD\strain(\wt\vv){+}\mathscr{T}(\Ee_{\EPS\DELTA\tau}^k,\theta_{\EPS\DELTA\tau}^k)
{-}\pp_{\EPS\DELTA\tau}^k{\otimes}\vv_{\EPS\DELTA\tau}^k\Big)\Colon\strain(\vv_{\EPS\DELTA\tau}^k{-}\wt\vv)
-\HYPER|\nabla^2\wt\vv|^{p-2}\nabla^2\wt\vv
\Vdots\nabla^2(\vv_{\EPS\DELTA\tau}^k{-}\wt\vv)\bigg)\,\d\xx
\\[-.7em]&\nonumber
\!\!\stackrel{\DELTA\to0}{\to}\!\!
 \int_\varOmega\bigg(
 \Big(\varrho_{\EPS\tau}^k\GRAVITY_{\tau}^k-\frac{\pp_{\EPS\tau}^k{-}
\pp_\tau^{k-1}\!\!\!}\tau-\tau\EPS|\vv_{\EPS\tau}^k|^{p-2}\vv_{\EPS\tau}^k\Big)
\Cdot(\vv_{\EPS\tau}^k{-}\wt\vv)
\\[-.5em]&\hspace{1em}
-\Big(\bbD\strain(\wt\vv){+}\mathscr{T}(\Ee_{\EPS\tau}^k,\theta_{\EPS\tau}^k)
{-}\pp_{\EPS\tau}^k{\otimes}\vv_{\EPS\tau}^k\Big)\Colon\strain(\vv_{\EPS\tau}^k{-}\wt\vv)
-\HYPER|\nabla^2\wt\vv|^{p-2}\nabla^2\wt\vv
\Vdots\nabla^2(\vv_{\EPS\tau}^k{-}\wt\vv)\bigg)\,\d\xx
\label{Euler-large-strong+k}\end{align}
for some $c_p>0$ and for any $\wt\vv\in W^{2,p}(\varOmega;\R^3)$.
We have also used the strong convergence
$\mathscr{T}(\Ee_{\EPS\DELTA\tau}^k,\theta_{\EPS\DELTA\tau}^k)\to
\mathscr{T}(\Ee_{\EPS\tau}^k,\theta_{\EPS\tau}^k)$ in $L^1(\varOmega;\Rsym)$.
Choosing $\wt\vv=\vv_{\EPS\tau}^k$, we obtain the mentioned strong convergence
$\nabla\vv_{\EPS\DELTA\tau}^k\to\nabla\vv_{\EPS\tau}^k$ in $W^{1,p}(\varOmega;\R^3)$.
Similarly, the monotonicity of the operator
$\Ee\mapsto-{\rm div}(\EPS|\nabla\Ee|^{s-2}\nabla\Ee)$ allows for making the
convergence in \eq{Euler-small-converge-E-k} strong. 

We thus obtain a solution 
$(\varrho_{\EPS\tau}^k,\vv_{\EPS\tau}^k,\Ee_{\EPS\tau}^k,\theta_{\EPS\tau}^k)$
and thus also $\pp_{\EPS\tau}^k$ and $\Ent_{\EPS\tau}^k$
of the system  \eq{Euler-small-therm-reg-discr} with
$\DELTA$ omitted. In particular, \eq{Euler-small-therm-reg-0discr}
turns into $\varrho_{\EPS\tau}^k+
{\rm div}(\tau\varrho_{\EPS\tau}^k\vv_{\EPS\tau}^k)=\varrho_\TAU^{k-1}$
so that, by the estimation as in
\eq{ineq-for-sparsity}--\eq{Euler-large-est-v-2-eps} when using also
the last estimate in \eq{Euler-small-est-regul2},
we obtain a bound for $\vv_{\EPS\tau}^k$ in $L^a(\varOmega;\R^3)$
for any $1\le a<2$ independent of $\EPS$.
Together with the first estimate in \eq{Euler-small-est-regul2},
we have $\vv_{\EPS\tau}^k$ bounded in $W^{2,p}(\varOmega;\R^3)$.

This allows us to pass to the limit with $\EPS\to0$ in terms of subsequences. The
$\EPS$-regularizing in \eq{Euler-small-therm-reg-1discr} with
$\DELTA$ omitted, i.e.\ now $\EPS|\vv_{\EPS\tau}^k|^{p-2}\vv_{\EPS\tau}^k$,
is $\mathscr{O}(\EPS)$ in $L^\infty(\varOmega;\R^3)$ so that it vanishes
in the limit for $\EPS$. The latter estimate in \eq{Euler-small-est-regul3}
makes the regularizing term 
${\rm div}(\EPS|\nabla\Ee_{\EPS\DELTA\tau}^k|^{s-2}\nabla\Ee_{\EPS\DELTA\tau}^k)$
as $\mathscr{O}(\sqrt\EPS)$ in $L^2(\varOmega;\Rsym)$, so that it vanishes in
the limit for $\EPS\to0$. Altogether, we obtain
$(\varrho_{\tau}^k,\vv_{\tau}^k,\Ee_{\tau}^k,\theta_{\tau}^k)$
and thus also $\pp_{\tau}^k$ and $\Ent_{\tau}^k$ solving
the boundary-value problem \eq{Euler-small-therm-ED+discr}--\eq{BC-disc}.

\medskip{\it Step 6: Convergence of \eq{Euler-small-therm-ED+dis}
for $\tau\to0$}.
By the Banach selection principle, we obtain a subsequence
converging weakly* with respect to the topologies indicated in
\eq{Euler-small-est}, \eq{Euler-small-est2}, \eq{est-of-rho-disc},
\eq{Euler-large-est-v-2-eps+}--\eq{Euler-small-est++},
\eq{Euler-entropy-nabla-theta++++}, \eq{Euler-small-E-Wr}, and
\eq{Euler-small-therm-ED-est3}, to some limit
$(\varrho,\pp,\vv,\Ee,\Ent)$. Specifically,
\begin{subequations}\label{Euler-small-converge}
\begin{align}
&&&\!\!\orT\to\varrho&&\text{weakly* in $\ L^\infty(I;W^{1,r}(\varOmega))$\,,}
\\&&&\!\!\varrho_\TAU\to\varrho&&\text{weakly* in $\ L^\infty(I;W^{1,r}(\varOmega))\,\cap\, 
 W^{1,p}(I;L^r(\varOmega))$\,,}
\\&&&\label{Euler-small-converge-bar-p}
\!\!\opT\to\pp&&\text{weakly\ \;in $\
L^p(I;W^{1,r}(\varOmega;\R^3))$\,,}
\\&&&\label{Euler-small-converge-p}
\!\!\pp_\TAU\to\pp&&\text{weakly\ \;in $\
L^p(I;W^{1,r}(\varOmega;\R^3))\,\cap\,W^{1,p'}(I;W^{2,p}(\varOmega;\R^3)^*)$\,,}
\\&&&\label{Euler-small-converge-bar-v}
\!\!\ovT\to\vv&&\text{weakly* in $\  L^\infty(I;L^a(\varOmega;\R^3))
\,\cap\, L^p(I;W^{2,p}(\varOmega;\R^3))$ with $1\le a<2$\,,}
\\\label{Euler-small-converge-bar-E}
&&&\!\!\oEeT\to\Ee\!\!\!&&\text{weakly* in $\ 
 L^\infty(I;W^{1,\EXS}(\varOmega;\Rsym))
\,\cap\,L^\infty(I;L^2(\varOmega;\Rsym))$\,,}\!\!
\\
&&&\!\!\Ee_\TAU\to\Ee\!\!\!&&\text{weakly*\ \;in $\ 
 L^\infty(I;W^{1,\EXS}(\varOmega;\Rsym))
\,\cap\,W^{1,p}(I;L^\EXS(\varOmega;\Rsym))$\,, and}&&
\\[-.2em]\label{Euler-thermo-Ch4-conv1}
&&&\!\!\ouT\to\Ent&&\text{weakly in }\
L^{(4+\alpha)\,\EXP/(3+3\alpha)}(I{\times}\varOmega)\,.
&&
\end{align}\end{subequations}
Notably, the limit{s}  of $\orT$ and $\varrho_\tau$  are 
indeed the same
due to the control of $\pdt{}\varrho_\TAU$ in \eq{Euler-small-est6}; cf.\
\cite[Sect.8.2]{Roub13NPDE}. The same holds true also for $\opT$
and $\pp_\TAU$, and for $\oEeT$ and $\Ee_\TAU$, too.

By the compact embedding $W^{1,r}(\varOmega)\subset C(\barOmega)$ for $r>3$
and the (generalized) Aubin-Lions theorem, cf.\
\cite[Corollary~7.9]{Roub13NPDE}, we have also 
\begin{subequations}\label{Euler-small-converge-strong}
\begin{align}\label{Euler-small-converge-strong-rho}
&&&\!\!\orT\to\varrho\!\!\!&&\text{strongly in }\
L^a(I;C(\barOmega))\ \ \text{ for any }\ 1\le a<\infty\,,&&
\\\label{Euler-small-converge-strong-p-}
&&&\!\!\opT\to\pp\!\!\!\!\!\!&&\text{strongly in }\ L^p(I;C(\barOmega;\R^3))\,,
\ \text{ and}
\\\label{Euler-small-converge-strong-E}
&&&\!\!\oEeT\to\Ee\!\!\!\!\!&&\text{strongly in }\
L^a(I;L^2(\varOmega;\Rsym))\ \text{ for any }\ 1\le a<\infty\,.
\intertext{Moreover, using the Arzel\`a-Ascoli-type theorem, cf.\
\cite[Lemma~7.10]{Roub13NPDE}, we have also}
\label{Euler-small-converge-strong-rho+}
&&&\!\!\varrho_\TAU\to\varrho\!\!\!&&\text{strongly in }\ C(I{\times}\barOmega)\,.
\end{align}\end{subequations}

Moreover, from (\ref{Euler-small-converge-strong}a,b), we know that, for a
subsequence of $\tau\to0$, $\opT/\orT$ converges to some $\wt\vv$
a.e.\ on $I{\times}\varOmega$. Simultaneously, we know
$\sup_{\tau>0}\int_0^T\int_\varOmega|\wt\vv-\opT/\orT|^p\,\d\xx\d t<\infty$
with $p>1$ so that, by the de\,la\,Vall\'ee\,Poussin theorem,
$\{|\wt\vv{-}\opT/\orT|^{q}\}_{\tau>0}$ is relatively
weakly compact in $L^1(I{\times}\varOmega)$ for any $1\le q<p$. Then,
by the Dunford-Pettis theorem, it is uniformly integrable and, by
the Vitali theorem, it converges a.e.\ to its limit which equals 0,
i.e.\ it converges to 0 in $L^q(I{\times}\varOmega)$. By
\eq{Euler-small-converge-bar-v}, we can identify
$\wt\vv=\vv$ so that all the (already chosen subsequence for
\eq{Euler-small-converge}) converges to $\vv$. This proves 
\begin{subequations}\label{Euler-small-converge-strong+}
\begin{align}
&\ovT=\frac{\opT}{\orT}\to\frac{\pp}{\varrho}=\vv\!\!\!
&&
\label{Euler-small-converge-strong-vv}
\hspace{-3em}\text{strongly in $\ L^q(I{\times}\varOmega;\R^3)\,$
for any $1\le q<p$\,,}
\intertext{Actually, by interpolation with
\eq{Euler-small-converge-bar-v}, we have the strong convergence $\ovT\to\vv$
even in a better space
$L^s(I;L^a(\varOmega;\R^3))\cap L^q(I;C(\barOmega;\R^3))$
with any $s<\infty$ and  $1\le a<2$. Thus,
by \eq{Euler-small-converge-strong-p-}, also}
&\opT{\otimes}\ovT\to\pp{\otimes}\vv=\varrho\vv{\otimes}\vv\!
&&\hspace{-3em}\text{strongly in }\,
L^q(I;L^a(\varOmega;\Rsym)),\ \,1\le q<p,\ 1\le a<2.\hspace{-5em}
\label{Euler-small-converge-strong-pxv}
\end{align}\end{subequations}

The bound \eq{Euler-entropy-nabla-theta++++} and, by comparison, the bound
for $\pdt{}\Ent_\TAU$ in $L^1(I;W^{2,6}(\varOmega)^*)$
allow us to use the (generalized) Aubin-Lions theorem, cf.\
\cite[Corollary~7.9]{Roub13NPDE}, so that the convergence
\eq{Euler-thermo-Ch4-conv1} is even strong.
Realizing $\otT=\ENT^{-1}(\ouT)$ and the estimate
\eq{Euler-entropy-nabla-theta+++}, we use the continuity of the Nemytski\u{\i}
mapping and obtain
\begin{subequations}\label{Euler-thermo-Ch4-conv-kappa-T-R}\begin{align}
\label{Euler-thermo-Ch4-conv2}
&\otT\!\to\theta=\ENT^{-1}(\Ent)
\ \text{ strongly in } L^a(I{\times}\varOmega)\
\text{for any $1\le a<\frac{(4+\alpha)\,\EXP}3$}\,.
\intertext{In view of \eq{Euler-hat-kappa-est}, we thus also have}
&\label{Euler-thermo-Ch4-conv-hat-kappa}
\intkappa(\otT)\!
\to\intkappa(\theta)\qquad\text{ strongly in }\
L^{(4+\alpha)\,\EXP/(3+3\beta)}(I{\times}\varOmega).
\intertext{Furthermore, for the Cauchy stress
and for the inelastic-rate mapping, we have}
&\mathscr{T}(\oEeT,\otT)\to\mathscr{T}(\Ee,\theta)
\quad\text{ strongly in }\ L^1(I{\times}\varOmega;\Rsym)\ \text{ and}
\label{Euler-thermo-Ch4-conv3}
\\
&\oPT=\mathscr{R}(\oEeT,\otT)\to\mathscr{R}(\Ee,\theta)=:\ZJEp
\ \text{ strongly in }\ L^2(I{\times}\varOmega;\Rdev)\,;
\label{Euler-thermo-Ch4-conv4}
\end{align}\end{subequations}
here we have used the growth conditions \eq{Euler-thermo-small-ass-Tth} and
\eq{Euler-thermo-small-ass-zeta+} and the continuity of $\mathscr{T}$ and
$\mathscr{R}$.

This already allows for the limit passage in
\eq{Euler-small-therm-ED+dis}. While the limit passage
in \eq{Euler-small-therm-ED+0dis} is simple
due to (\ref{Euler-small-converge}b,c), the limit passage in the quasilinear
momentum equation \eq{Euler-small-therm-ED+1dis} is a bit more technical.
To this aim, we use the uniform monotonicity of the operator
$\vv\mapsto{\rm div}({\rm div}(\HYPER|\nabla^2\vv|^{p-2}\nabla^2\vv)
-\bbD\strain(\vv))$ with the boundary conditions \eq{BC-disc+}
and, using \eq{Euler-small-therm-ED+1dis} tested by $\ovT{-}\wt\vv$,
we obtain
\begin{align}\nonumber
\big(&\inf_{|E|=1}{\bbD}E\Colon E\big)
\|\strain(\ovT{-}\wt\vv)\|_{L^2(I\times\varOmega;\R^{3\times3})}^2
+\mu c_p\|\nabla^2(\ovT{-}\wt\vv)\|_{L^p(I\times\varOmega;\R^{3\times3\times3})}^p
\\[-.4em]&\nonumber\ \ \le\int_0^T\!\!\!\int_\varOmega\!
\bbD\strain(\ovT{-}\wt\vv)
\Colon\strain(\ovT{-}\wt\vv)
 +\HYPER\Big(|\nabla^2\ovT|^{p-2}\nabla^2\ovT
-|\nabla^2\wt\vv|^{p-2}\nabla^2\wt\vv\Big)\Vdots\nabla^2(\ovT{-}\wt\vv)
\,\d\xx\d t
 \\[-.0em]&\ \ =\nonumber
 \int_0^T\!\!\!\int_\varOmega\bigg(
 \Big(\orT\overline\GRAVITY_{\tau}-\pdt{\pp_\tau}\Big)
 \Cdot(\ovT{-}\wt\vv)
-\Big(\mathscr{T}(\oEeT,\otT)
{-}\opT{\otimes}\ovT\Big)\Colon\strain(\ovT{-}\wt\vv)
\\[-.7em]&\nonumber\hspace{13em}
-\bbD\strain(\wt\vv)\Colon\strain(\ovT{-}\wt\vv)
-\HYPER|\nabla^2\wt\vv|^{p-2}\nabla^2\wt\vv\Vdots\nabla^2(\ovT{-}\wt\vv)
\bigg)\,\d\xx\d t
\\[-.7em]&\ \ \le\nonumber
 \int_0^T\!\!\!\int_\varOmega\bigg(\orT\GRAVITY_{\tau}\Cdot(\ovT{-}\wt\vv)
+\pdt{\pp_\tau}\Cdot\wt\vv-(\opT{\otimes}\ovT)\Colon\strain(\wt\vv)
-\big(\bbD\strain(\wt\vv){+}\mathscr{T}(\oEeT,\otT)\big)
\Colon\strain(\ovT{-}\wt\vv)
\\[-.4em]&\hspace{10em}
-\HYPER|\nabla^2\wt\vv|^{p-2}\nabla^2\wt\vv\Vdots\nabla^2(\ovT{-}\wt\vv)
\bigg)\,\d\xx\d t+
\int_\varOmega\!\frac{|\pp_0|^2}{2\varrho_0}-\frac{|\pp_\tau(T)|^2}{2\varrho_\tau(T)}\,\d\xx
\label{Euler-large-strong+}\end{align}
for some $c_p>0$ and for any $\wt\vv\in L^p(I;W^{2,p}(\varOmega;\R^3))$. 
The last inequality in \eq{Euler-large-strong+} has again exploited the
convexity of the kinetic energy $(\pp,\varrho)\mapsto\frac12|\pp|^2/\varrho$
in the calculus:
\begin{align}\nonumber
\int_\varOmega&\frac{|\pp_\tau(T)|^2\!}{2\varrho_\tau(T)}
-\frac{|\pp_0|^2}{2\varrho_0}\,\d\xx
\le\int_0^T\!\!\!\int_\varOmega\!\pdt{\pp_\tau\!}\Cdot\ovT
-\frac{|\ovT|^2\!}2\,\pdt{\varrho_\tau}\,\d\xx\d t
\\&\qquad
{\buildrel{\scriptsize\eq{Euler-small-therm-ED+0dis}}\over{=}}
\int_0^T\!\!\!\int_\varOmega\!\pdt{\pp_\tau\!}\Cdot\ovT
+\frac{|\ovT|^2\!}2\,{\rm div}\,\opT\,\d\xx\d t
=\int_0^T\!\!\!\int_\varOmega\!\pdt{\pp_\tau\!}\Cdot\ovT
+\ovT\Cdot{\rm div}\big(\opT{\otimes}\ovT\big)\,\d\xx\d t\,,
\label{Euler-large-strong-calc}
\end{align}
where, for the last equality, we have used \eq{calculus-convective}.
Now we want to pass to the limit in \eq{Euler-large-strong+} or, more
precisely, to estimate the limit superior from above. For this, we again
use the kinetic-energy convexity, which causes the weak lower semicontinuity
of $(\varrho,\pp)\mapsto\int_\varOmega|\pp|^2/\varrho\,\d\xx$ as a convex functional
$\{\rho{\in}L^1(\varOmega);\,\rho\ge0\}\times L^1(\varOmega;\R^3)\to[0,+\infty]$.
Here, we rely also on that $|\pp_\tau(T)|^2/\varrho_\tau(T)$ is bounded in
$L^1(\varOmega)$ due to the former estimate in \eq{Euler-small-est1} and on
that $\varrho_\tau(T)\to\varrho(T)$ even strongly in $C(\barOmega)$ due
to \eq{Euler-small-converge-strong-rho+}, and on that $\pp_\tau(T)$
converges weakly* in $C(\barOmega)^*$, i.e.\ as measures on $\barOmega$ due
to \eq{Euler-small-est3+++} to its limit which is $\pp(T)$ because
simultaneously $\pp_\tau(T)\to\pp(T)$ weakly in $W^{2,p}(\varOmega;\R^3)^*$
due to \eq{Euler-small-converge-p}.
For the term $(\opT{\otimes}\ovT)\Colon\strain(\wt\vv)$,
we use simply \eq{Euler-small-converge-strong-pxv}.
All of this allows us to estimate of the limit superior of the right-hand side in
\eq{Euler-large-strong+} from above, so that:
\begin{align}
\nonumber
&\!\!\!\limsup_{\tau\to0}\Big(\big(\inf_{|E|=1}{\bbD}E\Colon E\big)
\|\strain(\ovT{-}\wt\vv)\|_{L^2(I\times\varOmega;\R^{3\times3})}^2
+\mu c_p\|\nabla^2(\ovT{-}\wt\vv)\|_{L^2(I\times\varOmega;\R^{3\times3\times3})}^p\Big)
\\[-.1em]&\nonumber\ \
\le\int_0^T\!\!\bigg\langle\pdt{\pp\NOINDEX},\wt\vv\bigg\rangle
+\!\int_\varOmega\bigg(\varrho\NOINDEX\GRAVITY\NOINDEX\Cdot(\vv\NOINDEX{-}\wt\vv)
-(\pp\NOINDEX{\otimes}\vv\NOINDEX)\Colon\strain(\wt\vv)
  -\big(\mathscr{T}(\Ee,\theta){+}\bbD\strain(\wt\vv)\big)
  \Colon\strain(\vv\NOINDEX{-}\wt\vv)
\\[-.7em]&\nonumber\hspace{9em}
-\HYPER|\nabla^2\wt\vv|^{p-2}\nabla^2\wt\vv\Vdots\nabla^2(\vv\NOINDEX{-}\wt\vv)
\bigg)\,\d\xx\d t+
\int_\varOmega\frac{|\pp_0|^2}{2\varrho_0}-\frac{|\pp\NOINDEX(T)|^2}{2\varrho\NOINDEX(T)}\,\d\xx
\\[-.5em]&\ \ =\nonumber
 \int_0^T\!\!\bigg\langle\pdt{\pp\NOINDEX},\wt\vv-\vv\NOINDEX\bigg\rangle
 +\int_\varOmega\!\bigg(\!\varrho\NOINDEX\GRAVITY\NOINDEX\Cdot(\vv\NOINDEX{-}\wt\vv)
-\Big(\mathscr{T}(\Ee,\theta){-}\pp\NOINDEX{\otimes}\vv\NOINDEX\Big)
\Colon\strain(\vv\NOINDEX{-}\wt\vv)
\\[-.7em]&\hspace{11.5em}
-\bbD\strain(\wt\vv)\Colon\strain(\vv\NOINDEX{-}\wt\vv)
-\HYPER|\nabla^2\wt\vv|^{p-2}\nabla^2\wt\vv\Vdots\nabla^2(\vv\NOINDEX{-}\wt\vv)
\bigg)\,\d\xx\d t\,,
\label{Euler-large-strong+++}\end{align}
where $\langle\cdot,\cdot\rangle$ denotes the duality between
$W^{2,p}(\varOmega;\R^3)^*$ and $W^{2,p}(\varOmega;\R^3)$ and where, for
the last equality, we used the calculus like \eq{Euler-large-strong-calc}
but for the continuous-in-time limit which holds as an equality.
Choosing $\wt\vv=\vv\NOINDEX$ and reminding
also \eq{Euler-small-converge-strong-vv}, we obtained
\begin{align}\label{Euler-small-converge-strong-v}
\ovT\to\vv\ \ \text{ strongly in }\ L^p(I;W^{2,p}(\varOmega;\R^3))\,.
\end{align}
Thus, we can easily make the limit passage in \eq{Euler-small-therm-ED+1dis}
and thus prove that $\vv\NOINDEX$ satisfies the momentum equation
\eq{Euler-small-therm-ED-anal1} in the weak sense \eq{def-thermo-Ch4-momentum}.

 The limit passage in \eq{Euler-small-therm-ED+2dis} is even simpler.
Specifically, due to \eq{Euler-small-converge-bar-E} and
\eq{Euler-small-converge-strong-v}, we have
\begin{align}
\label{Euler-small-converge-strong+B}
&\bm B_\text{\sc zj}^{}(\ovT,\oEeT)\to\bm B_\text{\sc zj}^{}(\vv,\Ee)
\hspace{1em}\text{weakly in }\ L^p(I;L^s(\varOmega;\Rsym))\,.&&
\end{align}
The other terms in \eq{Euler-small-therm-ED+2dis} can be handled by 
\eq{Euler-small-converge-strong-E} with \eq{Euler-small-converge-strong-p-}.

Finally, for the limit passage in \eq{Euler-small-therm-ED+3dis}, we 
exploit the strong convergences in \eq{Euler-thermo-Ch4-conv-kappa-T-R}.

\medskip{\it Step 7: Energy balances}.
It is now important that the tests and then all the subsequent calculations
leading to the energy balance \eq{Euler-small-energy-balance-stress}
integrated over a current time interval $[0,t]$ are analytically legitimate.

First, note that, by \eq{Euler-small-est++}, we have that
$\pdt{}\varrho\in L^p(I;L^r(\varOmega))$ and
$\pdt{}\pp\in L^{p'}(I;W^{2,p}(\varOmega;\R^3)^*)$. By
\eq{Euler-small-therm-ED-est3}, we also have that
$\pdt{}\Ee\in L^p(I;L^\EXS(\varOmega;\Rsym))$. Then, for the calculus
\eq{Euler-small-divT.v++} integrated over in time $[0,t]$, we rely on that 
$\TT_0=\varphi'(\Ee)+\varphi(\Ee)\bbI\in L^\infty(I;L^3(\varOmega;\Rsym))$
is certainly in duality with
$\strain(\vv)\in L_{\rm w}^p(I;L^\infty(\varOmega;\Rsym))$ and
$\varphi'(\Ee)\in L^\infty(I;L^6(\varOmega;\Rsym))$ is in duality both with
$\ZJ\Ee\in L^p(I;L^2(\varOmega;\Rsym))$ as well as  with
${\rm dev}\,\TT_0\in L^\infty(I;L^6(\varOmega;\Rsym))$, too. Thus the calculus
\eq{Euler-small-divT.v++} is indeed legitimate when integrated over time.

Furthermore, $\pdt{}\pp\in L^{p'}(I;W^{2,p}(\varOmega;\R^3)^*)$ and
${\rm div}(\varrho\vv{\otimes}\vv)\in L^{1+p/2}(I{\times}\varOmega;\Rsym)$
as well as ${\rm div}\DD\in L^{p'}(I;W^{2,p}(\varOmega;\R^3)^*)$ are in duality
with $\vv\in L^p(I;W^{2,p}(\varOmega;\R^3))$, as used when testing the momentum
equation by $\vv$, in particular in \eq{calculus-convective}. The calculus
\eq{calculus-for-kinetic} also relies on that both $\pdt{}\varrho$ and
${\rm div}(\varrho\vv)=\vv{\cdot}\nabla\varrho+\varrho\,{\rm div}\vv$ live in
$L^p(I;L^r(\varOmega))$ and are thus certainly in duality with $|\vv|^2\in
L_{\rm w}^{p/2}(I;L^\infty(\varOmega))$ for $p\ge3$.

The right-hand side of \eq{Euler-small-therm-ED-anal3} is
in $L^1(I{\times}\varOmega)$, so that \eq{Euler-small-therm-ED-anal3}
bears legitimately integration over $[0,t]{\times}\varOmega$. Adding this
to \eq{Euler-small-energy-balance-stress} integrating over $[0,t]$ then yields
the total-energy balance \eq{Euler-engr-finite} integrated over $[0,t]$.
Thus the point (iv) is proven.  
\end{proof}

\begin{remark}[{\sl Heat bulk sources}]\label{rem-bulk-heat-sources}\upshape
The above analysis could easily be extended if the heat equation
\eq{Euler-small-therm-ED-anal3} had a non-negative right-hand side in
$L^1(I{\times}\varOmega)$. Such sources arise in engineering from some
exo- or endo-thermic chemical processes or from electric eddy
currents etc.\ or in geophysics in the planetary mantle from radiogenic heating
or absorption radiation from the Sun in the atmosphere.
\end{remark}

\begin{remark}[{\sl The restrictions on $\alpha$ and $\beta$ in
\eq{Euler-small-ass-alpha-beta}}]\label{rem-alpha-beta}\upshape
Depending on $\lambda>0$, the bounds \eq{Euler-small-ass-alpha-beta} restrict
possible exponents $\alpha$ and $\beta$. For $\lambda>2$, no pair is consistent
with \eq{Euler-small-ass-alpha-beta}. Depending on $0<\lambda\le2$, the
admissible pairs $(\alpha,\beta)$ lie within  the polyhedrons as in 
Figure~\ref{fig-restriction-alpha-beta}.
The usually considered situation that $\kappa$ bounded and the test for
$\lambda\to0$ leads here to the restriction $0\le\alpha<1/2$, as e.g.\ in 
\cite{KruRou19MMCM,Roub25TEFS}. 
Allowing for a growth of $\kappa=\kappa(\theta)$ as $\sim\theta^\beta$ with
$\beta>0$ and using a general test by some fixed positive $\lambda$ opens
possibilities for $1/2\le\alpha<2$.
\begin{figure}[ht]
\begin{center}
\psfrag{l=3/2}{\small$\lambda=3/2$}
\psfrag{l=1}{\small$\lambda=1$}
\psfrag{l=1/2}{\small$\lambda=1/2$}
\psfrag{l -> 0}{\small$\lambda\searrow0$}
\hspace*{-2em}\includegraphics[height=0.25\textwidth,width=0.65\textwidth]{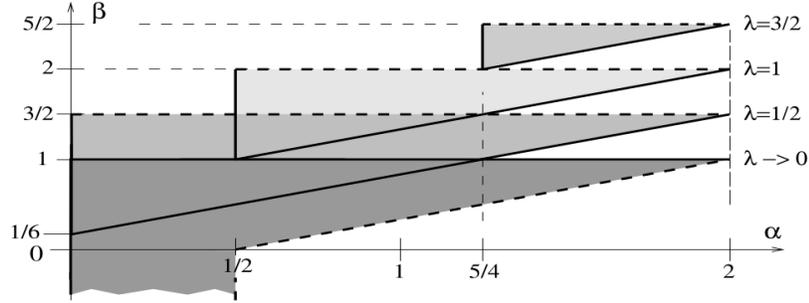}
\end{center}
\vspace*{-.6em}
\caption{\sl The $(\alpha,\beta)$-pairs complying with the
restrictions \eq{Euler-small-ass-alpha-beta} for four values of $0{<}\lambda{<}2$.}
\label{fig-restriction-alpha-beta}
\end{figure}
\end{remark}

\begin{remark}[{\sl The general  $\Ee$-dependent  heat conductivity}]
\label{rem-Euler-thermo-kappa}\upshape
When $\kappa=\kappa(\Ee,\theta)$, we can generalize $\intkappa$ in
\eq{def-thermo-Ch4-momentum3} as $\intkappa(\Ee,\theta):=
\int_0^\theta\kappa(\Ee,\vartheta)\,\d\vartheta$, so
that $\kappa(\Ee,\theta)\nabla\theta=\nabla\intkappa(\Ee,\theta)-
\intkappa'_{\Ee}(\Ee,\theta)\nabla\Ee$. Thus, the integral identity
\eq{def-thermo-Ch4-momentum3} augments by the term 
$\nabla\Ee\Vdots(\intkappa'_{\Ee}(\Ee,\theta){\otimes}\nabla\wt\theta)$
which, however, is not integrable if the Maxwellian viscosity $\zeta_{\rm p}$
is temperature dependent, as seen in the calculations in Step~4 above.
Therefore, instead of $\intkappa(\theta)\Delta\wt\theta$ in
\eq{def-thermo-Ch4-momentum3}, one should use 
$-\kappa(\Ee,\theta)\nabla\theta\Cdot\nabla\wt\theta$, for which we need the heat
flux to be integrable. Realizing the bound
$\kappa(\oEeT,\theta)=\mathscr{O}(\theta^\beta)$, we obtain the bound for
the (negative) heat flux $\kappa(\oEeT,\otT)\nabla\otT$, namely
\begin{align}\label{Euler-entropy-heat-flux}
\|\kappa(\oEeT,\otT)\nabla\otT
\|_{L^{(4+\alpha)\,\EXP/(4+\alpha+3\beta)}(I{\times}\varOmega;\R^3)}^{}\le C\,.
\end{align}
To legitimate the weak formulation of the heat equation and to facilitate the
limit passage, the exponent in \eq{Euler-entropy-heat-flux} must be greater than
(or equal to) 1. Recalling \eq{Euler-entropy-nabla-theta-mu}, this holds for
$\alpha\ge3\lambda-1$. This is a slightly stronger restriction on the exponent
$\alpha$ than $\alpha\ge\frac32\lambda-1$, which is needed for
\eq{Euler-hat-kappa-est}. For the convergence with $\tau\to0$, in view of
\eq{Euler-entropy-heat-flux}, we also have
\begin{align}&\label{Euler-thermo-Ch4-conv2+}
\kappa(\oEeT,\otT)\nabla\otT\!\to\kappa(\Ee,\theta)\nabla\theta\
\text{ weakly in }\ L^{(4+\alpha)\,\EXP/(4+\alpha+3\beta)}(I{\times}\varOmega;\R^3).\!
\end{align}
\end{remark}

\begin{remark}[{\sl Rate-dependent plasticity}]\label{rem-plasticity}\upshape
The inelastic-strain rate $\mathscr{R}:(E,\theta)$
involving the conjugate dissipation potential $\zeta_{\rm p}^*(\theta,\cdot)$,
considered here as continuously differentiable, can also cover the rate-dependent
plasticity. However, since plasticity is an activated process, it typically
involves a potential that is non-smooth at zero plastification rate. Nevertheless,
some applications (particularly in geophysics, where the activated plastification
is combined with so-called aseismic slip) combine non-smooth plastic-type and 
smooth creep-type (say, quadratic)
potentials ``in series''; such a serial arrangement is called an {\it extended
Maxwell model}. This means that the original dissipation potential
$\zeta_{\rm p}(\theta,\cdot)$ is a so-called infimal convolution of the two
aforementioned  potential and its conjugate $\zeta_{\rm p}^*(\theta,\cdot)$
is the sum of the conjugates of those two original potentials. To obtain
a smooth $\zeta_{\rm p}^*(\theta,\cdot)$, it should still be combined with 
the Stokes viscosity in parallel. This can
advantageously and directly be used in \eq{R-abbreviation}, avoiding
an explicit (and often nontrivial) evaluation of the infimal convolution.
 Cf.\ \cite{Roub25FNVR}. 
\end{remark}

\begin{example}[{\sl Creep in thermally expanding materials}]
\upshape\label{exa-temperature-creep}
The linear creep with the Maxwell  viscosity  modulus $M=M(\theta)>0$ is
governed by the quadratic  dissipation  functional
$\zeta_{\rm p}(\theta,\cdot)=\frac12M(\theta)|\cdot|^2$. Thus
$\mathscr{R}(\bm{E},\theta)=
M^{-1}(\theta){\rm dev}\,\psi'_{\bm{E}}(\bm{E},\theta)
=M^{-1}(\theta){\rm dev}\,\varphi'(\bm{E})+\theta M^{-1}(\theta){\rm dev}\,\phi'(\bm{E})$,
where the ansatz \eq{ansatz-for-entropy-base-estimate} has been taken into
account. When $\psi(\cdot,\theta)$ is convex, the positive semi-definiteness
of $\mathscr{R}'_{\bm{E}}$ holds. For the specific case that complies with
the ansatz \eq{ansatz-for-entropy-base-estimate}, we can consider
\begin{align}\nonumber
\psi(\bm{E},\theta)&=
\frac12K|{\rm tr}\,\bm{E}-\alpha_\text{\sc v}^{}\theta|^2
+G|{\rm dev}\,\bm{E}|^2
-\frac{c_{\rm v}}{\alpha(1{+}\alpha)}\theta^{1+\alpha}
-\tfrac12 K\alpha_\text{\sc v}^2\theta^2
\\&=\frac12K({\rm tr}\bm{E})^2+G|{\rm dev}\,\bm{E}|^2
-\alpha_\text{\sc v}^{}\theta K{\rm tr}\,\bm{E}
-\frac{c_{\rm v}}{\alpha(1{+}\alpha)}\theta^{1+\alpha}
\end{align}
with $K$ and $G$ as the bulk and the shear elastic moduli, respectively, and
with $\alpha_\text{\sc v}^{}$ denoting the {\it thermal volume expansibility}.
Thus, we obtain
\begin{align}
\nonumber
&\ENG(\bm{E},\theta)=\psi(\bm{E},\theta)
-\theta\psi_{\theta}'(\bm{E},\theta)
=\tfrac12K({\rm tr}\,E)^2+G|{\rm dev}\,\bm{E}|^2
+\tfrac{c_{\rm v}}{1{+}\alpha}\theta^{1+\alpha},
\\[-.1em]&\nonumber
\ENG'_{\!E}(E,\theta)=\psi'_{\bm{E}}(\bm{E},\theta)-\theta\psi''_{\bm{E}\theta}(\bm{E},\theta)
=K{\rm tr}\,\bm{E}+2G{\rm dev}\,\bm{E}\,,
\\&\nonumber
\ENG_\theta'(\bm{E},\theta)
=-\theta\psi''_{\theta\theta}(\bm{E},\theta)=c_{\rm v}\theta^\alpha\,,\ \
\text{ and}\ \ \eta(\bm{E},\theta)=-\psi_{\theta}'(\bm{E},\theta)
=\tfrac{c_{\rm v}}{\alpha}\theta^\alpha
+\alpha_\text{\sc v}^{}K\,{\rm tr}\,\bm{E}\,,
\\[-.1em]&\nonumber
\mathscr{T}(\bm{E},\theta)=\psi'_{\!\bm{E}}(\bm{E},\theta)
+\psi(E,\theta)\bbI
=\big(K({\rm tr}\,\bm{E}{-}\alpha_\text{\sc v}^{}\theta)
+\psi(\bm{E},\theta)\big)\,\bbI+2G{\rm dev}\,\bm{E}\,,\ \ \text{ and}
\\[-.2em]&\nonumber
\mathscr{R}(\bm{E},\theta)
=\big[{\zeta_{\rm p}(\theta,\cdot)^*}\big]'\big({\rm dev}\mathscr{T}(\bm{E},\theta)\big)
=\tfrac{2G}{M(\theta)}{\rm dev}\,\bm{E}\,.
\end{align}
This special case complies with (\ref{Euler-thermo-small-ass}b--d).
Also, the convexity of the function $1/[\ENT^{-1}](\cdot)^\lambda$,
used for \eq{calculus-to-entropy-disc-}, is satisfied for any
$0\le\lambda$. To see this, note that, from
$\Ent=\ENT(\theta)=\frac{c_{\rm v}}{1{+}\alpha}\theta^{1+\alpha}$, we can see 
$\theta=\ENT^{-1}(\Ent)=(\frac{1{+}\alpha}{c_{\rm v}}\Ent)^{1/(1+\alpha)}$ so that 
$1/\theta^\lambda
=1/[\ENT^{-1}](\Ent)^\lambda=(\frac{1{+}\alpha}{c_{\rm v}}\Ent)^{-\lambda/(1+\alpha)}$.
Also, note that $(\Ent^{-\lambda/(1+\alpha)})''=
-\frac\lambda{1+\alpha}(\Ent^{-1-\lambda/(1+\alpha)})'
=\frac\lambda{1+\alpha}(1+\frac\lambda{1+\alpha})\Ent^{-2-\lambda/(1+\alpha)}>0$
so that $1/[\ENT^{-1}](\cdot)^\lambda$ is convex on $\R^+$.
The function $\wt\eta_\lambda^{}=\wt\eta_\lambda^{}(\Ent)$ from 
\eq{calculus-to-entropy-disc-} is then
$\wt\eta_\lambda^{}(\Ent)=
\frac{1{+}\alpha}{1+\alpha-\lambda}(\frac{1{+}\alpha}{c_{\rm v}}\Ent)^{1-\lambda/(1+\alpha)}$.
\end{example}

{\small
\section*{{\large Acknowledgments}}

\vspace*{-1em}

This research has been partially supported also from the CSF (Czech Science
Foundation) project 23-06220S and by the institutional support RVO: 61388998
(\v CR).

\medskip
}

\end{sloppypar}
\end{document}